\documentclass[multi]{batemantight_vol2}
\usepackage{makeidx}
\usepackage{david}
\usepackage{latexsym}
\usepackage{amsfonts}
\usepackage{amsmath}
\usepackage{amsthm}
\usepackage{amssymb}
\usepackage[numbers,sectionbib]{natbib}
\usepackage{chapterbib}
\usepackage{txfonts}
\usepackage{url}
\usepackage{graphicx}
\usepackage{eucal}
\usepackage{enumitem}
\setlist{nosep}

\newtheorem{theorem}{Theorem}[section]
\newtheorem{lemma}[theorem]{Lemma}
\newtheorem{proposition}[theorem]{Proposition}
\newtheorem{corollary}[theorem]{Corollary}
\newtheorem{definition}[theorem]{Definition}
\newtheorem{example}[theorem]{Example}
\newtheorem{remark}[theorem]{Remark}
\newtheorem{conjecture}[theorem]{Conjecture}
\newtheorem*{acknowledgements*}{Acknowledgements}

\newcommand\theo[2]{\begin{theorem}\label{#1}{\sl\  #2}\end{theorem}}

\newcommand\prop[2]
   {\begin{proposition}\label{#1}{\sl\ #2}\end{proposition}}
\newcommand\corr[2]
    {\begin{corollary}\label{#1}{\sl\ #2}\end{corollary}}

\newcommand\acknow[1]
    {\begin{acknowledgements*}{\rm\ #1}\end{acknowledgements*}}

\newcommand{\hyp}[5]{\,\mbox{}_{#1}F_{#2}\left(
  \genfrac{}{}{0pt}{}{#3}{#4};#5\right)}

\def\KZ/{{\slshape KZ\/}}
\def\qKZ/{{\slshape qKZ\/}}
\def\KZB/{{\slshape KZB\/}}
\def\qKZB/{{\slshape qKZB\/}}

\allowdisplaybreaks

\newcommand\gobbleone[1]{}

\let\oldindex\index

\renewcommand{\index}[1]{\def\exptoindex{#1}\expandafter\oldindex\expandafter{\exptoindex}}

\begin{document}
\author[Yuan Xu]{Yuan Xu}
\chapter[Ch.~2, Orthogonal polynomials of several variables]{Orthogonal
polynomials of several variables}\footnote{This is a
preliminary version of
Chapter 2 in the book \emph{Encyclopedia of special functions:
The Askey--Bateman project, Vol.~2: Multivariable special functions},
T.~H. Koornwinder and J.~V. Stokman (eds.), Cambridge University Press, 2021.}

\contributor{Yuan Xu
\affiliation{Department of Mathematics, University of Oregon, Eugene,
OR 97403-1222, USA}}

\def\tr{{\triangle}}
\newcommand{\fs}{\triangle}
\newcommand{\bs}{\triangledown}
\def\ll{\lesssim}
\def\sub{\substack}
\def\sph{\mathbb{S}^{d-1}}
\def\f{\frac}

 \def\a{{\alpha}} 
 \def\b{{\beta}}
 \def\g{{\gamma}}
 \def\k{{\kappa}}
 \def\t{{\theta}}
 \def\l{{\lambda}}
 \def\d{{\delta}}
 \def\o{{\omega}}
 \def\s{{\sigma}}
 \def\la{{\langle}}
 \def\ra{{\rangle}}
 \def\ve{{\varepsilon}}

 \def\ab{{\mathbf a}}
 \def\bb{{\mathbf b}}
 \def\cb{{\mathbf c}}
 \def\eb{{\mathbf e}}
 \def\kb{{\mathbf k}}
 \def\jb{{\mathbf j}}
 \def\rb{{\mathbf r}}
 \def\sb{{\mathbf s}}
 \def\tb{{\mathbf t}}
 \def\ub{{\mathbf u}}
 \def\vb{{\mathbf v}}
 \def\xb{{\mathbf x}}
 \def\yb{{\mathbf y}}
 \def\Ab{{\mathbf A}}
 \def\Bb{{\mathbf B}}
 \def\Db{{\mathbf D}}
 \def\Eb{{\mathbf E}}
 \def\Fb{{\mathbf F}}     
 \def\Gb{{\mathbf G}}
 \def\Kb{{\mathbf K}}
 \def\CA{{\mathcal A}}
 \def\CB{{\mathcal B}}
 \def\CD{{\mathcal D}}
 \def\CF{{\mathcal F}}
 \def\CH{{\mathcal H}}
 \def\CI{{\mathcal I}}
 \def\CJ{{\mathcal J}}
 \def\CL{{\mathcal L}}
 \def\CO{{\mathcal O}}
 \def\CP{{\mathcal P}}     
 \def\CM{{\mathcal M}}     
 \def\CT{{\mathcal T}}     
 \def\CV{{\mathcal V}}
 \def\CW{{\mathcal W}}
 \def\BB{{\mathbb B}}
 \def\CC{{\mathbb C}}
 \def\HH{{\mathbb H}}
  \def\KK{{\mathbb K}}
 \def\NN{{\mathbb N}}
 \def\PP{{\mathbb P}}
 \def\QQ{{\mathbb Q}}
 \def\RR{{\mathbb R}}
  \def\SS{{\mathbb S}}
 \def\ZZ{{\mathbb Z}}
        \def\grad{\operatorname{grad}}
        \def\proj{\operatorname{proj}}
        \def\tan{\operatorname{tan}}
        \def\sspan{\operatorname{span}}
\newcommand{\sC}{{\mathsf {C}}}
\newcommand{\TC}{{\mathsf {TC}}}
\newcommand{\TS}{{\mathsf {TS}}}
\newcommand{\tran}{{\mathsf {tr}}}
 \def\LT{\operatorname{\sf LT}}

\newcommand{\wt}{\widetilde}
\newcommand{\wh}{\widehat}

\newcommand{\IM}{\operatorname{Im}}
\newcommand{\esssup}{\operatorname{ess\,sup}}
\newcommand{\meas}{\operatorname{meas}}
\newcommand{\rank}{\operatorname{rank}}

\newcommand{\half}{\frac12}
\newcommand{\thalf}{\tfrac12}
 
\section{Introduction} \label{Introduction}
Polynomials of $d$ variables are indexed by the set $\NN_0^d$ of
multi-indices, where $\NN_0:=\{0,1,2,\ldots\}$.
The standard multi-index notations will be used throughout this chapter. For 
$\alpha = (\alpha_1, \ldots, \alpha_d)\in \NN_0^d$ and $x = (x_1,\ldots,x_d)$, a
monomial $x^\a$ is defined by $x^\alpha = x_1^{\alpha_1} \ldots x_d^{\alpha_d}$.
The number $|\alpha| = \alpha_1 + \cdots + \alpha_d$ is called the (total) 
degree of $x^\alpha$. The space of homogeneous polynomials of degree $n$
is denoted by 
\begin{equation*} 
\CP_n^d :=
\textrm{span}\left\{x^\alpha\;\big|\; |\alpha| =n,\; \alpha \in \NN_0^d\right\}
\end{equation*} 
and the space of polynomials of degree at most $n$ is denoted by 
$$
\Pi_n^d := 
\textrm{span}\left\{x^\alpha\;\big|\; |\alpha| \le n,\;\alpha \in \NN_0^d\right\}. 
$$ 
Evidently, $\Pi_n^d$ is a direct sum of
$\CP_k^d$ for $k =0,1,\ldots, n$. Furthermore, 
\begin{equation} \label{dimPn}  
 \dim \CP_n^d = \binom{n+d-1}{n} \quad \hbox{and} \quad 
   \dim \Pi_n^d = \binom{n+d}{n}.   
\end{equation}

Let $\la \cdot, \cdot \ra$ be an inner product defined on the space of polynomials of $d$ variables.
A priori it may be indefinite or even degenerate.
Usually it will be given by
$$
      \la f, g \ra_\mu := \int_{\RR^d} f(x)g(x)\,d \mu(x),
$$
where the
{\it orthogonality measure}
\index{orthogonality measure}%
$\mu$ is a positive Borel measure on $\RR^d$ such that the
integral is well defined on polynomials. This inner product will
be nondegenerate, and hence positive definite if
$\mu$ is supported on a set $\Omega$ that has nonempty interior.
This will almost always be the case in this chapter,
apart from some exceptional cases in \S\ref{sect:discreteOP}.

A polynomial $P \in \Pi_n^d$ is said to be an
{\it orthogonal polynomial}
\index{orthogonal polynomials of several variables}%
of degree $n$ with respect to $\la \cdot, \cdot\ra$ if
$P$ is orthogonal to all polynomials of degree $<n$:
\begin{equation} \label{eq:ortho} 
\la P, Q \ra = 0, \qquad \forall Q \in \Pi^d \quad \hbox{with}
\quad \deg Q < \deg P. 
\end{equation}
However, two linearly independent orthogonal
polynomials of degree $n$ are not necessarily orthogonl to each other.
Let $\CV_n^d$ be the space of orthogonal polynomials of degree $n$, that is, 
\begin{equation} \label{CVn}
  \CV_n^d : = \left \{P \in \Pi_n^d\;\big|\; \la P, Q \ra = 0, \quad
  \forall Q \in \Pi_{n-1}^d \right \}.
\end{equation}
If the inner product is nondegenerate then
$\dim \CV_n^d = \dim \CP_n^d = \binom{n+d-1}{n} : = r_n^d$,
and $\Pi_n^d$ is a direct sum of
$\CV_k^d$ for $k =0,1,\ldots, n$. 

Given a nondegenerate inner product, we can assign to the set
$\big\{x^\a\mid \alpha \in \NN_0^d\big\}$ a linear order $\succ$
which is {\it graded}
\index{order!graded}%
(i.e.,  $x^\alpha \succ x^\beta$ if $|\alpha| > |\beta|$),
and apply the Gram--Schmidt orthogonalization process 
to generate a sequence of orthogonal polynomials. In contrast to $d =1$, however, 
there is no obvious natural graded order among monomials when $d > 1$. There are 
instead many well defined orders. One example is given by the
\emph{graded lexicographic order}:
\index{order!graded lexicographic}%
\\[\medskipamount]
$x^\alpha \succ x^\beta$ 
if $|\alpha| > |\beta|$ or if $|\alpha| = |\beta|$ and the first nonzero entry
in the difference $\alpha - \beta$ is positive.
\\[\medskipamount]
In general, different orderings will lead to different orthogonal systems. 
Consequently, orthogonal polynomials of several variables are not unique. Moreover, 
any system of orthogonal polynomials obtained by an ordering of the monomials is 
necessarily unsymmetric in the variables $x_1, \ldots, x_d$. These were recognized
as the essential difficulties in the study of orthogonal polynomials of several variables
in \cite[Ch.~XII]{Er2}, which contains a rather comprehensive account
of the results up to 1950.  

A sequence of 
polynomials $\{P_\alpha \} \in \CV_n^d$ is called
{\it orthogonal}
if $\la P_\alpha, P_\beta \ra = 0$ whenever $\a \ne \b$, and 
{\it orthonormal}
if moreover $\la P_\a, P_\a \ra = 1$ for all $\a$.
The space $\CV_n^d$ can have many 
different bases and a basis does not have to be orthogonal.
\index{orthogonal polynomials of several
variables!orthogonal basis}%
\index{orthogonal polynomials of several
variables!orthonormal basis}%
One way to extend the 
theory of orthogonal polynomials of one variable to several variables is to state the 
results in terms of $\CV^d_0, \CV^d_1, \ldots, \CV^d_n, \ldots$,
rather than in terms 
of a particular basis in each $\CV^d_n$. 

This point of view will be prominent in our next section, which contains a brief account
on the general properties of the orthogonal polynomials of several variables, mostly 
developed in the last two decades. In the later sections of this chapter, we will discuss
in more details specific systems of orthogonal polynomials in two and more variables
that correspond to, or are generalizations of, the classical orthogonal polynomials of one 
variable. Most of these systems are of separated type,
by which we mean that a basis 
of orthogonal polynomials can be expressed as products in some
separation of 
variables in terms of classical orthogonal polynomials
of one variable.

There are many points of contact with other chapters of this volume.
Some of the  orthogonal polynomials will be given in terms of Appell and
Lauricella hypergeometric functions, which are the subject of Chapter 3.
Orthogonal polynomials for weight function invariant under a 
reflection group are addressed in Chapter 7.
Orthogonal polynomials associated with root systems are treated in
Chapter 8. 
$q$-Analogues of such orthogonal polynomials are discussed in Chapter 9. 

\section{General properties of orthogonal polynomials of several variables}\label{sec:General}

The general properties of orthogonal polynomials of several variables
were studied as early as as 1936 by Jackson \cite{Jackson}.
Most earlier studies dealt with two variables, see references 
in \cite[Ch.~XII]{Er2} and \cite{Suetin}. The presentation below follows
the book \cite{DunklXu} by Dunkl and Xu.

\subsection{Moments and orthogonal polynomials} \label{sec:MomentsOP}
Associated with each multi-sequence $s\colon\NN_0^d \mapsto \RR$, 
$s = (s_\a)_{\a \in \NN_0^d}$, we can define a linear functional $\CL$, called 
{\it moment functional},
\index{moment functional}%
by 
$$
    \CL(x^\a) = s_\a, \qquad \a \in \NN_0^d. 
$$
A polynomial $P \in \Pi_n^d$ is called orthogonal with respect to $\CL$ if it is 
orthogonal with respect to the bilinear form $\la f, g \ra = \CL(f g)$,
which however is not necessarily a positive definite or
nondegenerate inner product.

For each $n \in \NN_0$ let $\xb^n$ denote the column vector
$$
  \xb^n : = (x^\alpha)_{|\alpha|=n} = (x^{\alpha_j})_{j=1}^{r_n^d},
$$
where $\alpha_1, \alpha_2,  \ldots, \alpha_{r_n^d}$, $r_n^d =  \dim \CP_n^d$,
is the arrangement of the elements in
$\big\{\alpha \in \NN_0^d\mid |\alpha|=n \big\}$ 
according to the lexicographical order. For $k, j \in \NN_0$ define 
a vector of moments $\sb_k$ and a matrix of moments
$\sb_{k,j}$ by
\begin{equation} \label{eq:sb-def}
 \sb_k : = \CL (\xb^k ) \quad \hbox{and} \quad
\sb_{k,j}= \sb_{\{k\}+\{j\}} := \CL (\xb^k (\xb^j)^\tran). 
\end{equation}
By definition, $\sb_{\{k\}+\{j\}}$ is a matrix of size $r_k^d \times r_j^d$, 
its elements are $s_{\alpha+\beta}$ for $|\alpha|=k$ and $|\beta|=j$.
Finally, for each $n \in \NN_0$, we define a
{\it moment matrix}\index{moment matrix}
by using $\sb_{\{k\}+\{j\}}$ as its building blocks,
\begin{equation}  \label{eq:Mnd-def}
  M_{n,d} : = \left(\sb_{\{k\}+\{j\}}\right)_{k,j=0}^n \quad \hbox{and} \quad
  \Delta_{n,d} : = \det M_{n,d}. 
\end{equation}
The elements of $M_{n,d}$ are $s_{\a+\b}$ for $|\alpha|\le n$ and $|\beta|\le n$. 

\theo{ch2_thk1} 
{\hskip-0.2cm{\rm\cite[Theorem 3.2.6]{DunklXu}}\\
Let $\CL$ be a moment functional. The corresponing inner product is
nondegenerate} if and only if $\Delta_{n,d} \ne 0$ for all $n \in \NN_0$.

From now on in this section we assume the above nondegeneracy
condition. Then orthogonal bases of $\CV_n^d$ esist.
A special, usually not orthogonal, basis
can be expressed in terms of moments
$\CL$ as follows. For $\a \in \NN_0^d$ we denote by
$\sb_{\a,k}$ the column vector 
$\sb_{\a,k}: = \CL (x^\a \xb^k)$; in particular, $\sb_{\a,0} = s_\a$. 

\theo{thm:OPinMoments}
{\hskip-0.2cm{\rm\cite[Theorem 3.2.12]{DunklXu}}\\
For $|\a| = n$, the polynomials
\goodbreak
\begin{equation} \label{eq:MonoOPbyMn}
P_\a^n (x): = \frac{1}{\Delta_{n-1,d}}  \det \left[ \begin{array}{c|c}
  M_{n-1,d} & \begin{matrix} \sb_{\a,0} \\ \sb_{\a,1} \\ \vdots \\ \sb_{\a,n-1} \end{matrix}\\ 
  \hline 1, \xb^\tran \hdots (\xb^{n-1})^\tran & x^\a
      \end{array} \right]
\end{equation}
form a basis for the space $\CV_n^d$ of
orthogonal polynomials of degree $n$ with respect to $\CL$.}

The polynomial $P_\a^n$ can also be characterized as
the unique polynomial in $\CV_n^d$ of the
form
$$
   P_\a^n(x) = x^\a + Q_{n-1}(x), \qquad Q_{n-1} \in \Pi_{n-1}^d. 
$$
It is evident that a polynomial thus characterized exists.
Because of the leading term $x^\a$ the basis $\{P_\a^n\}_{|\a|=n}$
is called a
{\it monic} or
{\it monomial} basis of orthogonal polynomials.
\index{orthogonal polynomials of several variables!monic basis}%
\index{orthogonal polynomials of several variables!monomial
basis|seeonly{monic basis}}%

If $\CL(p^2) > 0$ for all nonzero polynomials $p$, then
the moment functional $\CL$ is
called {\it positive definite}.
\index{moment functional!positive definite}%
In that case the corresponding inner product 
$\la f,g \ra = \CL (f g)$ is also positive definite and
orthonormal bases of polynomials with respect to $\CL$ will exist.

\theo{ch2_thk2}
{\hskip-0.2cm{\rm\cite[Lemma 3.2.8]{DunklXu}}\quad
If $\CL$ is positive definite, then $\Delta_{n,d} > 0$ for all $n \in \NN_0$.}

For a sequence of polynomials $\big\{P_\alpha\mid |\a| =n\big\}$,
we denote by $\PP_n$ 
the polynomial (column) vector 
\begin{equation}  \label{eq:PPdef} 
 \PP_n : = (P_\alpha^n)_{|\alpha|=n} = (P_{\alpha^{(1)}}^n, \ldots, P_{\alpha^{(r_n)}}^n)^\tran ,
\end{equation}
where $\alpha^{(1)}, \ldots, \alpha^{(r_n)}$ is the arrangement of elements in 
$\{\alpha \in \NN_0^d: |\alpha|=n\}$ according to a fixed monomial order. We also regard
$\PP_n$ as a set of polynomials
$\big\{P_\alpha^n\mid |\a| =n\big\}$. Many properties of orthogonal
polynomials of several variables can be expressed more compactly in
terms of $\PP_n$. For example, orthonormality of
$\big\{P_\a\mid  \a \in \NN_0^d\big\}$ with respect to $\la \cdot,\cdot \ra$
can be written as 
$$
      \la \PP_n, \PP_m \ra =  \delta_{m,n} I_{r_n^d},
$$
where $I_k$ denote the identity matrix of size $k \times k$. 

The rows of $M_{n,d}$ are indexed by $\{\a: |\a| \le n\}$. The row indexed by $\a$ is 
\begin{equation*} 
 (\sb_{\a,0}^\tran, \sb_{\alpha,1}^\tran, \ldots, \sb_{\alpha,n}^\tran) = \CL\left(x^{\alpha}, 
    x^{\alpha} \xb^\tran, \ldots, x^{\alpha}(\xb^n)^\tran \right).
\end{equation*}
For $|\a| =n$, let $\wt M_{\a}(x)$ be the matrix obtained from $M_{n,d}$ by replacing the above row of 
index $\a$ by $\left(1, \xb^\tran, \ldots, (\xb^n)^\tran \right).$ Define
$$
   \wt P_\a (x) : = \frac{1}{\Delta_{n,d}} \det \wt M_\a(x), \qquad |\a| = n, \quad \a \in \NN_0^d. 
$$
Let $N_{n,d}$ denote the principal submatrix of the inverse matrix
$M_{n,d}^{-1}$
of size $r_n^d \times r_n^d$ at the lower right corner.
Then $N_{n,d}$ is positive definite. 

\theo{thm:normalOP}  
{\hskip-0.2cm{\rm\cite[Theorem 3.2.13, (3.2.14)]{DunklXu}}\\
Let $\CL$ be a positive definite moment functional. Then  
\begin{equation} \label{eq:normalOP}  
  \PP_n(x) : =  (N_{n,d})^{-\half} \wt \PP_n(x) = 
       G_n \xb^n + \cdots 
\end{equation} 
consists of orthonormal polynomials.
Furthermore, the matrix $G_n$ is positive definite and
\[
G_n=(N_{n,d})^{\half},\qquad
\det G_n = (\Delta_{n-1,d} /\Delta_{n,d})^{1/2}.
\]}

As mentioned before, the space of orthogonal polynomials
$\CV_n^d$ of degree $n$
has many different bases. The orthonormal bases, however,
are unique up to   orthogonal transformations. 

\theo{thm:OrthonormalUnique}
{\hskip-0.2cm{\rm\cite[Theorem 3.2.14]{DunklXu}}\\
Let $\CL$ be positive definite and let $\{Q_{\alpha}^n\}$ be a sequence 
of orthonormal polynomials forming a basis of $\CV_n^d$.
Then there is an orthogonal
matrix $O_n$ such that $\QQ_n = O_n \PP_n$, where $\PP_n$
are the orthonormal
polynomials defined in \eqref{eq:normalOP}.}

Let $\CM = \CM(\RR^d)$ denote the set of nonnegative Borel measures on $\RR^d$
such that $\int_{\RR^d} |x^\a| d\mu(x) < \infty$ for all $\a \in \NN_0^d$. Each $\mu \in
\CM$ defines a positive linear functional 
\begin{equation} \label{eq:LinearFuncational}
     \CL f = \int_{\RR^d} f(x)\,d\mu(x), \qquad f \in \Pi^d.  
\end{equation}
A polynomial $p$ is orthogonal with respect to
the {\it orthogonality measure}
\index{orthogonality measure}%
\index{orthogonal polynomials of several variables!orthogonality
measure}%
$\mu$ if it is orthogonal with respect
to $\CL$ defined in \eqref{eq:LinearFuncational}. Thus, all theorems in this subsection
apply to the orthogonal polynomials with respect to $d \mu$ for $\mu \in \CM$. 

On the other hand, not every positive definite linear functional admits an integral 
representation \eqref{eq:LinearFuncational}.
The moment problem asks when a linear functional, defined via its 
moments, admits an integral representation \eqref{eq:LinearFuncational} for a $\mu \in \CM$
and, if so, when the measure will be determinate. Here a measure
$\mu$ is called
{\it determinate}
\index{orthogonality measure!determinate}%
if no other measure in $\CM$ has all its moments equal to
those of $\mu$. 
It is known that $\CL$ admits an integral representation if, and 
only if, $\CL$ is {\it positive}
\index{moment functional!positive}%
in the sense that $\CL p \ge 0$ for every nonnegative 
polynomial $p$. Evidently a positive linear functional is 
necessarily positive semidefinite (i.e, $\CL p^2 \ge 0$ for every
polynomial $p$), which 
also holds sufficiently for $d=1$. For $d > 1$, however, a positive
definite linear functional 
may not be positive:
there existnonnegative polynomials that cannot be written
as a sum of squared polynomials. For moment problems of several
variables, including various sufficient 
conditions on a measure being determinate, see
\cite{Berg, Fuglede, Schmu} and 
references therein.

\subsection{Three-term relations} \label{sec:ThreeTerm}

Every system of orthogonal polynomials of one variable satisfies a three term relation,
which can also be used to compute orthogonal polynomials recursively. For orthogonal
polynomials of several variables, an analogue of the three-term relation is stated in terms 
of $\CV_n^d$, or rather, in terms of $\PP_n$. 
\goodbreak
\theo{thm:3term}
{\hskip-0.2cm{\rm\cite[Theorem 3.3.1]{DunklXu}}\\
Let $\PP_n$ denote a basis of $\CV_n^d$, $n = 0,1,\ldots$ and
let $\PP_{-1}(x) :=0$. Then there exist unique matrices
$A_{n,i}\colon r_n^d\times r^d_{n+1}$,
$B_{n,i}\colon r_n^d\times r_n^d$, and $C_{n,i}\colon r_n^d \times r^d_{n-1}$
for $n =0,1,2,\ldots$, 
such that
\begin{equation} \label{eq:3term}
   x_i \PP_n(x) = A_{n,i} \PP_{n+1}(x) + B_{n,i} \PP_n(x) + C_{n,i} \PP_{n-1}(x),  \quad 1 \le i \le d.
\end{equation}
In fact, let $H_n: = \CL(\PP_n \PP_n^\tran)$; then the orthogonality shows that
\begin{equation*}
A_{n,i} H_{n+1}  = \CL (x_i   \PP_n  \PP_{n+1}^\tran),\quad 
 B_{n,i} H_n   = \CL (x_i   \PP_n  \PP_n^\tran), \quad
 A_{n,i} H_{n+1}   = H_n C_{n+1,i}^\tran. 
\end{equation*}}
\index{three-term relation}%
\index{orthogonal polynomials of several variables!three-term relation}%

The coefficient matrices $A_{n,i}$ and $C_{n,i}$ in \eqref{eq:3term} have full rank. Indeed, 
\begin{align}
   \rank A_{n,i} & = \rank C_{n+1,i} = r_n^d,  \label{eq:rank1}\\
  \rank A_n  & =  \rank C_{n+1}^\tran = r_{n+1}^d,  \label{eq:rank2}
\end{align}
where $A_n = (A_{n,1}^\tran, \ldots, A_{n,d}^\tran)^\tran$ and  
$C_{n+1}^\tran = (C_{n+1,1}, \ldots, C_{n+1,d})^\tran$ are matrices of 
size $d r_n^d \times r_{n+1}^d$. In the case of orthonormal polynomials, $H_n$ is 
the identity matrix and the three-term relation takes a simpler form.  

\theo{ch2_thk3}
{\hskip-0.2cm{\rm\cite[Theorem 3.3.2]{DunklXu}}\\
If $\PP_n$ is an orthonormal basis of $\CV_n^d$, $n = 0,1,\ldots$, then
\begin{equation} \label{eq:n3term}
x_i \PP_n(x) =A_{n,i}\PP_{n+1} (x) + B_{n,i}\PP_n (x)+ A_{n-1,i}^\tran \PP_{n-1}(x),
 \quad 1 \le i \le d, 
\end{equation} 
where $A_{n,i} = \CL (x_i   \PP_n  \PP_{n+1}^\tran)$, $B_{n,i}  = \CL (x_i   \PP_n  \PP_n^\tran)$ and
$B_{n,i}$ is symmetric.}

As an analogue of the classical Favard's theorem for orthogonal polynomials of one variable,
the three-term relation and the rank conditions characterize orthogonality. 

\theo{thm:nFavard}
{\hskip-0.2cm{\rm\cite[Theorem 3.3.8]{DunklXu}}\\
Let
$\{\PP_n\}_{n=0}^\infty = \big\{P_\alpha^n\mid |\alpha|=n, n\in \NN_0\big\}$,
$\PP_0 =1$, be an arbitrary sequence in $\Pi^d$ such that
$\big\{P_\a^m\mid |\a|=m \le n\big\}$ spans $\Pi_n^d$ for 
each $n$. Then the following statements are equivalent.
\begin{enumerate}[label=\arabic*.]
\item There exists a positive definite linear functional $\CL$ 
that makes $\{\PP_n\}_{n=0}^\infty$ an orthonormal basis for $\Pi^d$.
\item  For $n \ge 0$, $1\le i\le d$, there exist matrices 
$A_{n,i}$ and $B_{n,i}$  such that
\begin{enumerate}[label=\roman*.]
\item
$x_i \PP_n(x) =A_{n,i}\PP_{n+1} (x) + B_{n,i}\PP_n (x)+ A_{n-1,i}^\tran \PP_{n-1}(x),
 \quad 1 \le i \le d,$
\item  $\rank A_{n,i} =  r_n^d, \quad 1 \le i \le d$,\quad
and $\rank A_n  = r_{n+1}^d$.
\end{enumerate}
\end{enumerate}}
\index{Favard type theorem}%
\index{orthogonal polynomials of several variables!Favard type theorem}%

Unlike in one variable, the characterization does not conclude about
the existence of an orthogonality  measure. 

The coefficient matrices of the three-term relation for orthonormal polynomials 
satisfy a set of {\it commutativity conditions} \cite[Theorem 3.4.1]{DunklXu}:
for $1 \le i, j \le d$ and $k \ge 0$, 
\index{three-term relation!commutativity conditions}%
\begin{align} \label{eq:comm}
A_{k,i} A_{k+1,j} &= A_{k,j} A_{k+1,i}, \notag \\  
A_{k,i} B_{k+1,j} + B_{k,i} A_{k,j} &= B_{k,j} A_{k,i} + A_{k,j} B_{k+1,i},
 \\
A_{k-1,i}^\tran A_{k-1,j} + B_{k,i} B_{k,j} + A_{k,i} A_{k,j}^\tran  
& = A_{k-1,j}^\tran  A_{k-1,i} +  B_{k,i} B_{k,j} + A_{k,j} A_{k,i}^\tran, \notag
\end{align}
where $A_{-1,i} : =0$, which is derived from computing, say,  $\CL(x_ix_j \PP_k \PP_{k-1})$, 
by applying the three-term relation in two different ways.  These coefficient matrices also 
define a family of tri-diagonal matrices $J_i$, $1 \le i \le d$: 
\begin{equation} \label{Jmatrix}
J_i = \left[ \begin{matrix} B_{0,i}&A_{0,i}&&\bigcirc\cr
A_{0,i}^\tran &B_{1,i}&A_{1,i}&&\cr
&A_{1,i}^\tran &B_{2,i}&\ddots\cr
\bigcirc&&\ddots&\ddots \end{matrix}
\right],\qquad 1\le i\le d, 
\end{equation}
called {\it block Jacobi matrices}.
\index{block Jacobi matrix}%
\index{three-term relation!block Jacobi matrix}%
The entries of the $J_i$ are matrices that have increasing sizes going down the main diagonal.
The commutativity relations \eqref{eq:comm} is equivalent to the
formal commutativity of $J_i$ \cite[Lemma 3.4.4]{DunklXu}, that is, 
$$
    J_i J_j = J_j J_i, \qquad 1 \le i \ne j \le d. 
$$

These block Jacobi matrices can be viewed as the realization of the multiplication operators 
$X_1,\ldots, X_d$ defined on the space of polynomials by 
$$
      (X_ i f) (x) = x_i f(x),  \qquad  1 \le i \le d.
$$
The operators can be extended to a family of commuting self-adjoint operators on a 
$L^2$ space. This connection to the operator theory allows the use of the spectral theory 
of commuting self-adjoint operators, and helps to answer the question when the inner 
product with respect to which polynomials are orthogonal is defined by a measure. 
It gives, for example, the following theorem, which strengthens Theorem \ref{thm:nFavard}: 

\theo{thm:JFavard}    
{\hskip-0.2cm{\rm\cite[Theorem 3.4.7]{DunklXu}}\\
Let
$\{\PP_n\}_{n=0}^\infty = \big\{P_\alpha^n\mid |\alpha|=n, n\in \NN_0\big\}$,
$\PP_0 =1$, be an arbitrary sequence in $\Pi^d$ such that
$\{P_\a^m\mid |\a|=m \le n\}$ spans $\Pi_n^d$ for 
each $n$. Then the following
statements are equivalent.
\begin{enumerate}[label=\arabic*.]
\item There exists a determinate measure $\mu\in \CM$ with compact support in $\RR^d$ 
such that $\{\PP_n\}_{n=0}^\infty$ is orthonormal with respect to $d\mu$. 
\item  The statement (2) in Theorem \ref{thm:nFavard} holds together with
\begin{equation} \label{eq:Jbound} 
\sup_{k\ge 0} \|A_{k,i}\|_2   < \infty \quad  
  \hbox{and} \quad \sup_{k\ge0} \|B_{k,i}\|_2 < \infty, \quad 1 \le i \le d. 
\end{equation} 
\end{enumerate}}
\index{Favard type theorem}%
\index{orthogonal polynomials of several variables!Favard type theorem}%

If the measure $\mu \in \CM$ is given by
$d\mu(x) = W(x)\,dx$, $W$ being a nonnegative measurable function,
we call $W$ a
{\it weight function},
\index{orthogonal polynomials of several variables!weight function}%
\index{weight function}%
Let $\Omega \subset \RR^d$ be the support set of $W$. 
A function $W$ is called
{\it centrally symmetric},
\index{weight function!centrally symmetric}%
if 
$$
 x \in \Omega \Rightarrow - x \in \Omega, \quad \hbox{and} \quad 
W(x) = W(-x)\hbox{ a.e.\;.} 
$$
For example, the product weight function $\prod_{i=1}^d (1-x_i)^{a_i}(1+x_i)^{b_i}$ on the
cube $[-1,1]^d$ is centrally symmetric if and only if $a_i=b_i$. Furthermore, a linear functional 
$\CL$ is called {\it centrally symmetric}
\index{moment functional!centrally symmetric}%
if 
$$
 \CL(x^{\alpha})=0, \quad \alpha \in \NN^d, \quad
 |\alpha|\hbox{ is odd integer}. 
$$
The two notions are equivalent when $\CL $ is given by
$\CL f = \int f W\,dx$. 

\theo{thm:CenSymm}
{\hskip-0.2cm{\rm\cite[Theorem 3.3.10]{DunklXu}}\\
Let $\CL$ be a positive definite linear functional. Then $\CL$ is centrally symmetric if and 
only if $B_{n,i}=0$ for all $n \in \NN_0$ and $1 \le i \le d$, where $B_{n,i}$ are given 
in \eqref{eq:n3term}. Furthermore, if $\CL$ is centrally symmetric, then an orthogonal 
polynomial of degree $n$ is a sum of monomials of even degrees if $n$ is even, and 
a sum of monomials of odd degrees if $n$ is odd.}

In one variable, the three term relation can be used as a recurrence formula for 
computing orthogonal polynomials of one variable. For several variables, let 
$D_n^\tran = (D_{n,1}^\tran, \ldots, D_{n,d}^\tran)$, where
$D_{n,i}^\tran\colon r^d_{n+1} \times r_n^d$, be a matrix that satisfies 
\begin{equation*} 
   D_n^\tran A_n = \sum_{i=1}^d D_{n,i}^\tran A_{n,i} = I_{r_{n+1}^d}. 
\end{equation*} 
Such a matrix is not unique. 
The three-term relation \eqref{eq:n3term} implies 
\begin{equation} \label{eq:Recur}
 \PP_{n+1} = \sum_{i=1}^d x_i D^\tran_{n,i} \PP_n + E_n  \PP_n + F_n  \PP_{n-1} ,
\end{equation}
\index{recurrence relation}%
\index{orthogonal polynomials of several variables!recurrence relation}%
where $E_n: = - \sum_{i=1}^d D_{n,i}^\tran B_{n,i}$ and 
$F_n: = -\sum_{i=1}^d D_{n,i}^\tran A^\tran_{n-1,i}$. 

Given two sequences of matrices $A_{n,i}$ and $B_{n,i}$, \eqref{eq:Recur} can
be used as a recursive relation to generate a sequence of polynomials. These 
polynomials are orthogonal if the matrices satisfy certain relations:
\theo{thm:Recur}  
{\hskip-0.2cm{\rm\cite[Theorem 3.5.1]{DunklXu}}\\
Let $\{\PP_n\}_{n=0}^{\infty}$ be defined by \eqref{eq:Recur}. Then there is 
a positive definite linear functional $\CL$ that makes $\{\PP_n\}_{n=0}^{\infty}$ an 
orthonormal basis for $\Pi^d$ if and only if $B_{k,i}$ are symmetric, $A_{k,i}$ 
satisfy the rank conditions \eqref{eq:rank1} and \eqref{eq:rank2}, and together they 
satisfy the commutativity conditions \eqref{eq:comm}.}
\noindent
\emph{Further results and references}\quad
The idea of studying orthogonal polynomials of several variables in terms of $\CV_0^d,\CV_1^d,\,\ldots$\; goes back to \cite{KraShe}. The vector notion of three-term relation and 
Favard's theorem were initiated by M. Kowalski in \cite{Kowa1,Kowa2}. The versions in 
this section were developed by Xu in \cite{Xu93a, Xu94a} and subsequent papers. 
Another earlier work is \cite{GeKa}. See \cite[Chapter 3]{DunklXu} for references and further results. 

For further study of three-term relations, see \cite{CSS}.
For an approach based on matrix
factorization, see \cite{AM}. Three-term relations are used for evaluating orthogonal polynomials
of several variables in \cite{BPS}. 

\subsection{Zeros of orthogonal polynomials of several variables}\label{sec:Zeros}

An orthogonal polynomial of degree $n$ in one variable has $n$ distinct real zeros
and the zeros are nodes of a Gaussian quadrature formula. For a polynomial of 
several variables, its set of zero is an algebraic variety, an intrinsically difficult object. The
correct notion for orthogonal polynomials of several variables are the
common  zeros of a family of polynomials, such as~$\PP_n$. 

Let $\PP_n$ be an orthonormal basis of $\CV_n^d$. A point $x \in \RR^d$ is a {\it zero} 
of $\PP_n$ if it is a zero for every elements in $\PP_n$ (or all elements in $\CV_n^d$), 
and it is a {\it simple zero} if at least one partial derivative of $\PP_n$ is not zero at $x$.
\index{orthogonal polynomials of several variables!zeros}%
\index{orthogonal polynomials of several variables!simple zeros}%
Let $A_{n,i}$ and $B_{n,i}$ be matrices in \eqref{eq:n3term}. For each $n \in \NN_0$, 
define the {\it truncated block Jacobi matrices}
\index{block Jacobi matrix!truncated}%
\begin{equation*}
 J_{n,i} :  =  \left[ \begin{matrix} B_{0,i}&A_{0,i}&&&\bigcirc\cr
     A_{0,i}^\tran &B_{1,i}&A_{1,i}&&\cr  &\ddots&\ddots&\ddots&\cr
     &&A_{n-3,i}^\tran &B_{n-2,i}&A_{n-2,i}\cr
     \bigcirc&&&A_{n-2,i}^\tran &B_{n-1,i}\end{matrix} \right],
\qquad 1\le i\le d.
\end{equation*}
These are symmetric matrices of order $N= \dim \Pi_{n-1}^d$. An element $\l \in 
\RR^d$ is called a {\it joint eigenvalue} of $J_{n,1}, \ldots, J_{n,d}$, 
if there is a $\xi \ne 0$, $\xi \in \RR^N$, such that $J_{n,i}\xi = \l_i \xi$ 
for $i= 1, \ldots, d$; the vector $\xi$  is called a {\it joint eigenvector}. 

\theo{thm:JointEigen}
{\hskip-0.2cm{\rm\cite[Theorem 3.7.2]{DunklXu}}\\
A point $\l \in \RR^d$ is a common zero of $\PP_n$ if and only if it is a 
joint eigenvalue of $J_{n,1}, \ldots, J_{n,d}$; moreover, a joint eigenvector of 
$\l$ is $(\PP_0^\tran(\l), \ldots, \PP_{n-1}^\tran (\l))^\tran$.}

Many properties of zeros of $\PP_n$ can be derived from this characterization.

\theo{thm:ZeroSimple}
{\hskip-0.2cm{\rm\cite[Theorem 3.7.1, 3.7.5, Corollary 3.7.3, 3.7.4]{DunklXu}}\\
All zeros of $\PP_n$ are real, distinct and simple. $\PP_n$ has at most $\dim \Pi_{n-1}^d$ 
distinct zeros, and $\PP_n$ has $\dim \Pi_{n-1}^d$ zeros if and only if 
\begin{equation} \label{eq:AiAj=AjAi}
A_{n-1,i} A_{n-1,j}^\tran = A_{n-1,j}A_{n-1,i}^\tran, \quad 1 \le i,j \le d.
\end{equation}}

\theo{thm:ZeroCentral}
{\hskip-0.2cm{\rm\cite[Corollary 3.7.7]{DunklXu}, \cite[Theorem 3.1.1]{Xu94f}}\\
If $L$ is centrally symmetric and $d\ge2$, then $\PP_n$ has less than
$\dim \Pi_{n-1}^d$ common zeros.\\
If moreover $d=2$, then $\PP_n$ has no zero if $n$ is even and has one
zero $(x =0)$ if $n$ is odd.}

As in the case of one variable, zeros of orthogonal polynomials are closely related to 
cubature formulas, which are finite sums that approximate integrals.
A {\it cubature formula} $\CI_{\!n} (f)$ is said to have degree $2n-1$ if
\index{cubature formula}%
\begin{equation} \label{cubature}
  \int_{\RR^d} f(x)\,d\mu(x) = \sum_{k=1}^N \l_k f(x_k) =: \CI_{\!n}(f), \quad \forall f \in \Pi_{2n-1}^d,
\end{equation}
where $\l_k \in \RR$ and $x_k \in \RR^d$, and there is a
polynomial $f^* \in \Pi_{2n}^d$ for which 
the equality does not hold. The number of nodes $N$ of \eqref{cubature} satisfies
a lower bound
\begin{equation} \label{lwbd1}
     N \ge \dim \Pi_{n-1}^d.
\end{equation}
A cubature formula of degree $2n-1$ with $ \dim \Pi_{n-1}^d$ nodes is
called  {\it Gaussian}.
\index{cubature formula!Gaussian}%

\theo{thm:GaussCuba}  
{\hskip-0.2cm{\rm\cite[Theorem 3.8.4]{DunklXu}}\\
Let $\mu \in \CM$ and $\PP_n$ be an orthogonal basis of $\CV_n^d$ with 
respect to $d\mu$. Then the integral $\int f(x)\,d\mu(x)$ admits a
Gaussian cubature 
formula of degree $2n-1$ if and only if $\PP_n$ has $\dim \Pi^d_{n-1}$ 
common zeros.}
\index{orthogonal polynomials of several variables!Gaussian cubature}%

When combined with Theorem \ref{thm:ZeroCentral}, this shows the
following: 

\corr{ch2_thk4}
{If $\mu$ is centrally symmetric, then no Gaussian cubature formulas exist.} 

On the other hand, there are two families of weight functions, discussed in 
\S2.9.1, for which Gaussian cubature formulas do exist. 

The non-existence of Gaussian cubature means that the inequality \eqref{lwbd1} is not sharp. 
There is, in fact, an improved lower bound, which we only give here for
$d=2$.

\theo{ch2_thk5}
{\hskip-0.2cm{\rm\cite[(1.2.10) and  Ch.~5]{Xu94f}}\\
Let $\mu \in \CM(\RR^2)$. Then a cubature formula of degree $2n-1$ for $d\mu$ exists 
only if 
\begin{equation} \label{eq:lwbd2}
   N \ge \dim\Pi_{n-1}^2 + \half \rank \left(  A_{n-1,1} A_{n-1,2}^\tran 
             - A_{n-1,2}A_{n-1,1}^\tran \right).
\end{equation} 
In particular, if $\mu$ is centrally symmetric, then \eqref{eq:lwbd2}
specializes to 
\begin{equation} \label{eq:lwbd3}
N \ge \dim\Pi_{n-1}^2 +  \left\lfloor \frac{n}2 \right \rfloor. 
\end{equation}}

The bounds \eqref{eq:lwbd2}, \eqref{eq:lwbd3} were first obtained by
 M\"oller \cite{Moller73,Moller76}. See for general $d$
M\"oller  \cite{Moller79} and Xu \cite[Ch.~5]{Xu94f}. The bound presented
there for $\mu$ centrally symmetric
is, for $d>2$, sharper than the analogue of  \eqref{eq:lwbd2}.
 
The condition under which the lower bound in \eqref{eq:lwbd3}
is attained is determined as follows.

\theo{ch2_thk6}
{\hskip-0.2cm{\rm\cite[Theorems 4.1.4 and 5.3.1]{Xu94f}}\\
Let $\mu$ be centrally symmetric. A cubature formula of degree $2n-1$ attains the lower bound \eqref{eq:lwbd3}
if and only if its nodes are common zeros of $n+1 - \left \lfloor \frac{n}{2}\right \rfloor$ 
orthogonal polynomials in~$\CV_n^2$.}
\noindent
\emph{Further results and references}\quad
The lower bound \eqref{lwbd1} is classical, see \cite{Stroud, Myso81}. Theorem 2.2.15 was first proved in \cite{Myso76}. 
The first example of cubature formula 
that attains \eqref{eq:lwbd3} was a degree 5 formula on the square constructed by Radon \cite{Radon}. At the moment, 
the only weight functions for which \eqref{eq:lwbd3} is attained for all $n$ are given by the weight functions 
$W(x,y)= |x-y|^{2 \a+1} |x+y|^{2\b+1}(1-x^2)^{\pm 1/2} (1-y^2)^{\pm1/2}$, where $\a,\b \ge - 1/2$, of which the case 
$\a = \b = \-1/2$ is classical (\cite{MoPa, Xu94f}) and the general case is far more recent \cite{Xu12}.

\subsection{Reproducing kernel and Fourier orthogonal expansion}\label{sec:Kernel-Fourier}

Let $\mu \in \CM$ and assume that the space of polynomials is dense in 
$L^2(d \mu)$. Define the projection operator
$\proj_n\colon L^2(d\mu) \to \CV_n^d$  by 
\begin{equation} \label{eq:projection}
\big(\proj_n f\big)(x) := \int_{\RR^d} f(y) P_n(x,y)\,d\mu(y),
\end{equation}
where $P_n(\cdot,\cdot)$ is the \emph{reproducing kernel}
\index{reproducing kernel}%
\index{orthogonal polynomials of several variables!reproducing kernel}%
of $\CV_n^d$ satisfying
\begin{equation} \label{eq:reprod-def}
\int_{\RR^d} f(y) P_n(x,y)\,d\mu(y) = \begin{cases}
f(x),& f \in \CV_n^d,\\
0,& f\in\CV_m^d,\;m\ne n.
\end{cases}
\end{equation} 
Let $\big\{P_\a^n\big\}_{|\a| =n}$ be an orthonormal basis of $\CV_n^d$. Then $P_n(\cdot,\cdot)$
can be expressed as
\begin{equation} \label{eq:reprod}
  P_n(x,y) = \sum_{|\a| = n} P_\a^n(x) P_\a^n(y) = \PP_n^\tran(x) \PP_n(y).
\end{equation}
The projection operator is independent of a particular basis, and so
is the
reproducing kernel, as also seen by Theorem
\ref{thm:OrthonormalUnique}. The standard Hilbert space argument 
shows that $f \in L^2(d\mu)$ has a Fourier orthogonal expansion
\index{Fourier orthogonal expansion}%
\begin{equation} \label{eq:Foruier}
          f  =  \sum_{n=0}^\infty \proj_n f,
          \qquad \forall f \in  L^2(d\mu). 
\end{equation}
In terms of the orthonormal basis $\{P_\alpha^n\}$, the orthogonal expansion 
reads as
\begin{equation} \label{eq:OPexpansion}
  f = \sum_{n=0}^\infty \sum_{|\alpha|=n}a_\alpha^n(f) P_\alpha^n 
  \quad   \hbox{with} \quad 
 a_\alpha^n (f) : = \int_{\RR^d} f(x) P_\alpha^n(x)\,d\mu(x).
\end{equation}
For studying Fourier expansions, it is often important to have a
closed formula
for the kernel $P_n(\cdot,\cdot)$. Such formulas are often
available for classical orthogonal
polynomials in several variables.

The $n$-th \emph{partial sum} $S_{\!n} f$
\index{partial sum}%
\index{Fourier orthogonal expansion!partial sum}%
of the Fourier orthogonal expansion of $f\in L^2(d\mu)$
is defined by
\begin{equation} \label{eq:Sn(f)}
(S_{\!n} f)(x) : = \sum_{k=0}^n \big(\proj_k f\big)(x) =
\int_{\RR^d} K_n(x,y)f(y)\,d\mu(y),
\end{equation} 
where the kernel $K_n(\cdot, \cdot)$ is defined by 
\begin{equation} \label{eq:ReprodK}
 K_n(x,y) : = \sum_{k=0}^n \sum_{|\alpha|=k} P_\alpha^k(x)
      P_\alpha^k(y) = \sum_{k=0}^n  P_k(x,y).  
\end{equation} 
The kernel $K_n(\cdot,\cdot)$ is the reproducing kernel of
$\Pi_n^d$ in $L^2(d\mu)$. It 
satisfies a \emph{Christoffel--Darboux formula},
\index{Christoffel--Darboux formula}%
\index{orthogonal polynomials of several variables!Christoffel--Darboux
formula}%
deduced from the three-term relation.

\theo{Christoffel-Darboux}
{\hskip-0.2cm{\rm\cite[Theorem 3]{Xu93a}, \cite[Theorem 3.6.3]{DunklXu}}\\
Let $\CL$ be a positive definite linear functional, and let 
$\{\PP_k\}_{k=0}^\infty$ be a sequence of orthonormal polynomials with 
respect to $\CL$. Then, for any integer $n\ge 0$, $x,y \in \RR^d$,
\begin{align} \label{eq:Kn(x,y)}
 K_n(x,y) =  \frac{\bigl[ A_{n,i} \PP_{n+1}(x)\bigr]^\tran \PP_{n} (y) - 
  \PP_{n}^\tran \bigl[ A_{n,i} \PP_{n+1} (y)\bigr]}
  {x_i - y_i}, \quad 1 \le i \le d,
\end{align}
where $x = (x_1,\ldots, x_d)$ and $y = (y_1, \ldots,y_d)$.}

The right-hand side of \eqref{eq:Kn(x,y)}
can also be stated in 
terms of orthogonal, instead of orthonormal, polynomials.
A related function is
the \emph{Christoffel function}
\index{Christoffel function}%
\index{orthogonal polynomials of several variables!Christoffel
function}%
defined by  
\begin{equation}\label{eq:Christoffel}
   \Lambda_n(x) : = \left[K_n(x, x)\right]^{-1}. 
\end{equation} 

\theo{thm:Christoffel}
{\hskip-0.2cm{\rm\cite[Theorem 3.6.6]{DunklXu}}\quad
Let $\mu \in \CM$. Then for any $x \in \RR^d$, 
$$
\Lambda_n(x) = \min_{P(x) =1,\;P \in \Pi_n^d}\,\int_{\RR^d} P^2(y)
\,d\mu(y).$$}

\section{Orthogonal polynomials of two variables} \label{sec:OP-2d}
\index{orthogonal polynomials of two variables}%
Almost all that can be stated about the general properties of orthogonal polynomials 
of two variables also holds for orthogonal polynomials of more than two variables. 
This section contains results on various special systems of orthogonal polynomials of 
two variables and their properties.  

A basis of $\CV_n^2$ in two variables is often indexed by $\alpha = (k,n-k)$,
or by a single 
index~$k$, as $\{P_k^n\}_{k=0}^n$. Many examples below will be given in term of 
classical orthogonal polynomials of one variable, which are listed,
together with their associated weight function, orthogonality
interval and parameter constraints,
in Table \ref{table-classical}. The 
normalization given in the last column (here usually the value
attained at an endpoint of the orthogonality interval) makes the
definition precise.
The notation $(a)_n: = a(a+1) \ldots (a+n-1)$, shifted factorial or Pochhammer symbol,
will be used throughout the rest of the chapter.  
 
\begin{table}[h] 
\caption{Classical orthogonal polynomials of one variable}
\label{table-classical} 
\addtolength\tabcolsep{2pt}
\begin{tabular}{@{}c@{\hspace{25pt}}ccccl@{}} 
\hline \hline 
Name & notation & weight & interval & constraint & normalization \\ 
\hline 
Hermite & $H_n(t)$ & $e^{-t^2}$ & $(-\infty,\infty)$ & & $H_n(t)= 2^n t^n + \ldots$\\[3pt] 
Laguerre & $L_n^\a(t)$ & $t^\a e^{-t}$ & $[0,\infty)$ & $\a>-1$ & $L_n^\a(0) = \frac{(\a+1)_n}{n!}$ \\[3pt] 
Chebyshev 1st & $T_n(t)$ &  $(1-t^2)^{-1/2}$ & $[-1,1]$  & & $T_n(1)=1$ \\[3pt] 
Chebyshev 2nd & $U_n(t)$ &  $(1-t^2)^{1/2}$ & $[-1,1]$  & & $U_n(1) = n+1$\\[3pt] 
Gegenbauer & $C_n^{\l}(t)$ & $ (1-t^2)^{\l -1/2}$ & $[-1,1]$  & $\l>-1/2$ & $C_n^\l(1) = \frac{(2 \l)_n}{n!}$\\[3pt] 
Jacobi & $P_n^{(\a,\,\b)}(t)$ & $ (1-t)^\a (1+t)^\b$ & $[-1,1]$  & $\a,\b>-1$ & $P_n^{(\a,\,\b)}(1) = \frac{(\a+1)_n}{n!}$\\
\hline \hline 
\end{tabular} 
\index{orthogonal polynomials of one variable!classical}%
\index{orthogonal polynomials of one variable!classical!Hermite}%
\index{Hermite polynomials}%
\index{orthogonal polynomials of one variable!classical!Laguerre}%
\index{Laguerre polynomials}%
\index{orthogonal polynomials of one variable!classical!Chebyshev}%
\index{Chebyshev polynomials}%
\index{orthogonal polynomials of one variable!classical!Gegenbauer}%
\index{Gegenbauer polynomials}%
\index{orthogonal polynomials of one variable!classical!Jacobi}%
\index{Jacobi polynomials}%
\end{table} 

The Gegenbauer case $\l=0$ can be obtained for $n>0$
by the renormalization
\begin{equation} \label{limGegCheb}
\lim_{\l \to 0} \l^{-1} C_n^\l(x) = 2n^{-1}T_n(x),\quad n>0.
\end{equation}

\subsection{Product orthogonal polynomials} \label{sec:prodOP-2d}
\index{orthogonal polynomials of two variables!product orthogonal
polynomials}%
For the weight function $W(x,y) = w_1(x) w_2(y)$, where $w_1$ and $w_2$ are two weight 
functions of one variable, an orthogonal basis of $\CV_n^d$ is given by 
$$
    P_k^n (x,y) := p_{k}(x) q_{n-k}(y), \qquad 0 \le k \le n,  
$$
where $\{p_k\}$ and $\{q_k\}$ are sequences of orthogonal polynomials with 
respect to $w_1$ and $w_2$, respectively.
If $\{p_k\}$ and $\{q_k\}$ are orthonormal,
then so is $\{P_k^n\}$.
\vskip\medskipamount
\begin{enumerate}[label=\arabic*.]
\item {\it Product Hermite polynomials.}
\index{orthogonal polynomials of two variables!product Hermite
polynomials}%
For weight function
$W(x,y) = e^{-x^2- y^2}$, a possible orthogonal 
basis is given by $P_k^n(x,y) = H_k(x) H_{n-k}(y)$, $0 \le k \le n$. This
satisfies the differential equation
\begin{equation} \label{diff-Hermite-2d}
  \thalf (v_{xx} + v_{yy}) - (x v_x + y v_y)  = - n v.
\end{equation} 
\item {\it Product Laguerre polynomials.}
\index{orthogonal polynomials of two variables!product Laguerre
polynomials}%
For weight function $W(x,y) = x^\a y^\b e^{-x- y}$, a possible orthogonal 
basis is given by $P_k^n(x,y) = L^\a_k(x) L^\a_{n-k}(y)$, $0 \le k \le n$.
This satisfies the differential equation
\begin{equation} \label{diff-Laguerre-2d}
       x v_{xx} + y v_{yy} + (1+\a - x)  v_x + (1+\b-y) v_y  = - n v.
\end{equation} 
\item {\it Product Hermite--Laguerre polynomials.}
\index{orthogonal polynomials of two variables!product
Hermite--Laguerre polynomials}%
For weight function $W(x,y) =  y^\a e^{-x^2- y}$, a possible 
orthogonal basis is given by $P_k^n(x,y) = H_k(x) L^\a_{n-k}(y)$, $0 \le k \le n$. This satisfies the differential 
equation
\begin{equation} \label{diff-Hermite-Laguerre}
     \thalf v_{xx} + y v_{yy} - x v_x + (1+\a -y) v_y  = - n v.
\end{equation} 
\end{enumerate}

There are other bases and further results for these product weight functions. These three cases 
are the only product type orthogonal polynomials that are eigenfunctions of a second order 
differential operators with eigenvalues depending only on $n$. See \S\ref{sec:PDE-2d} 
for further results. 

\subsection{Orthogonal polynomial on the unit disk} \label{sec:OP-disc}
\index{orthogonal polynomials of two variables!on unit disk}%
On the unit disk $B^2 := \big\{(x,y)\in\RR^2\mid x^2 + y^2 \le 1\big\}$
consider the weight function
\begin{equation}\label{weightdisk}
   W_\mu(x,y) := \frac{2\,\mu + 1}{2\,\pi}\,(1-x^2-y^2)^{\mu-\half}, \qquad \mu > -1/2, 
\end{equation}
normalized such that its integral over $B^2$ is 1. There are several
distinct bases of $\CV_n^2$ that can be given explicitly. 
\vskip\medskipamount
\begin{enumerate}[label=\arabic*.]
\item {\it First orthonormal basis}\quad
This is the basis $\big\{P_k^n\big\}_{k=0}^n$ of $\CV_n^2$ defined by
\begin{align} \label{basis1_ball}
P_k^n(x,y) &:= (h_{k,n})^{-1} C_{n-k}^{k+\mu+\half}(x)\,
(1-x^2)^{\frac{k}{2}}\,C_k^\mu \left(\frac{y}{\sqrt{1-x^2}} \right),
\\
 \label{basis1_ball_norm}
[h_{k,n}]^2 &:=
\frac{(2k + 2 \mu +1)_{n-k} (2\mu)_k (\mu)_k  (\mu+\half)}
     {(n-k)! \,k ! \,(\mu+\half)_k  (n+ \mu+\half)},
\end{align}
where the case $\mu=0$ can be obtained as a limit for $\mu\to0$ after
dividing, for $k>0$,  $C_k^\mu$ and $h_{k,n}$ by $\mu$ and by using
\eqref{limGegCheb}.
\item {\it Second orthonormal basis}\quad
In polar coordinates $(x,y) = (r\cos \t, r\sin \t)$,~let
\begin{align} \label{basis2_ball}
\begin{split} 
 P_{j,1} (x,y) & := [h_{j,n}]^{-1} P_{j}^{(\mu-\half,n-2j)}(2r^2 -1)r^{n-2j}
 \cos\big((n-2j)\t\big),\quad 1 \le j \le n/2, \\
 P_{j,2} (x,y) & := [h_{j,n}]^{-1} P_{j}^{(\mu-\half,n-2j)}(2r^2 -1)r^{n-2j}
 \sin\big((n-2j) \t\big),\quad 1 \le j < n/2,
\end{split}\\
\label{basis2_ball_norm}
   [h_{j,n}]^2 &:= \frac{(\mu+\half)_j (n-j)!\,(n-j+\mu+\half)}{ j!\,(\mu+\frac32)_{n-j} (n+\mu + \half)} \times
       \begin{cases} 2, & n \ne 2j, \\ 1, & n = 2j. \end{cases}
\end{align}
Then $\{P_{j,1}\}_{j=0}^{\lfloor n/2\rfloor}\cup
\{P_{j,2}\}_{j=0}^{\lfloor(n-1)/2\rfloor}$
is an orthonormal basis of $\CV_n^2$.
\item {\it Appell's biorthogonal polynomials}\quad
\index{Appell biorthogonal polynomials!on unit disk}%
\index{biorthogonal polynomials|seealso{Appell biorthogonl polynomials}}%
These are two families
$\{U_k^n\}_{k=0}^n$ and $\{V_k^n\}_{k=0}^n$ of bases of
$\CV_n^2$ that satisfy 
$$
     \int_{B^2} U_k^n(x,y) V_j^n(x,y) W_\mu(x,y)\,dx\,dy = h_k^n \delta_{k,j}, \quad 0 \le k,j \le n.
$$
The first basis is defined via the \emph{Rodrigues type formula}
\index{Rodrigues type formula}%
\begin{align*}
    U_k^n(x,y) :=  \frac{(-1)^n (2\mu)_n}{2^n (\mu+\half)_n n!} (1-x^2-y^2)^{-\mu+\half}    \frac{\partial^n }{\partial x^k \partial y^{n-k}} \left(
             (1-x^2-y^2)^{n+\mu-\half}\right).
\end{align*}
The second basis {is monic, up to
constant factors}: $V_k^n(x,y)={\rm const.}\,x^k y^{n-k}+
\mbox{polynomial of degree $<n$}$.
See for an explicit expression of $V_k^n(x,y)$ and further
properties of these two bases \S\ref{sec:OPball}, where they 
are given in the setting of the $d$-dimensional ball. 
\item  {\it An orthonormal basis of ridge polynomials for $\mu = 1/2$}\quad
\index{ridge polynomials on unit disk}
Let
\begin{equation} \label{basis3_ball}
   P_k^n(x,y) =\pi^{-1/2} U_n\left(x \cos \tfrac{k \pi}{n+1} + y \sin \tfrac{k \pi}{n+1}\right), 
       \quad 0 \le k \le n. 
\end{equation}
Then $\{P_k^n\}_{k=0}^n$ is an orthonormal basis of $\CV_n^2$ for the weight function 
$W_{1/2} (x,y)= \frac{1}{\pi}$ on $B^2$.
\end{enumerate}
\vskip\medskipamount
\noindent
\emph{Differential equation}\quad
All orthogonal polynomials of degree $n$ for $W_\mu$ are eigenfunctions of 
a second order differential operator. For $n \ge 0$, 
\begin{equation} \label{diff-disc}
 (1-x^2) v_{xx}  - 2xy v_{xy}  +  (1-y^2) v_{yy}   
   -  (2 \mu+2) (x v_x+ y v_y)
= - n (n+ 2\mu + 1) v, \quad v \in \CV_n^2.
\end{equation}
\noindent
\emph{Further results and references}\quad
For further properties, such as a 
closed formula for the reproducing kernel and convergence of orthogonal expansions, 
see \S\ref{sec:OPball}, where the disk will be a special case ($d=2$) of the $d$-dimensional ball.   
If the complex plane $\CC$ is identified with $\RR^2$, then the basis \eqref{basis2_ball} can be written 
in variables $z = x+ i y$ and $\bar z = x- i y$; see \S \ref{sec:complex}. 

The first orthonormal basis goes back as far as Hermite and was studied in \cite{ApKa}. Biorthogonal polynomials 
were studied in detail in \cite{ApKa}, see also \cite{Er2}. The basis
of ridge polynomials in \eqref{basis3_ball} was first discovered 
in \cite{LoSh}, see also \cite{Xu00b}, and it plays an important role in computer tomography \cite{LoSh, Marr, Xu06a}. 
For further studies on the orthogonal polynomials on the disk, see \cite{Waldron1, Wunsche}. 

\subsection{Orthogonal polynomials on the triangle} \label{sec:OPtriangle}
\index{orthogonal polynomials of two variables!on triangle}%
On the triangle $T^2 : = \big\{ (x,y)\mid 0 \le x, y, x+y \le 1\big\}$
consider an analogue of the Jacobi weight function
\begin{equation}\label{weighttriangle}
    W_{\a,\,\b,\g}(x,y) := \frac{\Gamma(\a+ \b+\g + \frac{3}{2})}
  {\Gamma(\alpha+\half)\Gamma(\beta+\half)
     \Gamma(\gamma+\half)}\, x^{\alpha - \half}y^{\beta-\half}(1-x-y)^{\gamma-\half},\quad \a, \b, \g > -1/2,
\end{equation}
normalized such that its integral over $T^2$ equals 1.
Several distinct 
bases of $\CV_n^2$ can be given explicitly.
\vskip\medskipamount
\begin{enumerate}[label=\arabic*.]
\item
{\it An orthonormal basis} $\{P_k^n\}_{k=0}^n$ of $\CV_n^2$ with  
\begin{align} \label{triangle-basis}
P_k^n(x,y)&:= [h_{k,n}]^{-1} P_{n-k}^{(2k+\beta+\gamma,\alpha-\half)}
(2x-1) (1-x)^k  P_k^{(\gamma -\half,\,\beta -\half)}\biggl( \frac{2y}{1-x} -1\biggl),\\
 [h_{k,n}]^2   &:=   \frac{(\a+\half)_{n-k}(\b+\half)_{k}(\g+\half)_{k} (\b+\g+1)_{n+k}}
    {(n-k)!\, k!\, (\b+\g+1)_k (\a+\b+\g +\frac32)_{n+k}}\, \frac{(n+k + \a+\b+\g+\half)(k+\b+\g)}
   {(2n+\a + \b+\g+\half)(2k+\beta+\gamma)}.\notag
\end{align}
Parametrizing the triangle differently leads to two more orthonormal bases. Indeed, 
denote the $P_k^n$ in \eqref{triangle-basis} by $P_{k,n}^{\a,\,\b,\g}$ and define 
$$
Q_k^n (x,y):= P_{k,n}^{\g,\,\b,\a}(1-x-y,y) \quad \hbox{and} 
   \quad R_k^n (x,y) := P_{k,n}^{\a,\g,\,\b}(x,1-x-y).
$$
Then $\{Q_k^n\}_{k=0}^n$ and $\{R_k^n\}_{k=0}^n$ are also
orthonormal bases.
\item
{\it Biorthogonal polynomials including Appell polynomials}\quad
\index{Appell biorthogonal polynomials!on triangle}%
A basis $\{U_k^n\}_{k=0}^n$ of $\CV_n^2$ due to Appell is
defined via the Rodrigues type formula:
\index{Rodrigues type formula}%
\begin{equation}\label{Un-ball-d=2}
U_k^n(x,y):=x^{-\a+\half}y^{-\b+\half}(1-x-y)^{-\g+\half}
\frac{\partial^n }{\partial x^k \partial y^{n-k}}
\left(x^{k+\a-\half}y^{n-k+\b-\half}(1-x-y)^{n+\g-\half}\right). 
\end{equation}
Biorthogonal to this basis is a basis $\{V_k^n\}_{k=0}^n$ of $\CV_n^2$:
$$
  \int_{B^2} U_k^n(x,y) V_j^n(x,y) W_{\a,\,\b,\g}(x,y)\,dx\,dy = h_k^n \delta_{k,j}, 
   \quad 0 \le k,j \le n.
$$
This is a monic basis, up to
constant factors: $V_k^n(x,y)={\rm const.}\,x^k y^{n-k}+
\mbox{polynomial of degree $<n$}$.
See for an explicit expression of $V_k^n(x,y)$ and further
properties of these two bases \S\ref{sec:OPsimplex}, where they are
given in the setting of the $d$-dimensional simplex.
\end{enumerate}
\vskip\medskipamount
\noindent
\emph{Differential equation}\quad
Orthogonal polynomials of degree $n$ with respect to $W_{\a,\b,\g}$ are eigenfunctions of a 
second order differential operator. For $n\ge0$,
\begin{multline} \label{diff-triangle}
   x(1-x)  v_{xx}  -2xy v_{xy} + y(1-y) v_{yy} -
   \left( (\a + \b + \g + \tfrac32) x - (\a+\thalf) \right)v_x\\
   - \left( (\a + \b + \g 
      + \tfrac32) y - (\b+\thalf)\right) v_y   =
      - n \left(n+ \a+\b+\g + \thalf \right) v, \qquad v \in \CV_n^2.
\end{multline}

\noindent
\emph{Further results and references}\quad
For further properties, such as biorthogonal polynomials, a closed formula for the reproducing 
kernel and convergence of orthogonal expansions, see
\S\ref{sec:OPsimplex},  
where the triangle will be the special case $d=2$ of the
$d$-dimensional simplex.   

The orthonormal polynomials in \eqref{triangle-basis} were first introduced in \cite{Proriol},
see \cite{Koorn75}; the case $\alpha = \beta = \gamma = 0$ became known as Dubiner's 
polynomials \cite{Dub} much later in the finite elements community. Appell polynomials were 
studied in detail in \cite{ApKa}. See \S\ref{sec:OPsimplex} for further references on orthogonal
polynomials on the simplex. 

\subsection{Differential equations and orthogonal polynomials of two variables} \label{sec:PDE-2d} 
\index{orthogonal polynomials of two variables!differential equation}%
A  linear second order partial differential operator 
\begin{equation} \label{admissL}
L : = A(x,y)\partial_1^2 + 2 B(x,y) \partial_1 \partial_2 + C(x,y) \partial_2^2 
    + D(x,y)\partial_1 + E(x,y) \partial_2,
\end{equation} 
where $\partial_1 : = \frac{\partial}{\partial x}$ and $\partial_2 : = \frac{\partial}{\partial y}$, is
called \emph{admissible}
\index{admissible partial differential operator}
if for each nonnegative integer $n$ there exists a number $\l_n$ such that  
the equation 
$$
        L u = \l_n u     
$$
has $n+1$ linearly independent solutions of polynomials of degree $n$ and has no nonzero  
solutions of polynomials of degree less than $n$. For $L$ in \eqref{admissL} to be admissible, its coefficients 
must be of the form
\begin{align*}
  A(x,y) &= A x^2 + a_1 x +b_1 y + c_1, \quad
  B(x,y) = A xy + a_2 x +b_2 y + c_2,\\
  C(x,y) &= A y^2 + a_3 x +b_3 y + c_3,\quad
  D(x,y) = B x + d_1, \quad E(x,y) = B y + d_2,
\end{align*}
and, furthermore, for each $n = 0,1,2,\ldots$, 
$$
    n A + B \ne 0,  \quad \hbox{and} \quad
    \l_n = - n \big((n-1)A + B\big). 
$$

A classification of the admissible equation that have orthogonal polynomials as eigenfunctions 
was given by Krall and Sheffer \cite{KraShe}. Up to affine transformations, there are only nine equations. Five of them 
admit classical orthogonal polynomials. These are

\smallskip\noindent\smallskip
\quad (1) product Hermite polynomials, see \eqref{diff-Hermite-2d}, \\ \smallskip
\quad (2) product Laguerre polynomials, see \eqref{diff-Laguerre-2d}, \\ \smallskip
\quad (3) product Hermite and Laguerre polynomials, see \eqref{diff-Hermite-Laguerre}, \\ \smallskip
\quad (4) orthogonal polynomials on the disk, see \eqref{diff-disc}, \\ \smallskip
\quad (5) orthogonal polynomials on the triangle, see \eqref{diff-triangle}. \\ \smallskip
The other four admissible differential equations are listed below.

\noindent \smallskip
\quad (6) $3y  v_{xx} + 2 v_{yy} - x v_x - y v_y = \l u$, \\ \smallskip
\quad (7) $(x^2+y+1) v_{xx} + (2xy+2x)v_{xy} + (y^2+2y+1) v_{yy} + g(x v_x + y v_y) = \l u$, \\ \smallskip
\quad (8) $x^2 v_{xx} + 2 xy v_{xy} + (y^2-y) v_{yy} + g[ (x-1) v_x + (y -\a) v_y] = \l u$, \\ \smallskip
\quad (9) $(x+\a) v_{xx} + 2(y+1)v_{yy} + x v_x + y v_y = \l u$. \\
The solutions for the last four equations are weak orthogonal polynomials in the sense that 
the polynomials are orthogonal with respect to a linear functional that is not positive definite. 

Another classification in \cite{Suetin}, based on \cite{Engelis}, listed fifteen cases, some of them 
are equivalent  under affine transformations in \cite{KraShe} but are treated separately because of 
other considerations. The orthogonality of the cases (6) and (7) is determined in \cite{KraShe}, 
while the cases (8) and (9) are determined in \cite{BSX1, KLL}. For further results, including 
solutions of the last four cases and further discussion on the impact of affine transformations, 
see \cite{Litt, Lyskova, KLL} and references therein.  Classical orthogonal polynomials in 
two variables were studied in the context of hypergroups in \cite{CoSc}. 

By the definition of the admissibility, all orthogonal polynomials of degree $n$ are 
eigenfunctions of an admissible differential operator for the same eigenvalue. In other words, 
the eigenfunction space for each eigenvalue is $\CV_n^2$. This requirement excludes, for 
example, the product Jacobi polynomial
$P_k^{(\a,\,\b)}(x)P_{n-k}^{(\g,\d)}(y)$, which satisfies
a second order equation of the form $L u = \l_{k,n} u$, where $\l_{k,n}$ depends on both 
$k$ and $n$. The product Jacobi polynomials, and other classical orthogonal polynomials of 
two variables, satisfy a second order matrix differential equation, see \cite{FPP} and
the references therein, and they also satisfy a matrix form of Rodrigue's type formula 
\cite{AFPP}. 

\subsection{Orthogonal polynomials of complex variables} \label{sec:complex}
\index{orthogonal polynomials of complex variables}%
Orthogonal polynomials of two real variables can be given
as polynomials of complex variables $z$ and $\bar z$ by
identifying $\RR^2$ with the complex plane $\CC$ and setting $z = x + i y$. For a real weight function
$W$ on $\Omega \in \RR^2$, we consider polynomials in $z$ and $\bar z$ that are orthogonal with 
respect to the inner product
\begin{equation} \label{eq:ipdC}
   \la f, g \ra_W^\CC: = \int_\Omega f(z,\bar z )\,
   \overline{g(z, \bar z)}\,w(z)\,dx\,dy, 
\end{equation}
where $w(z) = W(x,y)$. Let $\CV_n^2(W,\CC)$ denote the space of orthogonal polynomials in $z$ and $\bar z$ 
with respect to the inner product \eqref{eq:ipdC}. 
In this subsection, we denote by $P_{k,n}(x,y)$ real orthogonal polynomials 
with respect to $W$ and denote by $Q_{k,n}(z,\bar z)$ orthogonal polynomials in $\CV_n^2(W,\CC)$. 

\prop{ch2_thk7}
{The space $\CV_n^2(W,\CC)$ has a basis $Q_{k,n}$ that satisfies 
\begin{equation} \label{eq:Qconjugate}
Q_{k,n} (z,\bar z) = \overline{Q_{n-k,n}(z,\bar z)},\quad 0 \le k \le n.
\end{equation}
Furthermore, this basis is related to the basis of $\CV_n^2(W)$ by 
\begin{align} \label{eq:PvsQ}
\begin{split}
  P_{k,n}(x,y) & =
  \tfrac1{2^{1/2}} \left(Q_{k,n} (z, \bar z) + Q_{n-k,n} (z, \bar z) \right),
  \quad 0 \le k \le \thalf n, \\
  P_{k,n}(x,y) & =
  \tfrac{1}{2^{1/2} i } \left(Q_{k,n} (z, \bar z) - Q_{k-k,n} (z, \bar z) \right),
  \quad\thalf n < k \le n. 
\end{split}
\end{align}}
Writing orthogonal polynomials in terms of complex variables
often leads to more symmetric formulas. Below are two
examples. 
\vskip\medskipamount
\begin{enumerate}[label=\arabic*.]
\item \emph{Complex Hermite polynomials}\quad
\index{orthogonal polynomials of complex variables!complex Hermite
polynomials}%
For $j,k\in\NN_0$ define 
\begin{equation*}
  H_{k,j} (z,\bar z) := z^k {\bar z}^j
  \hyp20{-k,-j}{-}{\f{1}{z \bar z}}. 
\end{equation*}
Then $H_{k,j} \in \CV_{k+j}^2(w_H,\CC)$, where 
$w_H(z) = \pi^{-1} e^{-x^2-y^2} := \pi^{-1} e^{-|z|^2}$ with
$z = x+ i y \in \CC$.
These polynomials satisfy:
\begin{enumerate}[label=\roman*.]
\item $H_{k,j} (z,\bar z) = \overline{H_{j,k}(z,\bar z)}$;
\item  $H_{k,j} (z,\bar z) =
(-1)^j j!\, z^{k-j} L_j^{k-j}(|z|^2)$,\quad $k \ge j$\quad
($L_j^\a$ Laguerre polynomial, for $k\le j$ use (i)); 
\item $\displaystyle { \frac{\partial} {\partial z} H_{k,j}  = \bar{z} H_{k,j}  - H_{k,j+1}, \quad
   \frac{\partial} {\partial \bar z} H_{k,j}   = z H_{k,j}   - H_{k+1,j} }$; 
\item $z H_{k,j}   = H_{k+1,j}  + j H_{k,j-1}, \quad
\bar{z} H_{k,j}  = H_{k,j+1}  + k H_{k-1,j}$; 

\item $\displaystyle\int_\CC H_{k,j}(z,\bar z)\,
\overline {H_{m,l}(z,\bar z)}\, w_H(z)\,dx\,dy = j!\, k!\,
\delta_{k,m}\,\delta_{j,l}$.
\end{enumerate}
\vskip\medskipamount
\item \emph{Disk polynomials}\quad
\index{orthogonal polynomials of complex variables!disk polynomials}%
\index{disk polynomials}%
For $j,k \in \NN_0$ define 
$$
   P^\l_{k,j}(z,\bar z) := \frac{(\l+1)_{k+j}}{(\l+1)_k(\l+1)_j} 
    z^k \bar z^j 
    \hyp21{-k,-j}{-\l-k-j}{\frac{1}{z \bar z}},
$$
normalized by $P^\l_{k,j}(1,1) =1$.
Then $P^\l_{\k,j} \in \CV_{k+1}^2(w_\l,\CC)$,
where $w_\l(z) := \frac{\l+1}{\pi}(1-|z|^2)^\l$, $\l > -1$.
These polynomials satisfy 

\begin{enumerate}[label=\roman*.]
\item
$P_{k,j}^\l(z,\bar z)=\overline{P_{j,k}^\l(z,\bar z)}$;
\item
$\displaystyle P^\l_{k,j}(z,\bar z) = 
   \frac{j!}{(\l+1)_j} P_{j}^{(\l,k-j)}(2 |z|^2 -1) z^{k-j}$,
   \quad $k\ge j$\quad
($P_j^{(\a,\b)}$ Jacobi polynomial,\\
for $k \le j$ use (i));
\item
$|P^\l_{k,j}(z,\bar z)| \le 1$ for $|z| \le 1$ and $\l \ge 0$; 
\item
$z P^\l_{k,j} = \dfrac{\l + k + 1}{\l + k + j +1}
  P^\l_{k+1,j}+ \dfrac{j}{\l + k + j +1}
  P^\l_{k,j-1}$ \\
and a similar relation holds for $\bar z P^\l_{k,j}$ upon using (i);

\item $\displaystyle{ \int_{D} P^\l_{k,j}(z,\bar z)\,\overline{
 P^\l_{m,l}(z,\bar z)}\,w_\l(z)\,dx\,dy 
 = \frac{\l+1}{\l+k+j+1}\,\frac{k! j!}{(\l+1)_k(\l+1)_j} }\,
  \delta_{k,m} \delta_{j,l}.$
\end{enumerate}
\end{enumerate}
\vskip\medskipamount

The complex Hermite polynomials were introduced by It\^o \cite{Ito}.
They have been widely studied by 
many authors, see \cite{Ghanmi, Int-int, Ismail} for some recent studies and the references therein. 
Disk polynomials were introduced by Zernike \cite{Ze} for
$\a = \b = 1/2$, and in a subsequent paper together with Brinkman
\cite{ZeBr} in general. They
are also called \emph{Zernike polynomials}
\index{Zernike polynomials}%
and they
have applications in optics. They were used in \cite{Folland} to expand 
the Poisson--Szeg\H{o} kernel for the ball in $\CC^d$. A Banach algebra related to disk polynomials was studied 
in \cite{Kanjin}. For further properties of disk polynomials, including the fact that for 
$\l = d-2$, $d= 2, 3, \ldots$, they are spherical functions for
$\mathrm{U}(d)/\mathrm{U}(d-1)$, see 
\cite{Ikeda,Koorn75} and \cite{ViKl2, Wunsche}. 

The structure of complex orthogonal polynomials of two variables and its connection and contrast to 
its real counterpart was studied in \cite{Xu15a}.

\subsection{Jacobi polynomials associated with root systems and related orthogonal polynomials}\label{sec:Koornwinder}
\begin{figure}[b] 
\includegraphics[scale=0.31]{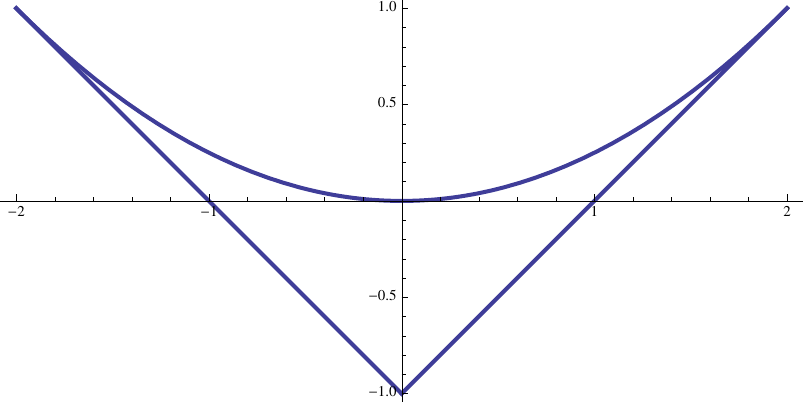} \qquad\qquad
\includegraphics[scale=0.31]{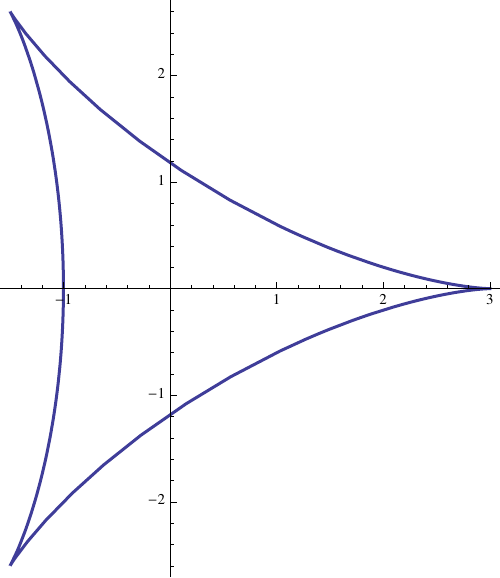} 
\caption[Domains of Koornwinder orthogonal polynomials] 
{Regions for Koornwinder orthogonal polynomials. Left: $\mathrm{BC}_2$.
Right: $\mathrm A_2$.} 
\label{figure:region} 
\end{figure}

\index{Jacobi polynomials!associated with root systems}%
There are two families of Jacoobi polynomials of two variables
associated with a
root system and a related family of orthogonal polynomials. 
\\[\medskipamount]
\emph{Jacobi polynomials associated with root system $\mathrm{BC}_2$}\quad
\index{Jacobi polynomials!$\mathrm{BC}_2$}%
\index{orthogonal polynomials of two variables!$\mathrm{BC}_2$}%
Consider the weight function 
\begin{equation} \label{1stKoorn-weight}
  W_{\a,\b,\g}(x,y): = (1-x+y)^\a (1+ x + y)^\b (x^2 - 4 y)^\g, 
\end{equation}
where $\a,\b,\g > -1$, $\a+\g > - \tfrac32$ and $\b+ \g > - \tfrac32$, defined on the domain
\begin{equation} \label{1stKoorn-domain}
   \Omega: = \big\{(x,y)\mid |x| < y+1, x^2 > 4 y\big\}, 
\end{equation}
which is depicted in the left part of Figure \ref{figure:region}.
After a change of variables $x = u+v$, $y= u v$
the domain and weight function become
\begin{equation}\label{eq:W*(u,v)}
\begin{split}
    \Omega^* &= \big\{(u,v)\mid -1 < u < v < 1\big\},  \\
     W_{\a,\b,\g}^*(u,v)  & = (1-u)^\a (1+u)^\b (1-v)^\a (1+v)^\b (v-u)^{2 \g +1}. 
\end{split}
\end{equation}
Let ${\mathcal N} :=\big\{(n,k)\mid 0 \le k \le n\big\}$. In ${ \mathcal N}$ define $(j,m) \prec (k,n)$ 
if $m < n$ or $m = n$ and $j \le k$ (graded lexicographic order).
Then an orthogonal polynomial 
$P_k^n$ that satisfies 
\begin{equation}\label{eq:jacobi-biangle}
   P_k^n(x,y) = x^{n-k} y^k + \sum_{(j,m) \prec (k,n)} a_{j,m} x^{m-j} y^j 
\end{equation}
and is orthogonal to all $x^{m-j} y^j$ for $(j,m) \prec (k,n)$ is uniquely determined. 

In the case $\g =   \pm \half$, a basis of orthogonal polynomials can 
be given explicitly. In fact, such a basis can be given
in the more general case
where we twice replace in 
\eqref{eq:W*(u,v)} the Jacobi weight function by an arbitrary weight 
function $w$. Then
\begin{equation} \label{1stKoorn-weight-2}
  W_{\pm \half} (x,y)  = w(u)w(v) (x^2 - 4 y)^{\pm \half}\quad\hbox{with}\quad x = u+v, \,\, y = uv, 
\end{equation}
defined on the domain
$\big\{(x,y)\mid  x^2  > 4 y, \; u, v \in \operatorname{supp}  w\big\}$ 
for any weight function $w$ on $\RR$. Let $\{p_n\}$ denote an orthonormal 
basis with respect to $w$. Then an orthonormal basis of polynomials for $W_{\pm\half}$ is given by
\index{orthogonal polynomials of two variables!symmetrized products}%
\index{orthogonal polynomials of two variables!anti-symmetrized products}%
\begin{align} \label{1stKoorn-poly-1}
P_{k}^{n, (-\half)}(x,y)&=\begin{cases} 
  p_n(u)p_k(v)+p_n(u)p_k(v),& k<n, \\
  2^{1/2} p_n(u)p_n(v),&  k=n,\end{cases}\\
\label{1stKoorn-poly-2}
 P_{k}^{n,(\half)}(x,y)&= \frac{p_{n+1}(u)p_k(v)-p_{n+1}(v)p_k(u)}
 {u-v},   \quad 0 \le k \le n,
\end{align}
where in both cases $(u,v)$ is related to
$(x,y)$ by $x = u+v, \,\, y = uv$. 
\\[\medskipamount]
\emph{A related family of orthogonal polynomials}\quad\\
Consider the family of weight functions 
defined by 
\begin{equation} \label{CWabc}
   \CW_{\a,\b,\g} (x,y) : =  |x-y|^{2\a +1}  |x + y|^{2\b +1}
   (1-x^2)^\g (1-y^2)^\g,\quad (x,y) \in [-1,1]^2,
\end{equation}
 where $\a, \b, \g > -1$, $\a + \g + \frac32> 0$ and $\b+ \g + \frac32> 0$.
These weight functions are related to those in \eqref{1stKoorn-weight} by
$$
  \CW_{\a,\b,\g}(x,y) = 4^\g\,\big|x^2 - y^2\big| \,W_{\a,\b,\g}(2xy, x^2+y^2 -1). 
$$
Let $\big\{P_{k,n}^{\,\a,\b,\g}\big\}_{k=0}^n$ denote an orthogonal basis of
$\CV_n^2$ for $W_{\a,\b,\g}$. Then an orthogonal basis of $\CV_{2n}^2$ 
for $\CW_{\a,\b,\g}$ is given by
\begin{align} \label{Qeven}
\begin{split}
      & P_{k,n}^{\a,\b,\g}(2xy, x^2+y^2 -1), \qquad 0 \le k \le n, \\
     & (x^2-y^2)  P_{k,n-1}^{\a+1,\b+1,\g}(2xy, x^2+y^2 -1),  \,\, 0 \le k \le n-1, 
 \end{split}
\end{align}
and an orthogonal basis of $\CV_{2n+1}^2$ for $\CW_{\a,\b,\g}$ is given by
\begin{align} \label{Qodd}
 \begin{split}
       &  (x+ y) P_{k,n}^{\a,\b+1,\g}(2xy, x^2+y^2 -1), \quad 0 \le k \le n, \\
       &  (x-y) P_{k,n}^{\a+1,\b,\g}(2xy, x^2+y^2 -1), \quad  0 \le k \le n.
\end{split}
\end{align}
In particular, when $\g = \pm \f12$, the basis can be given in terms of the Jacobi  
polynomials of one variable by using \eqref{1stKoorn-poly-1} and \eqref{1stKoorn-poly-2}.
\paragraph{Jacobi polynomials associated with root system $\mathrm A_2$}\quad\\
\index{Jacobi polynomials!$\mathrm A_2$}%
\index{orthogonal polynomials of two variables!$\mathrm A_2$}%
These polynomials are orthogonal 
with respect to the weight function 
\begin{equation}
W_\a(x,y) : = \left [-3(x^2+y^2 + 1)^2 + 8 (x^3 - 3 xy^2)  +4\right]^{\a}
\end{equation}
on the region bounded by $-3(x^2+y^2 + 1)^2 + 8 (x^3 - 3 xy^2)  + 4 =0$, which is
called \emph{Steiner's hypocycloid}
\index{Steiner's hypocycloid}%
and can be described as the curve  
$$
  x + i y = (2 e^{ i \t} + e^{- 2 i \t})/3. \qquad  0 \le \t \le 2 \pi.
$$ 
This three-cusped region is depicted in the right part of
Figure \ref{figure:region}. 
Apart from $\a = \pm \half$, orthogonal polynomials with respect to $W_\a$ are not 
explicitly known. In the case of $\a = \pm \half$, a basis of orthogonal polynomials 
can be given in homogeneous coordinates as follows. Let 
$$
\RR^3_H: = \big\{ \tb = (t_1,t_2,t_3) \in \RR^3\mid  t_1 +t_2 +t_3 =0\big\}.
$$ 
For $\tb \in \RR^3_H$ and $\kb = (k_1, k_2, k_3) \in \RR^3_H \cap \ZZ^3$, define 
$\phi_\kb(\tb)= e^{\tfrac{2 \pi i}{3} \kb \cdot \tb}$ and 
\begin{equation*} 
 \sC^\pm_\kb(\tb) : = 
    \tfrac{1}{6} \Big( \phi_{k_1,k_2,k_3}(\tb)+ \phi_{k_2,k_3,k_1}(\tb)+
     \phi_{k_3,k_1,k_2}(\tb) 
  \pm \phi_{-k_1,-k_3,-k_2}(\tb) \pm \phi_{-k_2,-k_1,-k_3}(\tb) 
      \pm \phi_{-k_3,-k_2,-k_1}(\tb)\Big), 
\end{equation*}
which are analogues of cosine and sine functions. 
The region bounded by Steiner's hypocycloid is the image of the
triangle $\Delta = \big\{(t_1,t_2)\in\RR^2\mid t_1, t_2 \ge 0, \, t_1+t_2 \le 1\big\}$ under
the change of variables $(t_1,t_2) \mapsto (x,y)$, defined by 
\begin{equation} \label{z} 
z := x+ i y  =  \sC^+_{0,1,-1}(\tb)  = \tfrac{1}{3} \Big(\phi_{0,1,-1}(\tb) + \phi_{1,-1,0}(\tb) 
     + \phi_{-1,0,1}(\tb)\Big),
\end{equation}
Under the change of variables \eqref{z}, define 
\begin{align*}
T_k^n(z,\bar z) : & = \sC^+_{k, n-k, - n}(\tb), \qquad\qquad\;\;
0 \le k \le n, \\
U_k^n(z,\bar z) : & = \frac{\sC^-_{k+1,n-k+1,-n-2}(\tb)}{\sC^-_{1,1,-2}(\tb)},
                        \qquad 0 \le k \le m.
\end{align*}
Then $\{T_k^n\}_{k=0}^n$ and $\{U_k^n\}_{k=0}^n$ are bases of
$\CV_n^2$ with respect to $W_{-\half}$ and $W_{\half}$, 
respectively. Both families of polynomials satisfy the relation 
$P_{n-k}^n(z, \bar z)  = \overline{P_k^n(z, \bar{z})}$, so that real-valued orthogonal
bases can be derived from their real and imaginary parts.  These polynomials 
are analogues of Chebyshev polynomials of the first and the second kind.
\index{Chebyshev polynomials!$\mathrm A_2$}%
They  satisfy  the following three-term relations: 
\begin{align} \label{recurT}
     P _k^{m+1} (z,\bar z) = 3 z P_{k}^m(z,\bar z) -
                 P_{k+1}^m(z,\bar z) - P_{k-1}^{m-1}(z,\bar z) 
\end{align}
for $ 0 \le k \le m$ and $m\ge 1$, where  
\begin{align*}
&T_{-1}^m(z,\bar z ) : = T_1^{m+1}(z,\bar z ), \quad T_{m+1}^m(z,\bar z )
  : = T_{m}^{m+1}(z,\bar z ),\\
&U_{-1}^m(z,\bar z ) : =  0, \quad  U_m^{m-1}(z,\bar z ): =0,
\end{align*}
and, moreover, 
\begin{align*}
& T_0^0(z,\bar z)=1, \quad T_0^1(z,\bar z)= z,  \quad T_1^1(z,\bar z)= \bar z,\\
& U_0^0(z,\bar z)=1, \quad U_0^1(z,\bar z)= 3z,  \quad U_1^1(z,\bar z)= 3 \bar z.
\end{align*}
\emph{Further results and references} \quad
The Jacobi polynomials associated with $\mathrm{BC}_2$ and $\mathrm A_2$ were
initially
studied by Koornwinder \cite{Koorn74a, Koorn74b}
(see also \cite{Koorn75}), where it was shown that the orthogonal polynomials for $W_{\a,\b,\g}$ in 
\eqref{1stKoorn-weight} are eigenfunctions of two commuting differential operators of second and fourth order, 
whereas the orthogonal polynomials associated with $\mathrm A_2$ are
eigenfunctions of two commuting differential
operators of second and third order.
The two families are rank two cases of the Jacobi polynomials
associated with root systems for general rank, the study of which was
initiated by Heckman and Opdam, see Ch.~7.
The Jacobi polynomials associated with $\mathrm A_2$ can be identified
with \emph{Jack polynomials} of two variables.
\index{Jack polynomials}%

The special case of $P_0^n(x,y)$ of the first family when $\g = 1/2$ was 
studied also in \cite{EiLi74}. For further results on the first family, see \cite{KoornSpr} and \cite{Spr-K}, including 
explicit formula for $P_k^n$ in \eqref{eq:jacobi-biangle} given in terms of power series, and Rodrigues
type formulas, and \cite{Xu12} where an explicit formula for the reproducing kernel in the case 
of $W$ in \eqref{1stKoorn-weight-2} with $\g = \pm 1/2$ was given in terms of the reproducing 
kernels of the orthogonal polynomials of one variable. For further results on the second family, 
see \cite{Shishkin, Suetin} and \cite{LSX1}, the latter one includes
a connection with translation 
tiling and convergence of orthogonal expansions for $\mu = \pm \half$. 

Orthogonal polynomials with respect to $W_{\pm \f12}$ in \eqref{1stKoorn-weight-2} and
the Jacobi polynomials associated with $\mathrm A_2$ when $\g = \f12$
are remarkable for having the maximum
number of common zeros, i.e.,
$\PP_n =\big \{P_k^n\mid 0 \le k \le n\big\}$ has $\dim \Pi_{n-1}^2$ distinct 
real zeros (\cite{LSX1, SchmXu}). By Theorem~\ref{thm:GaussCuba}, Gaussian cubature formulas exist for 
these weight functions. For their generalizations to higher dimension, see \S\S\ref{sec:SymmFuncOP}, \ref{sec:AdOP}.   

The orthogonal polynomials with respect to $\CW$ in \eqref{CWabc} were studied in \cite{Xu12}.
The reproducing kernel of $\Pi_n^2$ in $L^2(\CW)$ can also be expressed in terms of the 
reproducing kernels of $\Pi_n^2$ in $L^2(W)$ for $W$ in \eqref{1stKoorn-weight}. In particular,
in the case of $\g = \pm \half$, the kernel can be expressed in terms of the reproducing kernels
of the Jacobi polynomials. When $\g = \pm \half$, these weight functions admit minimal cubature 
rules that attain the lower bound \eqref{eq:lwbd3}. 

\subsection{Methods of constructing orthogonal polynomials of two variables
from one variable} \label{sec:MethodOP}
\index{orthogonal polynomials of two variables!constructed from
one-variable polynomials}%
Let $w_1$ and $w_2$ be weight functions defined on the intervals $(a,b)$ and $(c,d)$,
respectively. Let $\rho$ be a positive function on $(a,b)$. For the 
weight function 
\begin{equation}\label{eq:KoWeight}
 W(x,y) := w_1(x) w_2(\rho^{-1}(x)y), \qquad (x,y) \in \Omega,
\end{equation}
where the domain $\Omega$ is defined by
\begin{equation} \label{eq:KoDomain}
 \Omega := \{(x,y)\in\RR^2\mid a < x <b, \,\,\, c\rho(x) <y < d\rho(x)\}, 
\end{equation}
a basis of orthogonal polynomials of two variables can be given in terms of orthogonal 
polynomials of one variable whenever either one of the following additional assumptions
is satisfied: 
\\[\medskipamount]
{\it Case 1.}\; $\rho$ is a polynomial of degree $1$; \\ 
{\it Case 2.}\; $\rho = \sqrt{q}$ with $q$ a nonnegative polynomial of degree at most $2$,
and further assume that $c = - d >0$ and $w_2$ is an even function on
$(-c,c)$. 
\\[\medskipamount]\indent
For each $k \in \NN_0$ let $\{p_{n,k}\}_{n=0}^\infty$ denote the system of 
orthonormal polynomials with respect to the weight function
$\big(\rho(x)\big)^{2k+1}w_1(x)$ 
on $(a,b)$. And let $\{q_n\}$ be the system of orthonormal polynomials with respect 
to $w_2(x)$ on $(c,d)$. Define polynomials of two variables by 
\begin{equation} \label{eq:KoOP}
P_k^n(x,y) : = p_{n-k,k}(x)\,\big(\rho(x)\big)^k\,q_k\left(\frac{y}{\rho(x)}\right),
\quad 0\le k \le n. 
\end{equation}
In Case 2 we see that $P_k^n$ are polynomials of degree $n$
because $q_k$ has the same parity as $k$ by evenness of $w_2$.
Then  $\{P_k^n\}_{k=0}^n$ is an orthonormal basis of $\CV_n^2$
with respect to 
$W$ on $\Omega$.

Examples of orthogonal polynomials constructed by this method
include product 
orthogonal polynomials, for which $\rho(x) =1$,
\index{orthogonal polynomials of two variables!product orthogonal
polynomials}%
and also the
following cases:
\\[\medskipamount]
{\it Jacobi polynomials on the disk}\quad
\index{orthogonal polynomials of two variables!on unit disk}%
Let $w_1(x) = w_2(x):=(1-x^2)^{\mu-1/2}$ on 
$[-1,1]$ and $\rho(x) := (1- x^2)^{1/2}$. Then the weight function
\eqref{eq:KoWeight} 
and the basis \eqref{eq:KoOP} 
coincide up to constant factors with  \eqref{weightdisk}
and \eqref{basis1_ball}, respectively,
\\[\medskipamount]
{\it Jacobi polynomials on the triangle}\;\;
\index{orthogonal polynomials of two variables!on triangle}%
Let $w_1(x):= x^{\a-\half}(1-x)^{\b+ \g -1}$ 
and $w_2(x):= x^{\b-\half}(1-x)^{\g-\half}$, both defined on the interval
$(0,1)$, 
and let $\rho(x) := 1- x$.
Then the weight function \eqref{eq:KoWeight}
and the basis \eqref{eq:KoOP}
coincide up to constant factors with \eqref{weighttriangle}
and \eqref{triangle-basis}, respectively.
\\[\medskipamount]
{\it Orthogonal polynomials on the parabolic domain}\quad
\index{orthogonal polynomials of two variables!on parabolic domain}%
Let $w_1(x) := x^a(1-x)^b$ on $[0,1]$, $w_2(x) := (1-x^2)^a$ on $[-1,1]$, 
and $\rho(x) := \sqrt{x}$. Then the weight function \eqref{eq:KoWeight} becomes 
\begin{equation} \label{parabolic-domaiin}
  W_{a,b}(x,y) : = (1-x)^b (x-y^2)^b, \qquad y^2 < x< 1.
\end{equation}
The domain $\big\{(x,y)\mid y^2 < x < 1\big\}$ is bounded by a straight line and a 
parabola. The orthogonal polynomials $P_k^n$ in \eqref{eq:KoOP} are 
\begin{equation} \label{parabolic-domaiin-OP}
  P_k^n(x,y) = P_{n-k}^{(a,b+k+1/2)}(2x-1)\,x^{k/2}\,
  P_k^{(b,b)}\big(yx^{-1/2}\big), 
    \quad 0 \le k \le n.
\end{equation} 

\paragraph{Further results and references}\quad
This method of generating orthogonal polynomials of two variables first appeared in 
\cite{Larcher} and was used in \cite{Aga} in certain special cases. It was presented 
systematically in \cite{Koorn75}, where the two cases for $\rho$ were stated. 
For further examples of explicit bases constructed in various domains,
such as 
$\big\{(x,y)\mid x^2 + y ^2 \le 1,\; -a \le y \le b\big\}$  ($0< a,b < 1$), see
\cite{Suetin}. 

The sequence of polynomials in \eqref{parabolic-domaiin-OP} satisfies a product 
formula \cite{KoornSchw} that generates a convolution structure for $L^2(W_{\a,\b})$, 
which was used to study the convergence of orthogonal expansions \cite{CFX}. 

\subsection{Other orthogonal polynomials of two variables} \label{sec:OtherOP2d}   

This subsection contains miscellaneous results on orthogonal polynomials of two 
variables.
 
\begin{enumerate}[label=\arabic*.]
\item
\emph{Orthogonal polynomials for radial weight function}.\quad
\index{orthogonal polynomials of two variables!for radial weight function}%
A weight function $W$ is  called \emph{radial}
\index{weight function!radial}%
if it is of the form $W(x,y) = w(r)$,
where $r=\sqrt{x^2 +y^2}$. For such
a weight function, an orthonormal basis can be given in polar coordinates
$(x,y) = (r \cos \t, r \sin \t)$. Let $p_m^{(k)}$ denote the orthogonal polynomial of
degree $m$ with respect to the weight function $r^{k+1} w(r)$ on $[0,\infty)$. Define
\begin{align} \label{OP-radial-d2}
\begin{split}
  P_{j,1} (x,y)  & = p_{2j}^{(2n-4j+1)}(r) r^{(n-2j)}\cos\big((n-2j)\t\big),
  \quad 0 \le j \le n/2, \\
  P_{j,2} (x,y)  & = p_{2j}^{(2n-4j+1)}(r) r^{(n-2j)}\sin\big((n-2j)\t\big), 
  \quad 0 \le j < n/2. 
\end{split}
\end{align}
Then $\{P_{j,1}\}_{j=0}^{\lfloor n/2\rfloor}\cup
\{P_{j,2}\}_{j=0}^{\lfloor(n-1)/2\rfloor}$
 is an orthogonal
basis of $\CV_n^2$ with respect to $W$. For $W(x,y) = (1-r^2)^{\mu-1/2}$
this is the basis given in \eqref{basis2_ball}.
\index{orthogonal polynomials of two variables!on unit disk}
Another example is the following:

\medskip\noindent
{\it Product Hermite weight function $W(x,y) = e^{-x^2-y^2}$}.
\index{orthogonal polynomials of two variables!for product Hermite weight
function}%
The basis \eqref{OP-radial-d2} 
is given by
\begin{align} \label{Hermite-radial-d2}
\begin{split}
  P_{j,1} (x,y)  & = L_{j}^{n-2j}(r^2) r^{(n-2j)}\cos\big((n-2j)\t\big),
  \quad 0 \le j \le n/2, \\
  P_{j,2} (x,y)  & = L_{2j}^{n-2j}(r^2) ar^{(n-2j)}\sin\big((n-2j)\t\big), 
  \quad 0 \le j < n/2 
\end{split}
\end{align}
in terms of Laguerre polynomials. 

\item
\emph{Bernstein--Szeg\H{o} weight function of two variables}.\quad
\index{orthogonal polynomials of two variables!for Bernstein--Szeg\H{o}
weight function}%
\index{Bernstein--Szeg\H{o} weight function (of two variables)}%
 Let $m\in\NN$. 
For any $i=0,\dots,m$ let $h_i(y)$ be polynomials in $y$ with real
coefficients of degree at most 
$\frac{m}{2} - \big|\frac{m}{2}-i\big|$, with $h_0(y)=1$, such that, for all
$y\in[-1,1]$,  
\begin{equation}\label{Bernstein-szego1}
    h(z,y)=\sum_{i=0}^{m}h_i(y)z^i, \qquad z \in \CC, 
\end{equation}
is nonzero whenever $|z|\leq 1$. Consider the two variable weight function
\begin{equation}\label{Bernstein-szego2}
   W(x,y):=\frac{4}{\pi^2}\,\frac{\sqrt{1-x^2}\sqrt{1-y^2}}{|h(z,y)|^2}, \qquad
   x = \thalf(z+z^{-1}).
\end{equation}
For $\lceil \frac{m -2}{2} \rceil\le k \le n$, define the polynomials 
\begin{equation}\label{Bernstein-szego-OP}
   P_k^n(x,y) = U_{n-k}(y) \sum_{i=0}^{m} h_i(y)U_{k-i}(x), 
\end{equation} 
where it is understood that $U_n(x)=-U_{-n-2}(x)$ if $n<0$. 
Then $P_k^n$ is an orthogonal polynomial of degree $n$ with respect to $W$. In particular, 
if $m \le 2$, then $\{P_k^n\}_{k=0}^n$ is an orthogonal basiof $\CV_n^2$.  

\medskip\noindent
{\it Example}\quad 
$h(z,y) =  1- 2 a y z + a^2 z^2$, where $|a| <1$ and $a \in \RR$. Then
\begin{equation}\label{Bernstein-szego3}
   W(x,y)  =    \frac{4}{\pi^2}\,\frac{(1-a^2) \sqrt{1-x^2}\sqrt{1-y^2}}
{4a^2(x^2+y^2)-4a(1+a^2)xy+(1-a^2)^2},
 \end{equation}
The orthogonal polynomials in \eqref{Bernstein-szego-OP} are given,
up to a constant, by
$P_0^n(x,y) = U_n(y)$ and 
$$
P_k^n(x,y) = \left(U_k(x)-2ayU_{k-1}(x)+a^2U_{k-2}(x)\right)U_{n-k}(y), 
\quad 1 \le k \le n.
$$
For further examples and other properties of such polynomials, see
\cite{DGIX}.  

\item
\emph{Orthogonal polynomials on the regular hexagon}.\quad
Orthogonal polynomials with respect to the constant weight
function on the regular hexagon were studied in \cite{Dunkl87}, 
An algorithm was given there for generating an orthogonal basis.
No closed form of such a basis is known. 

\section{Spherical harmonics} \label{sec:harmonics}
Here and later we will use notation $\|x\|:=(x_1^2+\cdots+x_d^2)^{1/2}$
and $\la x,y\ra:=x_1y_1+\cdots+x_dy_d$ ($x,y\in\RR^d$).
Spherical harmonics are an essential tool for Fourier analysis on the 
unit sphere 
$\SS^{d-1}:=\big\{x\in\RR^d\mid\|x\|=1\big\}$ in $\RR^d$ ($d\ge2$).
They are also building blocks for
families of orthogonal polynomials with respect to
radial weight functions on $\RR^d$.
\end{enumerate}
 
 \subsection{Ordinary spherical harmonics}
 
Let $\Delta := \frac{\partial^2}{\partial x_1^2} + \cdots +
\frac{\partial^2}{\partial x_d^2}$ 
be the \emph{Laplace operator} on $\RR^d$. A polynomial $Y$ on $\RR^d$
is called \emph{harmonic}
if $\Delta  Y =0$. For $n=0,1,2,\ldots $ let $\mathcal{H}_{n}^{d}$ denote the linear space 
of homogeneous harmonic polynomials of degree $n$ in $d$ variables,
i.e.,
$$
\CH_{n}^{d} :=\left\{ P\in \CP_{n}^{d}\mid\Delta P=0\right\}. 
$$
By definition, a \emph{spherical harmonic} is the restriction of a
homogeneous harmonic
polynomial to the unit sphere.
\index{spherical harmonics}%
If $Y \in \CH_n^d$, then $Y(x) = \|x\|^n Y(x')$ where $x ' = x /\|x\| \in \SS^{d-1}$. 
We shall also use $\CH_n^d$ to denote the space of spherical harmonics of degree $n$. 
For $n \in \NN_0$, 
\begin{equation}
\label{dimspherharm}
\dim \CH_{n}^{d}= \dim \CP_n^d - \dim \CP_{n-2}^d =
\binom{n+d-1}{d-1}-\binom{n+d-3}{d-1}=
\frac{(2n+d-2)\,(n+d-3)!}{n!\,(d-2)!}.
\end{equation}

Spherical harmonics of different degrees are orthogonal in
$L^2(\SS^{d-1}, d\sigma)$,
where $d\sigma$ denotes the normalized spherical measure on
$\SS^{d-1}$.
There is a unique decomposition 
$$
   \CP_n^d = \bigoplus_{ j=0}^{\lfloor \frac{n}{2} \rfloor} \|x\|^{2j} \CH_{n-2j}^d: \qquad
    P(x) = \sum_{j=0}^{\lfloor \frac{n}{2} \rfloor} \|x\|^{2j} Y_{n-2j}^d(x'). 
$$ 
\paragraph{Orthonormal basis}\quad
In terms of \emph{spherical polar coordinates}
\index{spherical polar coordinates}%
\begin{equation} \label{eq:sph-polar}
\begin{cases}
\quad x_{1} = r\sin \theta _{d-1}\ldots \sin \theta _{2}\sin \theta _{1},  \\
\quad x_{2} = r\sin \theta _{d-1}\ldots \sin \theta _{2}\cos \theta _{1}, \\
\quad  \cdots  \\
x_{d-1} = r\sin \theta _{d-1}\cos \theta _{d-2}, \\
\quad x_{d} =r\cos \theta _{d-1},
\end{cases}
\end{equation}
where $r\geq 0$, $0\leq \theta_{1}<2\pi$, $0\leq \theta_{i}\leq \pi$ for 
$i =2, \ldots, d$, the normalized measure $d\sigma$ on $\SS^{d-1}$ is
given by 
\begin{align} \label{eq:d-omega}  
d\sigma   :=  \o_d^{-1} \prod_{j=1}^{d-2}\big( \sin \theta_{d-j}\big)^{d-j-1}
 d\theta_{1}\,d\theta_{2}\ldots d\theta _{d-1},  
\end{align}
where $\o_{d} :=\frac{2 \pi^{d/2}}{\Gamma(d/2)}$ is the surface area of $\SS^{d-1}$. For $d =2$, 
$\dim \CH_n^2 = 2$. Then an orthogonal basis for $\CH_n^2$ is given
in polar coordinates by
$$
   Y_n^{(1)} (x) := r^n \cos(n \theta), \qquad Y_n^{(2)}(x):= r^n \sin(n \theta). 
$$
For $d >  2$ and $\a \in \NN_0^d$, define 
\begin{equation}
Y_{\alpha}(x) := [h_\alpha] ^{-1} r^{|\alpha|} g_{\alpha}(\theta_1)
 \prod_{j=1}^{d-2}  ( \sin \theta _{d-j})^{|\alpha^{j+1}|}\,
  C_{\alpha_j}^{\l_{j}} (\cos \theta _{d-j}),
\end{equation}
where $g_{\alpha}(\theta_1) := \cos ( \alpha_{d-1}\theta _1)$
(if $\alpha_d=0$) and $:=\sin ( (\alpha_{d-1}+1) \theta_1)$
(if $\alpha_d=1$),  where
$|\alpha^j| := \alpha_j+ \ldots+ \alpha_d$, $\l_j := |\alpha^{j+1}| + (d-j-1)/2$,
and where
$$
  [ h_\alpha]^{2}:=  b_{\alpha}
    \prod_{j=1}^{d-2}\frac{\alpha_j!\,\big(\thalf(d-j+1)\big)_{|\alpha^{j+1}|}\,( \alpha_j+\l _{j}) }
       { ( 2\l _j )_{\alpha_j} \,\big( \thalf(d-j)\big)_{|\alpha^{j+1}|}\,\l _{j} }
$$
with $b_{\alpha}= 2 $ if $\alpha_{d-1}+\alpha_{d}>0$, else $=1.$ Then
$\{Y_\a\}_{| \alpha| = n, \, \alpha_d = 0, 1}$
is an orthonormal basis of $\CH_n^d$. 
\paragraph{Projection opeator}\quad 
The operator $\proj_n\colon L^2(\SS^{d-1}, d\sigma) \to \CH_n^d$ satisfies
\begin{align}
\proj_n P (x) = & \sum_{j=0}^{\lfloor \frac{n}{2} \rfloor}
\frac{1}{ 4^j j!\,\big(- n + 2 - \thalf d\big)_j}\,\|x\|^{2j } \Delta^h P(x) \notag \\
    = & \frac{(-1)^n} {2^n \big( \thalf d -1\big)_n}\,\|x\|^{2 \|\alpha\| + d-2} P\left(\frac{\partial}{\partial x} \right)  \left\{ \|x\|^{-d+2} \right\}.
\end{align}
In particular, the projection of $x^\alpha$ is, up to a constant,
\emph{Maxwell's representation},
\index{Maxwell's representation}%
defined by 
\begin{equation}\label{Maxwell}
 H_\alpha(x):= \|x\|^{2|\alpha| + d-2} \frac{\partial^\alpha} {\partial x^\alpha} \left\{ \|x\|^{-d+2} \right\}, \quad \alpha \in \NN_0^d.
\end{equation}
The set $\{H_\a\}_{| \alpha| = n, \, \alpha_d = 0, 1}$
is a basis of $\CH_n^d$. Furthermore, $H_\alpha$ satisfy a recursive relation 
\begin{equation}
    H_{\alpha+ e_i}(x) = - (2|\alpha| + d-2) x_i H_\alpha(x) + \|x\|^2 \frac{\partial}{\partial x_i} H_\alpha(x). 
\end{equation}
Let $P_n(\cdot, \cdot)$ denote the reproducing kernel of $\CH_n^d$.
Then the projection operator can be written as 
an integral operator
$$
   \proj_n f(x) = \frac{1}{\o_d} \int_{\SS^{d-1}} f(y) P_n(x,y)\,d\sigma(y). 
$$
\paragraph{Reproducing kernel and zonal spherical harmonics}\quad\;
In terms of an orthonormal basis
$\{Y_j\}_{j=1}^{\dim \CH_n^d}$ of $\CH_n^d$, the reproducing kernel,
\index{reproducing kernel}%
by definition, can be written as 
\begin{equation}\label{zonal1}
P_{n}(x,y) =\sum_{1 \le j \le \dim \CH_n^d}Y_j(x) Y_j(y). 
\end{equation}
The kernel is invariant under the action of the orthogonal group
\index{orthogonal group}%
$O(d)$ and it depends only on the 
distance between $x$ and $y$ on the sphere. Moreover,
\begin{equation} \label{zonal2}
  P_n (x,y) = \frac{n + \l} {\l}\,C_n^{\l}\left( \la x,y \ra\right),
\quad x, y \in \SS^{d-1}, \quad \l =\thalf d-1. 
\end{equation}
For $y \in \SS^{d-1}$ fixed, both sides of \eqref{zonal2} as a function of
$x \in \SS^{d-1}$ are \emph{zonal spherical harmonics},
\index{spherical harmonics!zonal}%
i.e., spherical
harmonics which are invariant under an orthogonal transformation leaving
$y$ fixed.
The corresponding homogeneous polynomial
$\|x\|^n C_{n}^{(d-2)/2} \big(\la x, y \ra/\|x\|\big)$ is 
then called a \emph{zonal harmonic polynomial}
\index{zonal harmonic polynomial}%
in $x\in\RR^d$.
The combination of  \eqref{zonal1} and \eqref{zonal2} is known
as the {\it addition formula of spherical harmonics}.
\index{spherical harmonics!addition formula}%
\index{addition formula}%
The reproducing property of the kernel leads to 
the \emph{Funk--Hecke formula}
\index{Funk--Hecke formula}%
\begin{equation}
       \int_{\SS^{d-1}} f(\la x,y \ra) Y(x)\,d\sigma(x)  =  \l_n Y(y), \quad Y \in \CH_n^{d-1}, \quad
           y \in \SS^{d-1}, 
\end{equation}
for all functions $f$ for which the left-hand side is finite, where 
$$
   \l_n = \frac{1}{\o_{d-1}} \int_{-1}^1 f(t)\, \frac{C_n^\l(t)}{ C_n^\l(1)}\,
   (1-t^2)^{\l -\half}\,dt, \quad
    \l =\thalf d-1. 
$$
The \emph{Poisson summation kernel} satisfies, for $ x, y \in \SS^{d-1}$,
\index{Poisson summation kernel}%
\begin{equation}
    \sum_{n=0}^\infty P_n(x,y) r^n =
    \frac{1- r^2} { \big(1- \la x,  y  \ra r + r^2\big)^{\half d}} = 
       \frac{1- r^2} {\|\,r x - y\,\|^d}. 
\end{equation}
\paragraph{Laplace--Beltrami operator}\quad
\index{Laplace--Beltrami operator}%
This is the operator $\Delta_0$ defined by
\begin{equation} \label{Laplace-Beltrami}
    (\Delta_0 f)(x) := (\Delta F)(x), \quad x \in\SS^{d-1},
\end{equation}
where $F(y):=f(y/\|y\|)$ is the extension of $f$ to $\RR^d\backslash\{0\}$
which is homogeneous of degree 0. In terms of the spherical polar 
coordinates $x = r x'$, $r > 0$ and $x' \in \SS^{d-1}$, the usual Laplace operator $\Delta$ is 
decomposed as                                                    
\begin{equation}
   \Delta = \frac{d^2}{d r^2} + \frac{d-1}{r} \frac{d}{dr} + \frac{1}{r^2} \Delta_0. 
\end{equation} 
The spherical harmonics are eigenfunctions of $\Delta_0$\,:
\begin{equation}
   \Delta_0 Y = - n(n+d-2) Y, \qquad Y \in \CH_n^d. 
\end{equation} 
In terms of the spherical coordinates \eqref{eq:sph-polar}, $\Delta_0$ is given by 
\begin{align}
   \Delta_0 = & \frac{1}{(\sin \theta_{d-1})^{d-2} }
     \frac{\partial} {\partial \theta_{d-1}} \left( (\sin \theta_{d-1})^{d-2}  \frac{\partial} {\partial \theta_{d-1}} \right)  \\
      &  +    \sum_{j=1}^{d-2} 
      \frac{1}{(\sin \t_{d-1})^2 \cdots (\sin \theta_{j+1})^2 \sin^{j-1} \theta_j 
      } \frac{\partial} {\partial \theta_j} \left( (\sin\theta_j )^{j-1}  \frac{\partial} {\partial \theta_j} \right).       \notag
\end{align}    
Furthermore, it satisfies a decomposition
\begin{equation} \label{Laplace-Beltrami2}
  \Delta_0 = \sum_{1 \le i < j \le d}  \left(x_i \frac{\partial}{\partial x_j} -  x_j \frac{\partial}{\partial x_i}\right)^2
\end{equation}
for $x \in \SS^{d-1}$. The operator $A_{i,j} = x_i \frac{\partial}{\partial x_j} -  x_j \frac{\partial}{\partial x_i}$ 
is the derivative with respect to the angle (Euler angle) in the polar coordinates of $(x_i,x_j)$-plane, 
and it is also the infinitesimal operator of the regular representation $f \mapsto f(Q^{-1}x)$ of the
rotation group ${\rm SO}(d)$. 
\paragraph{Further results and references}\quad
A number of books contain chapters or sections on spherical harmonics, treating the subject 
from various points of view. For earlier development, especially on $\SS^2$, see \cite{Hobson}. 
A well circulated early introductory is \cite{Muller66}. The connection to Fourier analysis in
Euclidean space is treated in \cite{SW}, see also \cite{Muller}. 
For the connection to the Radon transform, see \cite{Helg}. For applications in integral 
geometry, see \cite{Groemer}. For the point of view of group representations, see \cite{Vile}. 
The fact that the zonal polynomial is of the form $p_n (\la x, y\ra)$ can be used as a starting 
point to study properties of Gegenbauer polynomials, see \cite{Muller66, Vile} as well as 
\cite{AAR}. Spherical harmonics are used as building blocks for orthogonal families on 
radial symmetric measures, see \cite{DunklXu, Xu05} and the next section. 

\subsection{$h$-Harmonics for product weight functions on the sphere}\label{sec:h-harmonics}

A far reaching extension of spherical harmonics
are Dunkl's $h$-harmonics associated with 
reflection groups, see Chapter 7.
We consider the case $\ZZ_2^d$, since explicit 
formulas are available mostly in this case. Let 
\begin{equation} \label{wk-h-harmonic}
   w_\kappa(x) := c_\kappa \prod_{i=1}^d |x_i|^{2 \kappa_i} \quad \hbox{with}\quad c_\kappa 
   :=  \frac{\pi^{\half d}}{\Gamma\big(\thalf d\big)}\,
   \frac {\Gamma\big(|\kappa| + \thalf d\big)} 
 {\Gamma\big(\kappa_1 + \half\big) \ldots
 \Gamma\big(\kappa_d + \half\big)}, 
 \quad |\kappa| := \kappa_1 + \ldots + \kappa_d,
\end{equation}
normalized such that $\int_{\SS^{d-1}} w_\kappa(x)\,d\sigma = 1$. 
This weight function is invariant under the group $\ZZ_2^d$, for which the results for ordinary 
spherical harmonics can be extended in explicit formulas. 
\paragraph{Definition}\quad
\index{h-harmonics@$h$-harmonics}%
The \emph{$h$-harmonics} are homogeneous polynomials that satisfy
$\Delta_h Y =0$, 
where 
\begin{equation}\label{DunklLaplacian}
\Delta_h := \CD_1^2 + \ldots + \CD_d^2
\end{equation}
is the \emph{Dunkl Laplacian}
\index{Dunkl Laplacian}%
and $\CD_j$, $1 \le j \le d$, are 
the \emph{Dunkl operators}
\index{Dunkl operator}%
associated with $\ZZ_2^d$, 
\begin{equation} \label{eq:DunklZ2}
 \CD_j f(x) := \frac{\partial}{\partial x_j } f(x) + \kappa_j 
    \frac{f(x) - f(x_1,\ldots, - x_j, \ldots, x_d)} { x_j}. 
\end{equation}
The \emph{spherical $h$-harmonics}
\index{spherical h@spherical $h$-harmonics}%
are the restriction of
$h$-harmonics to the sphere. 
Let $\CH_n^d(w_\kappa)$ denote the space of $h$-harmonics
of degree $n$. Then
\begin{equation} \label{eq:hBetrami}
  \Delta_{h,0} Y = - n (n+ 2 |\kappa|+ d-2) Y,  \qquad Y \in \CH_n^d(w_\kappa), 
\end{equation}
where $\Delta_{h,0}$ is the spherical $h$-Laplacian operator, and 
$\dim \CH_n^d(w_\kappa) = \dim \CH_n^d$. 
\paragraph{Orthonormal basis}\quad
A basis of $\CH_n^d(w_\kappa)$ can be given in 
spherical coordinates \eqref{eq:sph-polar} and in terms of
\emph{generalized Gegenbauer polynomials}
\index{Gegenbauer polynomials!generalized}%
which are defined by 
\begin{align}
\label{generGegenb}
\begin{split}
C_{2n}^{(\l ,\mu )}(x) &:=\frac{\left( \l +\mu \right) _{n}}
{\big( \mu +\half\big)_{n}}\,
P_{n}^{(\l -1/2,\mu-1/2)}(2x^{2}-1) ,\\
C_{2n+1}^{(\l ,\mu )}(x) &:=\frac{\left( \l +\mu \right)_{n+1}}
{\big( \mu +\half\big)_{n+1}}\,x
P_{n}^{(\l -1/2,\mu+1/2)}(2x^{2}-1).
\end{split}
\end{align}
The polynomials $C_n^{(\l,\mu)}$ are orthogonal with respect to the
weight function
$$
   w_{\l,\mu}(x) = |x|^{2\l} (1-x^2)^{\,\mu-1/2}, \quad x \in [-1,1]. 
$$
Let $h_n^{(\l,\mu)} := c \int_{-1}^1 \big(C_n^{(\l,\mu)}(t)\big)^2
w_{\l,\mu}(t)\,dt$,
normalized such that $h_0^{(\l,\mu)} =1$. Then 
\begin{align}
\begin{split}
     h_{2n}^{(\l,\mu)} &  =
     \frac{\big( \l +\half\big)_{n}(\l+\mu)_{n}}
     {n!\,\big( \mu +\half\big)_{n}}\,
     \frac{\l +\mu}{\l +\mu +2n}, \\
h_{2n+1}^{(\l,\mu)}  & =
\dfrac{\big( \l +\half\big)_{n}\left( \l +\mu \right) _{n+1}}
{n!\,\big( \mu +\half\big)_{n+1}}\,
\frac{\l +\mu}{\l +\mu +2n+1}.
\end{split}
\end{align}
For $d \ge 2$ and $\alpha \in \NN_0^d$, define 
\begin{align} \label{h-harmonic-basis}
Y_{\alpha} (x) := [h_{\alpha}]^{-1} r^{|\alpha|}  g_\alpha(\theta_1) \prod_{j=1}^{d-2}  
   (\sin \theta_{d-j})^{|\alpha^{j+1}|}\,C_{\alpha_j}^{(\l_j,\kappa_j)} (\cos \theta_{d-j}),
\end{align}
where
$g_\alpha(\theta) :=
C_{\alpha_{d-1}}^{(\kappa_d,\kappa_{d-1})}(\cos \theta)$ (if $\alpha_d =0$)
and
$:=\sin\theta\,C_{\alpha_{d-1}-1}^{(\kappa_d+1,\kappa_{d-1})}(\cos\theta)$
if ($\alpha_d=1$),
and where
$|\alpha^j| := \alpha_j + \ldots +\alpha_d$,\;\;
$|\kappa^j| := \kappa_j + \ldots +\kappa_d$,\;\; 
$\l_j := |\alpha^{j+1}|+|\kappa^{j+1}|+\thalf(d-j-1)$,\;\;
and  
$$
\left[h_{\alpha}^n\right]^2 :=
\frac{a_\alpha} {\big(|\kappa| + \thalf d\big)_n}
\prod_{j=1}^{d-1} h_{\alpha_i}^{(\l_i, \kappa_i)} (\kappa_i + \l_i)_{\alpha_i}, \quad
 a_\alpha := \begin{cases} 1 & \hbox{if $\alpha_d =0$} \\
        \kappa_d +\half & \hbox{if $\alpha_d =1$} \end{cases}.
$$  
Then $\{Y_\alpha\}_{|\alpha| = n,\,\alpha_d = 0, 1}$ is an
orthonormal basis of $\CH_n^d(w_\kappa)$. 
\paragraph{Reproducing kernel}\quad
The reproducing kernel of $\CH_n^d(w_\kappa)$
\index{reproducing kernel}%
is given by 
$$
    P_n(w_\kappa; x,y ) = \sum_{1 \le j \le \dim \CH_n^d(w_\kappa)}  Y_j (x)Y_j(y)
$$
where $\{Y_j\}_{j=1}^{\dim \CH_n^d(w_\kappa)}$ is an orthonormal
basis of $\CH_n^d(w_\kappa)$. It 
satisfies a closed formula
\begin{multline} \label{h-harmonic-kernel}
 P_n(w_\kappa;x,y) =  \frac{n + |\kappa| + \half d-1}
     {|\kappa| + \half d-1} \int_{[-1,1]^d} 
  C_n^{|\kappa| +\half d-1}  \big(x_1 y_1 t_1 + \cdots + x_d y_d t_d\big)\\
\times \left({\textstyle\prod\limits_{i=1}^d} c_{\kappa_i}(1+t_i)
    (1-t_i^2)^{\kappa_i-1}\right) dt.
\end{multline}

\paragraph{Further results and references}\quad
The $h$-harmonics associated with a finite reflection group
were first studied by Dunkl in \cite{Dunkl88}. Next he defined his
Dunkl operators in \cite{Dunkl89}.
For an overview of the extensive theory of $h$-harmonics, 
see \cite{DunklXu} and Chapter 7. The 
case $\ZZ_2^d$ was studied in detail in \cite{Xu97b}, which contains \eqref{h-harmonic-basis} and
\eqref{h-harmonic-kernel}, as well as a closed formula for an analog
of the Poisson integral.
A Funk--Hecke type formula was given in \cite{Xu00b}.
For a monic $h$-harmonic basis and a
biorthogonal basis, see \cite{Xu05c}. For a connection to products
of Heine--Stieltjes polynomials, see \cite{Volkmer}. 

\section{Classical orthogonal polynomials of several variables}
General properties of orthogonal polynomials of several variables
were given in \S\ref{sec:General}.
This section contains results for specific weight functions. 

\subsection{Classical orthogonal polynomials on the unit ball} \label{sec:OPball}
\index{orthogonal polynomials of several variables!on classical domains!
on unit ball}%
On the unit ball $B^d := \{x \in \RR^d:\mid\|x\| \le 1\}$, 
consider the weight function
\begin{equation}\label{ball-weight}
  W_\mu(x,y) := \frac {\Gamma\big(\mu + \half(d+1)\big)} 
  {\pi^{d/2}\Gamma\big(\mu + \half\big)}\,
  (1-\|x\|^2)^{\mu-\half}, \quad \mu > -\thalf,
\end{equation} 
normalized such that its integral over $B^d$ is 1.
For $d =2$, see \S\ref{sec:OP-disc}.
\paragraph{Differential operator}\quad
Orthogonal polynomials of degree $n$ with respect to $W_\mu$
are eigenfunctions of  a second order differential operator:
\begin{equation} \label{eq:Bdiff}
 \left( \Delta -{\textstyle\sum\limits_{j=1}^d}\,\frac{\partial}{\partial x_j}\,x_j
 \left((2 \mu -1)  + {\textstyle\sum\limits_{i=1}^d} x_i
 \frac{\partial }{\partial x_i} \right)\right)P   =
  -(n+d) (n + 2 \mu -1)P, \quad P \in \CV_n^d, 
\end{equation} 
where $\Delta$ is the Laplace operator.  
\paragraph{First orthonormal basis}\quad
Associated with $x = (x_1,\ldots,x_d) \in
\RR^d$, define $\xb_j := (x_1, \ldots, x_j)$ for $1 \le j \le d$ and $\xb_0 := 0$.
For  $\a \in \NN_0^d$ and $1 \le j \le d$, let $\a^j : = (\a_j ,  \ldots, \a_d)$ and
$\a^{d+1} : =0$.
An orthonormal basis $\{P_\a\}_{|\a| =n}$ of $\CV_n^d$ is given by
\begin{equation}\label{ball-base1}
  P_\alpha(x) : = [h_\a]^{-1} \prod_{j=1}^d 
  \big(1-\|\xb_{j-1}\|^2\big)^{\half\alpha_j}\,
      C_{\alpha_j}^{\l_j} \biggl(\frac {x_j} {(1-\|\xb_{j-1}\|^2)^{1/2}}\biggl),
\end{equation} 
where $\l_j := \mu+|\alpha^{j+1}|+ \half(d-j)$ and $h_\alpha$ is given by
\begin{equation}\label{ball-base1-norm}
[h_{\alpha}]^2 :=
\frac{\big(\mu + \half d\big)_{|\a|}} {\big(\mu + \half(d+1)\big)_{|\a|}}
\prod_{j=1}^d  \frac{\big(\mu + \half(d-j)\big)_{|\a^j|}\,
\big(2 \mu + 2 |\alpha^{j+1}|+ d-j\big)_{\a_j}} 
  {\big(\mu + \half(d-j+1)\big)_{|\a^j|}\,\a_j!}, 
\end{equation} 
and where the case $\mu=0$ can be obtained as a limit for $\mu\to0$
by using \eqref{limGegCheb}, similarly as for
\eqref{basis1_ball}, \eqref{basis1_ball_norm}.
\paragraph{Second orthonormal basis}\quad\,
Let $r_{k}^d := \dim \CH_{k}^d$.
For $0 \le j \le n/2$ let $\big\{Y_{\ell, n-2j}\mid 1\le\ell\le r_{n-2j}^d\big\}$
be an orthonormal 
basis of $\CH_{n-2j}^d$, the space of spherical harmonics of degree $n-2j$,
with respect to 
the normalized surface measure. For $0\le j \le n/2$, define
\begin{equation}\label{ball-base2}
P_{\ell,j}(x) := [h_{j,n}]^{-1}
     P_{j}^{(\mu-\half, n-2j + \half d-1)}\big(2\|x\|^2-1\big)\,Y_{\ell,n-2j}(x),
\end{equation}
where 
$$
[h_{j,n}]^2 := \frac{\big(\mu +\half\big)_j\,\big(\half d\big)_{n-j}}
{ j!\,\big(\mu+\half(d+1)\big)_{n-j}}\,
\frac{\big(n-j+\mu+ \half(d-1)\big)}{\big(n+\mu+ \half(d-1)\big)}.
$$
Then $\big\{P_{\ell,j}\mid 1 \le \ell \le r_{n-2j}^d,\;0\le j \le n/2\big\}$
is an orthonormal basis of $\CV_n^d$.
\paragraph{Appell's biorthogonal polynomials}\quad
\index{Appell biorthogonal polynomials!on unit ball}%
These are two families of polynomials yielding bases 
$\big\{U_\a\}_{|\a|=n}$ and
$\big\{V_\a\}_{|\a|=n}$ of $\CV_n^d$, where the second
basis is the monic basis, up to constant factors, and the first basis is
biorthogonal to it.
\vskip\medskipamount
\begin{enumerate}[label=(\roman*)]
\item
The family $\{U_\a\}$ is defined by the generating function
\index{generating function}%
\begin{equation} \label{eq:Ugenerating}
\Big(\big(1-\la b, x \ra\big)^2 + \|b\|^2\big(1-\|x\|^2\big)\Big)^{-\mu}
   = \sum_{\alpha \in \NN_0^d} b^\alpha U_\alpha(x),\qquad
   b \in \RR^d,\; \|b\| < 1.
\end{equation}
It satisfies the Rodrigues type formula,
\index{Rodrigues type formula}%
\begin{equation}\label{Rodrigues-ball}
   U_\alpha(x) = \frac {(-1)^{|\alpha|}\,(2\mu)_{|\alpha|}}
   {2^{|\alpha|}\,\big(\mu+\half\big)_{|\alpha|}\,\alpha !}\,
    \big(1-\|x\|^2\big)^{-\mu+\half}\,
 \frac{\partial^{|\alpha|}} {\partial x^\alpha}\,
 \big(1-\|x\|^2\big)^{|\alpha|+\mu-\half}.
\end{equation} 
where $\frac{\partial^{|\a|}}{\partial x^\a} :=  \frac{\partial^{|\alpha|}} {\partial x_1^{\alpha_1} 
\ldots \partial x_d^{\a_d}}$. Furthermore, it can be explicitly given as
\begin{align}
U_\alpha(x) & = \frac {(2\mu)_{|\alpha|}} {\alpha!} 
 \sum_{\beta \le \alpha}  \frac{(-1)^{|\beta|}\,(-\alpha)_{2\beta}}
 {2^{2 |\beta|}\,\beta!\,\big(\mu+\half\big)_{|\beta|}}\,
 x^{\alpha - 2 \beta}\big(1-\|x\|^2\big)^{|\beta|} \\
& = \frac{(2 \mu)_{|\alpha|} x^\alpha}{\alpha!}\,
F_B\Big(-\thalf\alpha,  -\thalf(\alpha-\mathbf{1}); \mu+\thalf;
  x_1^{-2}\big(1-\|x\|^2\big),\ldots,x_d^{-2}\big(1-\|x\|^2\big)\Big), \notag 
\end{align}
where $(-\a)_{2\b} := (-\alpha_1)_{2\beta_1} \ldots  (-\alpha_d)_{2\beta_d}$, 
$\mathbf{1} := (1,\ldots,1)$, $\beta \le \alpha$ means  $\beta_i \le \alpha_i$
($1\le i \le d$), 
and $F_B$ is Lauricella's hypergeometric series 
of type $B$ (see \cite{Ex} and Chapter 3).
\index{Lauricella hypergeometric function}%
\index{hypergeometric function!Lauricella}%
\vskip\medskipamount
\item
The family $\{V_\a\}$ is defined by the  generating function
\index{generating function}%
\begin{equation} \label{eq:Vgenerating}
  \big(1- 2 \la b, x \ra + \|b\|^2\big)^{-\mu - (d-1)/2} = 
    \sum_{\alpha \in \NN_0^d} b^\alpha V_\alpha(x),\qquad
    b \in \RR^d,\;\|b\| < 1.
    \end{equation} 
The generating function implies that 
\begin{equation}\label{eq:Vpoly-Ball1}  
     \|b\|^n\,C_n^{\mu+\half(d-1)} \left(\frac{ \la b, x \ra}{\|b\|}\right)
          = \sum_{|\a| = n } b^\a V_\a(x). 
\end{equation}
The polynomial $V_\a$ can be written explicitly as 
\begin{align} \label{eq:Vpoly-Ball2}
 V_\alpha(x) & = 2^{|\alpha|}x^\alpha \sum_{\gamma < \alpha}
   \frac{\big(\mu+\thalf(d-1)\big)_{|\alpha| - |\gamma|}\,
   (-\alpha+\gamma)_{\gamma}}
   {(\alpha-\gamma)!\,\gamma!}\,2^{-2 |\gamma|} x^{- 2\gamma} \\
    & = \frac{2^{|\alpha|}\big(\mu+\half(d-1)\big)_{|\alpha|}} {\alpha!}\,x^\alpha 
F_B\Big(-\thalf\alpha,   -\thalf(\alpha-\mathbf{1}); 
  -|\alpha|- \mu - \thalf(d-3); x_1^{-2},\ldots,x_d^{-2}\Big). \notag
\end{align}
\end{enumerate}

Neither $\{U_\a\}_{\a| =n}$ nor $\{V_\a\}_{|\a| =n}$ is an orthogonal
basis of $\CV_n^d$, but the two bases are biorthogonal to each other:
\begin{equation} \label{eq:biorth-ball}
   \int_{B^d} V_\alpha(x) U_\beta(x) W_\mu(x)\,dx = 
\frac{\mu+\thalf(d-1)}{|\alpha|+\mu+\thalf(d-1)}\,
\frac{(2\mu)_{|\alpha|}} {\alpha!}\,\delta_{\alpha,\beta}. 
\end{equation} 
{\it The monic basis}\quad For each $\a$, define 
\begin{equation}  \label{eq:monic-ball}
R_\a(x) :=  \frac{\alpha!} {2^{|\alpha|}\big(\mu+\thalf(d-1)\big)_{|\alpha|}}
\,V_\a(x). 
\end{equation}
Then the polynomials $R_\a$ ($|\a|=n$)
form a monic basis of $\CV_n^d$, i.e.,
 $R_\a(x) = x^\a - Q_\a(x)$ with $Q_\a \in \Pi_{n-1}^d$. 
The $L^2(W_\mu, B^d)$ norm, $\|\cdot\|_{2,\mu}$, of $R_\a$ is the error of the best approximation 
of $x^\a$, which satisfies a closed formula 
\begin{equation}  \label{eq:L2norm-ball}
   \min_{P \in \Pi_{n-1}^d} \|x^\a - P(x) \|_{2,\mu}^2  = \|R_\a\|_{2,\mu}^2 =
     \frac{\l \, \alpha!} {2^{n-1}(\l)_{n}} \int_0^1
     {\textstyle\prod\limits_{i=1}^d}  P_{\alpha_i}(t)\,t^{n+2\l -1}\,dt,     
\end{equation}
where $\l = \mu + (d-1)/2 >0$, $n = |\alpha|$, and $P_{\a_i}$ is the
Legendre polynomial. 
\\[\medskipamount]
{\it Reproducing kernel}\quad
\index{reproducing kernel}%
The kernel $P_n(\cdot, \cdot)$  of $\CV_n^d$
with respect to $W_\mu$, defined in \eqref{eq:reprod},
satisfies a compact formula. For $\mu > 0$, $x, y \in B^d$, 
\begin{align} \label{reprod-ball1}
P_n(x,y)= c_{\mu} \frac{2n+2\mu+d-1}{2\mu+d-1}
    \int_{-1}^1 C_n^{\mu+\half(d-1)}\left(\la x,y
    \ra + t {\textstyle\sqrt{1-\|x\|^2} \sqrt{1-\|y\|^2}}\,\right)   (1-t^2)^{\mu-1} \,dt, 
\end{align}
where $[c_\mu]^{-1} := \int_{-1}^1 (1-t^2)^{\mu-1}  d t$. 
For $\mu = 0$, $x,y \in B^d$ this degenerates to
\begin{equation} \label{reprod-ball2}
 P_n(x,y) =  \frac{n + \half(d-1)}{d-1}  \sum_{\epsilon=0,1}
    C_n^{ \half(d-1)} \left(\la x,y\ra +
    (-1)^\epsilon {\textstyle\sqrt{1-\|x\|^2} \sqrt{1-\|y\|^2}}\,\right).
 \end{equation}
These formulas are essential for obtaining sharp results for convergence of orthogonal expansions.
When $\mu = \half (m-1)$, \eqref{reprod-ball1} can also be written as
\begin{align} \label{reprod-ball2A}
P_n(x,y) =\frac{2n+m+d-2}{m+d-2} \int_{\SS^{m-1}}
C_n^{\half(m+d)-1} \left (  \la x, y \ra 
            +{\textstyle\sqrt{1-\|x\|^2}\sqrt{1-\|y\|^2}} \, \la \xi,e_1 \ra \right ) d\s_m(\xi), 
\end{align}
where $d\sigma_m$ is the normalized surface measure on $\SS^{m-1}$.

For the constant weight $W_{1/2}(x)=(d+1)/\pi$, there is another formula for the reproducing kernel,
\begin{align} \label{reprod-ball3}
P_n(x,y)= (2nd^{-1}+1)
      \int_{\SS^{d-1}} C_n^{d/2} \big(\la x, \xi \ra\big)\,
      C_n^{d/2} \big(\la  \xi, y \ra\big)\,d\sigma_d(\xi).
\end{align}
\paragraph{Further results and references}\quad
For orthogonal bases on the ball, see \cite{ApKa, DunklXu, Er2}. There are further
results on biorthogonal bases, see \cite{ApKa, Er2}. For
the monic basis, see \cite{Xu05c}. Orthogonal bases 
consisting of ridge polynomials
were discussed in \cite{Xu00b},
together with a Funk--Hecke type formula for orthogonal polynomials.
The compact formulas \eqref{reprod-ball1} and  \eqref{reprod-ball2} for the reproducing kernels were proved 
in \cite{Xu99a} and used to study expansion problems, whereas the compact formula \eqref{reprod-ball2A}
was proved in \cite[Theorem 2.6]{Xu01a} (there take $H_2(\eta)=1$).
Formula  \eqref{reprod-ball3} was proved in \cite{Petru} 
in the context of approximation by ridge functions, and in \cite{Xu07} in connection with Radon 
transforms. 
In \cite{ACH} three-term relations are used to develop an efficient
numerical algorithm for the evaluation of orthogonal 
polynomials in \eqref{ball-base1} with $d =2,3$. 
For convergence and summability of orthogonal expansions,
see \S\ref{OrthoExpan}.
  
\subsection{Classical orthogonal polynomials on the simplex}\label{sec:OPsimplex}
\index{orthogonal polynomials of several variables!on classical domains!
on simplex}%
For $x \in \RR^d$, let $|x|: =  x_1+\cdots +x_d$.
Let $T^d:=\{x \in \RR^d\mid x_1,\ldots, x_d, 1- |x| \ge0\}$
be the simplex in $\RR^d$.
The classical weight function on $T^d$ is defined by
\begin{equation} \label{eq:weight-Td}
 W_{\kappa}(x) : = \frac{\Gamma\big(|\kappa|+\thalf(d+1)\big)}{\prod_{i=1}^{d+1}\Gamma\big(\kappa_i+\thalf\big)}\,
         x_1^{\kappa_1-\half} \ldots x_d^{\kappa_d-\half}
         (1-|x|)^{\kappa_{d+1} -\half},\quad
\kappa_1,\ldots, \kappa_{d+1} > -\thalf.
\end{equation}

\paragraph{Differential operator}\quad
Orthogonal polynomials of degree $n$ with respect to $W_\kappa$ are eigenfunctions of 
a second order differential operator, 
\begin{multline} \label{eq:Tdiff}
\left(\sum_{i=1}^d x_i(1-x_i)\,\frac{\partial^2} {\partial x_i^2}- 
2\sum_{1\le i<j\le d}x_ix_j\,\frac{\partial^2}{\partial x_i\partial x_j} 
+\sum_{i=1}^d \left( \big(\kappa_i +\thalf\big) -
\big(|\kappa|+\thalf(d+1)\big)x_i\right)
\frac{\partial}{\partial x_i}\right)P\\ 
= -n  \left(n+|\kappa| + \thalf(d-1) \right) P, \qquad 
P \in \CV_n^d,
\end{multline} 
where $|\kappa| = \kappa_1 + \cdots +\kappa_{d+1}$.  
\paragraph{An orthonormal basis}\quad
To state this basis we use the notation of $\xb_j$ and $\alpha^j$ as in 
the first orthonormal basis on $B^d$ of \S\ref{sec:OPball}.
We also put  $\kappa^j := (\kappa_j,\ldots, \kappa_{d+1})$
($j = 0,1,\ldots, d+1$) if $\kappa =(\kappa_1, \ldots, \kappa_{d+1})$.
Then an orthonormal basis
$\{P_\a\}_{|\a| =n}$  of $\CV_n^d$ is given by
\begin{equation} \label{eq:OPonT} 
P_\alpha(x) := [h_{\alpha}]^{-1} \prod_{j=1}^d \big(1-|\xb_{j-1}|\big)^{\alpha_j}   
\,P_{\alpha_j}^{(a_j,\kappa_j -\half)}\left(\frac{2x_j}{1-|\xb_{j-1}|} -1\right), 
\end{equation} 
where $a_j:=2|\alpha^{j+1}|+|\kappa^{j+1}|+\thalf(d-j-1)$ and  
$h_\alpha$ is given by
\begin{equation*}
[h_{\alpha}]^2 :=  \prod_{j=1}^d \frac{\big(\kappa_j +\half\big)_{\a_j}\,
\big(|\kappa^{j+1}|+\half(d-j+1)\big)_{|\alpha^j|+|\alpha^{j+1}|}}
{\alpha_j! \,\big(|\kappa^j|+\half(d-j+2)\big)_{|\alpha^j|+|\alpha^{j+1}|}}\,
\frac{2(a_j+\kappa_j+\a_j) +1}{2(a_j+\kappa_j+2 \a_j)+1}.
 \end{equation*} 

\paragraph{Appell's biorthogonal polynomials}\quad
\index{Appell biorthogonal polynomials!on simplex}%
These are two families of polynomials yielding bases 
$\big\{U_\a\}_{|\a|=n}$ and
$\big\{V_\a\}_{|\a|=n}$ of $\CV_n^d$, where the second
basis is the monic basis and the first basis is
biorthogonal to it.\\[\medskipamount]
(i)  The family $\{U_\a\}$ is defined by the Rodrigues type formula
\index{Rodrigues type formula}%
\begin{equation} \label{eq:UonT}
U_\alpha(x) :=  x_1^{-\kappa_1+ \half} \ldots x_d^{-\kappa_d+\half}
    \big(1-|x|\,\big)^{-\kappa_{d+1} +\half}\,
\frac {\partial^{|\alpha|}} {\partial x^\alpha}
\left(x_1^{\alpha_1+\kappa_1-\half} \ldots 
     x_d^{\alpha_d+\kappa_d-\half} (1-|x|)^{|\alpha|+\kappa_{d+1}-\half}\right). 
\end{equation}
(ii) The family $\{V_\a\}$ is explicitly defined by
\begin{align} \label{eq:VonT}
  V_\alpha(x) & := \sum_{\beta \le \alpha} (-1)^{n+|\beta|}
\left({\textstyle \prod\limits_{i=1}^d}\,\binom{\alpha_i}{\beta_i}\,
\frac{\big(\kappa_i+ \half\big)_{\alpha_i}}
{\big(\kappa_i+\half\big)_{\beta_i}}\right)\,
\frac{\big(|\kappa|+\half(d-1)\big)_{n+|\beta|}}
   {\big(|\kappa|+\half(d-1)\big)_{n+|\alpha|}}\, x^\beta  \notag \\ 
  & =   \frac{(-1)^n \big(\kappa+\frac{\bf 1}2\big)_\alpha}
  {\big(n+|\kappa|+\half(d-1)\big)_{|\alpha|}}\,
        F_A\left (n+|\kappa|+\thalf(d-1), - \alpha; \kappa+ \tfrac{\bf 1}{2}, x\right) 
\end{align}
where $F_A$ denotes Lauricella's hypergeometric series of type A
(see \cite{Ex} and Chapter 3).
\index{Lauricella hypergeometric function}%
\index{hypergeometric function!Lauricella}%

Neither $\{U_\a\}_{\a| =n}$ nor $\{V_\a\}_{|\a| =n}$ is an orthogonal
basis of $\CV_n^d$, but the two bases are biorthogonal to each other:
\begin{equation}\label{biorthosimplex}
  \int_{T^d} V_\beta(x) U_\alpha(x) W_\kappa(x)\,dx 
= \frac {\big(\kappa +\half\big)_{\alpha}
\big(\kappa_{d+1}+\half\big)_{|\alpha|}}
 {\big(|\kappa|+\half(d+1)\big)_{2|\alpha|}}\, \alpha!\,\delta_{\alpha,\beta}.
\end{equation}

\paragraph{Monic orthogonal basis}\quad
By \eqref{eq:VonT} and \eqref{biorthosimplex} the polynomials
$V_\a$ ($|\a|=n$)
form a monic basis: $V_\a(x) = x^\a - Q_\a(x)$, $Q_\a \in \Pi_{n-1}^d$. 
Such polynomials can be defined more generally in view of the observation
that the simplex $T^d$ is associated with the permutation 
group of $X : = (x_1,\ldots,x_d, x_{d+1})$, where $x_{d+1} = 1-|x|$.
For $x \in T^d$ and  $\a \in \NN_0^{d+1}$, define
$X^\a := x_1^{\a_1} \ldots x_d^{\a_d} (1-|x|)^{\a_{d+1}}$. 
For $|\a|=n$ let
$R_\alpha(x) = X^\alpha - Q_\alpha(x)$ be the element of $\CV_n^d$
such that $Q_\alpha$ is a polynomial of degree at 
most $|\a| - 1$. When
$\a_{d+1} =0$, $R_\alpha(x)$ agrees with $V_\alpha(x)$ in
\eqref{eq:VonT} up to a constant factor.
The $L^2(W_\kappa, T^d)$ norm $\|\cdot\|_{2,\kappa}$ of $R_\alpha$ is 
the error of the best approximation of $X^\a$, which satisfies a closed
formula: 
\begin{equation} \label{monic-Td-norm}
\min_{P \in \Pi_{n-1}^d} \|X^\a - P(x)\|_{2,\kappa}^2   =
\|V_\a\|_{2,\kappa}^2
=    \frac{\big( |\kappa| + \half(d-1)\big)
\big(\kappa +\frac{\mathbf{1}}2\big)_{\alpha}} 
{\big( |\kappa| + \half(d-1\big )_{2 |\alpha|}}  
     \int_0^1 \left({\textstyle\prod\limits_{i=1}^{d+1}}
        P_{\alpha_i}^{(0, \kappa_i - \half)}(2 t-1)\right)
        t^{|\alpha|+  |\kappa| + \frac{d-3}{2}}\,dt.
\end{equation}
If $\a_{d+1} = 0$, \eqref{monic-Td-norm} gives the error of best
approximation to $x^\a$. 
\paragraph{Reproducing kernel}\quad
\index{reproducing kernel}%
The kernel $P_n(\cdot, \cdot)$ of $\CV_n^d$
with respect to $W_\kappa$, defined in \eqref{eq:reprod},
satisfies a compact formula.
For $\kappa_i > 0$, $1 \le i \le d$, $x, y \in T^d$, 
\begin{multline} \label{reprod-simplex1}
 P_n(x,y) = c_\k\, 
 \frac{2(2n+|\k|)+d-1}{2|\k|+d-1}\,
 \frac{\big(|\k|+\half(d-1)\big)_n}{\big(\half\big)_n}\\
 \times\int_{[-1,1]^{d+1}}
 P_n^{(|\k|+(d-2)/2,-1/2)} \left( 2z(x,y,t)^2-1\right)
 \left({\textstyle \prod\limits_{i=1}^{d+1}} (1-t_i^2)^{\k_i-1}\right) dt, 
\end{multline}
where
$z(x,y,t):=\sqrt{x_1 y_1}\, t_1 + \cdots + \sqrt{x_d y_d}\, t_d+
\sqrt{1-|x|} \sqrt{1-|y|}\, t_{d+1}$ and
\[
[c_\kappa]^{-1} :=
\int_{[-1,1]^{d+1}}
 \left({\textstyle \prod\limits_{i=1}^{d+1}} (1-t_i^2)^{\k_i-1}\right) dt.
\]
If some
$\kappa_i =0$, then the formula  holds under the limit relation 
\begin{equation*}
 \lim_{\l \to 0} c_\l \int_{-1}^1 g(t) (1-t)^{\l -1}\,dt   =
 \thalf\big(g(1) + g(-1)\big).
\end{equation*}

\paragraph{Further results and references}\quad
For $d=2$ the polynomials $U_\alpha$ on $T^d$ were defined in
\cite{ApKa} and the polynomials $V_\alpha$ were studied in \cite{FaLi}.
For $d>2$ they appeared in \cite{GrMo} when $W_\kappa(x) =1$,
and in \cite{DunklXu}  in general.
The monic basis of polynomials $R_\a$ was studied in \cite{Xu05c}. 
The formula \eqref{reprod-simplex1} of the reproducing kernel appeared in
\cite{Xu98c}. 
A product formula for orthogonal polynomials on the simplex was
established in \cite{KoornSchw}. 
The polynomials in \eqref{eq:OPonT} serve as generating functions
for the Hahn polynomials 
in several variables \cite{KaMc, Xu15b}. For convergence and summability
of orthogonal expansions see \S\ref{OrthoExpan}.
   
\subsection{Hermite polynomials of several variables} \label{sec:ProdHermite}
\index{orthogonal polynomials of several variables!on classical domains!
Hermite}
These are orthogonal polynomials with respect to the product weight
\begin{equation} \label{Hermite-weight}
    W_H(x) := \pi^{-d/2} e^{- \|x\|^2}, \qquad x \in \RR^d.
\end{equation}
Many properties are inherited from Hermite polynomials of one variable. 
\paragraph{Differential operator}\quad
Orthogonal polynomials of degree $n$ with respect to
$W_H$ are eigenfunctions of a second order differential operator:
\begin{equation} \label{Hermite-pde}
\left(\Delta  - 2 \sum_{i=1}^d x_i\,\frac{\partial }{\partial x_i}\right)
P = - 2n P, \qquad P \in \CV_n^d. 
\end{equation}

\paragraph{Product orthogonal basis}\quad
For $\alpha \in \NN_0^d$, define 
\begin{equation} \label{Hermite-product}
 H_\alpha(x) := H_{\a_1}(x_1) \ldots H_{\a_d}(x_d). 
\end{equation}
Then $\big\{ [c_\alpha]^{-1} H_\a\big\}_{|\a| =n}$ is an orthonormal basis of
$\CV_n^d$, where
$[c_\a]^2 := 2^{|\alpha|} \alpha!$.  As products of Hermite polynomials of one variable, 
they inherit a generating function and Rodrigues type formula. Furthermore, they satisfy
\begin{align*}
\sum_{|\alpha|=n}  \frac{H_\alpha(y)}{\alpha!}\,x^\alpha
  =  \frac{1}{n!\,\pi^\half} \int_{-\infty}^\infty 
  H_n\Big(\la x,y\ra+ s \sqrt{1-\|x\|^2}\,\Big)\,e^{-s^2}\,ds,\quad
 x \in \BB^d,\;y \in \RR^d.
\end{align*}
In particular, 
\begin{equation}
\sum_{|\alpha|=n}  \frac{H_\alpha(y)}{\alpha!}\,x^\alpha
  =  \frac{1}{n!}\,H_n\big(\la x,y\ra)\big), \qquad \|x\|=1,\;y \in \RR^d. 
\end{equation}
 
\paragraph{Second orthonormal basis}\quad
For $0 \le j \le n/2$ and $r_k^d := \dim \CH_k^d$
let $\{Y_{\ell, n-2j}\}_{\ell=1}^{r_{n-2j}^d}$
be an orthonormal basis of
$\CH_{n-2j}^d$ as in \eqref{ball-base2}.
Define 
\begin{equation}  \label{Hermite-2nd-basies} 
   P_{\ell,j}^n(x) := [c_{j,n}]^{-1} 
           L_j^{n-2j + \half(d-2)}\big(\|x\|^2\big)\,Y_{\ell,n-2j}(x), \qquad [c_{j,n}]^2 := \frac{\big(\half d\big)_{n-j} } {j!}.
\end{equation}
Then $\big\{P_{\ell,j}^n\mid 1 \le \ell \le r_{n-2j}^d,\;
0 \le j \le \half n\big\}$ is an orthonormal basis of $\CV_n^d$.
\paragraph{Mehler formula}\quad
\index{Mehler formula}%
The reproducing kernel $P_n(\cdot,\cdot)$ of $\CV_n^d$, as defined in 
\eqref{eq:reprod}, satisfies, for $0<z<1$ and $x,y \in \RR^d$, 
\begin{align}
\sum_{n=0}^\infty  P_n(x,y)\,z^n = 
\frac{1}{(1-z^2)^{d/2} }\,
   \exp \left( -\frac{z^2 (\|x\|^2+\|y\|^2) -2z \la x, y\ra}{1-z^2} \right)  .  
\end{align}
 
\paragraph{Further results and references}\quad
The study of Hermite polynomials of several variables
was started by Hermite and followed by many other authors, see \cite{ApKa, Er2} for references. 
Analogues of Hermite polynomials can be defined more generally for the
weight function 
\begin{equation} \label{Hermite-general}
    W(x) =  (\det A)^{\half}\,\pi^{- d/2}\,
    \exp \left(- x^\tran A x \right),
\end{equation}
where $A$ is a positive definite matrix. Two families of biorthogonal polynomials can be defined for $W$
in \eqref{Hermite-general}, which coincide when $A$ is an identity matrix. These were studied in \cite{ApKa}, 
see also \cite{Er2}. Since $A$ is positive definite, it can be written as $A = B^\tran B$. Thus, orthogonal 
polynomials for $W$ in \eqref{Hermite-general} can be derived from Hermite polynomials
for $W_H$ by a change of variables. 

\subsection{Laguere and generalized Hermite polynomials} \label{sec:ProdLaguerre} 
\paragraph{Laguerre polynomials of several variables}\quad
\index{orthogonal polynomials of several variables!on classical domains!
Laguerre}%
Put $|x| := x_1 + \cdots x_d$ and
$\RR_+^d := \big\{x \in \RR^d\mid x_1,\ldots,x_d \ge 0\big\}$.
\emph{Laguerre polynomials of several variables}
are orthogonal polynomials with respect to 
the weight function
\begin{equation} \label{Laguerre-weight}
    W_L(x) :=  \frac{1}{\prod_{k=1}^d \Gamma(\kappa_i + \thalf)}\,
    x^\kappa e^{- |x|}, \qquad x \in \RR_+^d,
\end{equation}
which can be written as a product of weight functions in one variable.
Many properties are inherited from Laguerre polynomials of one variable. 
 
The polynomials in $\CV_n$ are eigenfunctions of a differential
operator:
\begin{equation} \label{Laguere-pde}
 \left(\sum_{i=1}^d x_i\,\frac{\partial^2 }{\partial x_i^2} + 
    \sum_{i=1}^d (\kappa_i+1-x_i)\,\frac{\partial }{\partial x_i}\right)P 
   = - n P, \qquad P \in \CV_n^d. 
\end{equation}
An orthonormal basis of $\CV_n^d$ is given by
$\{ L_\a^\kappa\}_{|\a| =n}$, where
\begin{equation} \label{Laguerre-product}
\textstyle L^\kappa_\alpha(x) := \binom{\alpha+\kappa}{\alpha}^{-1/2}
 L_{\a_1}^{\kappa_1}(x_1) \ldots L_{\a_d}^{\kappa_d}(x_d)
\end{equation}
The reproducing kernel $P_n(\cdot,\cdot)$ of $\CV_n^d$, defined in
\eqref{eq:reprod}, satisfies, for $0<z<1$ and $x,y \in \RR_+^d$, 
\begin{align}
\sum_{n=0}^\infty  P_n(x,y)\,z^n =(1-z)^{-1}
    \prod_{i=1}^d  \exp\left(-\frac{z (x_i+y_i)}{1-z}\right)\,
    (x_iy_iz)^{- \half\kappa_i}\,
    I_{\kappa_i}\left(\frac{2\sqrt{x_iy_iz}}{1-z}\right),  
\end{align}
where $I_\kappa$ denote the modified Bessel function of order $\kappa$. 
\paragraph{Generalized Hermite polynomials of several variables}\quad
\index{orthogonal polynomials of several variables!on classical domains!
generalized Hermite}%
These are orthogonal polynomials with respect to the product weight function
\begin{equation}
       W_\kappa(x) :=   \prod_{i=1}^d 
       |x_i|^{2 \kappa_i} e^{-x_i^2}
       \end{equation}
Let $\Delta_h$ be the Dunkl Laplacian \eqref{DunklLaplacian}.
The orthogonal 
polynomials in $\CV_n^d$ are eigenfunctions of a
differential-difference operator:
\begin{equation}
   \Big(\Delta_h  - 2 \la x,\nabla \ra\Big) P = - 2 n P, \qquad
   P \in \CV_n^d.
\end{equation}
For $\kappa > 0$ let the
\emph{generalized Hermite polynomials} $H_n^\kappa$
in one variable be defined by
\index{orthogonal polynomials of one variable!generalized Hermite}%
\begin{align}
\begin{split}
H_{2n}^\kappa(x) & := (-1)^n 2^{2n} n!\,L_n^{\kappa -\half}(x^2), \\
H_{2n+1}^\kappa(x) & := (-1)^n 2^{2n+1} n!\,x\,L_n^{\kappa +\half}(x^2). 
\end{split}
\end{align}
Let $h_n^\kappa = \Gamma\big(\kappa+\half\big)^{-1} \int_\RR [H_n^\kappa(x)]^2 x^\kappa e^{-x^2}\,dx$,
normalized such that $h_0^\kappa =1$. Then
\begin{align}
  h_{2n}^\kappa = 2^{4n} n!\,(\mu+\thalf)_n \quad \hbox{and} \quad  h_{2n+1}^\kappa = 2^{4n+2} n!\,(\mu+\thalf)_{n+1}.
\end{align}
An orthonormal basis of $\CV_n^d$ for $W_\kappa$ is given by
$\{P_\alpha\}_{|\alpha| =n}$, where 
\begin{equation} \label{GeneralHermitebasis}
P_\alpha(x) :=
\big(h_{\alpha_1}^{\kappa_1} \ldots
h_{\alpha_d}^{\kappa_d}\big)^{-\half}
H_{\alpha_1}^{\kappa_1}(x_1) \cdots H_{\alpha_d}^{\kappa_d}(x_d). 
\end{equation}
Another orthonormal basis, analogous to \eqref{Hermite-2nd-basies},
can be given, in polar coordinates, 
by Laguerre polynomials and $h$-harmonics.  
\paragraph{Further results and references}\quad 
As products of Laguerre polynomials of one variable,
the polynomials $L^\kappa_\alpha(x)$ in \eqref{Laguerre-product} have a generating function, a 
Rodrigues type formula, as well as a product formula that induces a convolution structure. 
The generalized Hermite polynomials can be defined for weight functions under other 
reflection groups.
For those and further properties of these functions,
including a Mehler type formula, see \cite{DunklXu, Rosl98, Xu01b}. 
 
\subsection{Jacobi polynomials of several variables}\label{sec:ProdJacobi}
\index{orthogonal polynomials of several variables!on classical domains!
product Jacobi}%
These polynomials are orthogonal with respect to the weight function
\begin{equation}\label{Jacopbi-weight}
   W_{a,b}(x) := \prod_{i=1}^d (1-x_i)^{a_i}(1+x_i)^{b_i}, \quad x \in [-1,1]^d.
\end{equation}
An orthogonal basis is formed by the \emph{product Jacobi polynomials}, 
$$
   P_\alpha(x) :=  P_{\a_1}^{(a_1,b_1)}(x_1) \ldots P_{\a_d}^{(a_d,b_d)}(x_d), \quad 
\a \in \NN_0^d.    
$$
Most results for these orthogonal polynomials follow from
properties of Jacobi polynomials of 
one variable. The reproducing kernel $P_n(\cdot, \cdot)$ of $\CV_n^d$ satisfies 
\begin{equation}  
  \sum_{n=0}^\infty P_n(x, {\bf 1}) r^n  = \prod_{i=1}^d  \frac{(1-r)(1+r)^{a_i-b_i+1}}{(1-2r x_i + r^2)^{a_i + 3/2}} \; 
  \hyp21{\half(b_i - a_i), \half(b_i-a_i -1)}
   {b_i+1} {\frac{2r (1+x_i)}{(1+r)^2}}, 
\end{equation}
where $0 < r< 1$ and $\mathbf{1}= (1,\ldots, 1)$. 

In the case of the product Chebyshev weight function, i.e.,
$a_i = b_i = -\half$, $i =1,\ldots,d$, 
the reproducing kernel
\index{reproducing kernel}%
$P_n(\cdot,\cdot)$ satisfies a closed formula
\cite[Theorem 9.6.3]{DunklXu}. This is
given in terms of
a \emph{divided difference} $[x_1,\ldots,x_d] f$ defined by
\index{divided difference}%
$$
[x_0]\,f := f(x_0) \quad\hbox{and}\quad [x_0,\ldots,x_n]\,f :=
\frac{[x_0,\ldots,x_{n-1}]\,f - [x_1,\ldots,x_n]\,f}{x_0 - x_n},
$$
which is a symmetric function in $x$.
For $W(x) := \prod_{i=1}^d (1-x_i^2)^{-\half}$,  
\begin{equation} \label{prodChebykernel}
    P_n(x,\mathbf{1}) = [x_1,\ldots,x_d]\,G_{n}, 
\end{equation}
where 
$$
G_{n}(t) := 2(-1)^{\left[\half(d+1)\right]}\,(1-t^2)^{\half(d-1)} 
\begin{cases} T_n(t) & \hbox{for $d$ even}, \cr
U_{n-1}(t) & \hbox{for $d$ odd}. \end{cases}
$$

The product Jacobi polynomials inherite a product formula from 
the one-variable case \cite[Lemma 9.6.1]{DunklXu}, which
allows one to define a convolution structure for orthogonal expansions. 

\section{Relation between orthogonal polynomials on classical domains} \label{sec:Sphere-ball-simplex}

By \emph{classical domains} we mean the sphere, ball, simplex, $\RR^d$
and $\RR_+^d$. Orthogonal polynomials
on these domains are closely related.
\index{orthogonal polynomials of several variables!on classical domains!
relations between OPs on different domains}%

\subsection{Orthogonal polynomials on the sphere and on the ball}\label{sec:ball-sphere}
\index{orthogonal polynomials of several variables!on classical domains!
on unit ball}%
\index{orthogonal polynomials of several variables!on classical domains!
on sphere}%
A nonnegative weight function $H$ defined on $\RR^{d+1}$ is called
{\it $S$-symmetric}
\index{weight function!S-symmetric@$S$-symmetric}%
if $H(x', x_{d+1}) = H(x', - x_{d+1})= H(-x',x_{d+1})$, where $x' \in \RR^d$, 
and if the restriction $W_H$
of $H$ on the sphere $\SS^d$
is a non-trivial weight function. Let 
$\CH_n^{d+1}(H)$ denote the space of homogeneous polynomials of degree $n$ that are 
orthogonal in $L^2(\SS^d, H)$ to polynomials of lower degrees.  Then,
just as in the case \eqref{dimspherharm} of $H=1$,
\begin{equation}
   \dim \CH_n^{d+1}(H) = \dim\CP_n^{d+1} - \dim \CP_{n-2}^{d+1} = 
       \binom{n+d}{d} - \binom{n+d-2}{d}.
\end{equation}

Associated with an $S$-symmetric weight function $H$, define 
$$
   W_H(x) := H\Big(x,\big(1-\|x\|^2\big)^\half\Big), \quad x \in \BB^d,
$$
which is a centrally symmetric weight function on the ball.
Let $\{P_\alpha\}_{|\alpha| =n,\,\alpha \in \NN_0^d}$
be an orthogonal basis for $\CV_n^d$ with 
respect to $(1-\|x\|^2)^{-\half}W_H(x)$, and let
$\{Q_\beta\}_{|\beta| = n-1,\,\beta \in \NN_0^d}$ 
be an orthogonal basis for $\CV_{n-1}^d$ with respect to the weight function
$(1-\|x\|^2)^\half W_H(x)$. 
For $y \in \RR^{d+1}$ define 
$$
    y  = r\,(x, x_{d+1}), \qquad x \in \BB^d, \quad (x,x_{d+1}) \in \SS^d, \quad r \ge 0.  
$$
For $\alpha, \beta \in \NN_0^d$, $|\alpha| =n$ and $|\beta| = n-1$, define
$$
     Y_\alpha^{(1)}(y) = r^n P_\alpha(x) \quad \hbox{and} \quad   Y_\beta^{(2)}(y) = r^n x_{d+1} Q_\beta(x). 
$$

\theo{ch2_thk8}
{The functions $Y_\alpha^{(1)}$ and $Y_\beta^{(2)}$ are homogeneous polynomials of degree $n$ 
in the variable $y$.
Furthermore,
$\big\{Y_\alpha^{(1)}\big\}_{|\alpha| =n} \cup
\big\{Y_\beta^{(2)}\big\}_{|\beta| =n-1}$
is an orthogonal basis for $\CH_n^{d+1}(H)$.}

Let $P_n^H(\cdot,\cdot)$ denote the reproducing kernel of $\CH_n^d(H)$ and let $P_n(\cdot,\cdot)$ denote
the reproducing kernel of $\CV_n^d$ with respect to
$(1-\|x\|^2)^{-\half}W_H(x)$. Then
\begin{equation}
P_n (x,y) = \thalf \left(P_n^H \big((x,x_{d+1}),(y, y_{d+1}) \big)+
        P_n^H\big((x,x_{d+1}),(y,- y_{d+1})\big)  \right),
\end{equation}
where $x_{d+1} =(1-\|x\|^2)^\half$ and $y_{d+1} =(1-\|y\|^2)^\half$. 
The relation is based on 
\begin{equation}
 \int_{\SS^d} f(y)\,d\omega(y) = \int_{\BB^d} \left( f\big(x,(1-\|x\|^2)^\half\big)+
 f(x, - (1-\|x\|^2)^\half\big) \right)\,(1-\|x\|^2)^{-\half}\,dx, 
\end{equation}
where $d\omega$ is the Lebesgue measure (not normalized) on $\SS^{d}$. A further relation between
orthogonal polynomials on $\BB^d$ and those on $\SS^{d+m-1}$ follows from 
$$
  \int_{\SS^{d+m-1}} f(y)\,d\omega(y) =
  \int_{\BB^d} (1-\|x\|^2)^{\half(m-1)} 
    \left( \int_{\SS^m} f\big(x, (1-\|x\|^2)^\half \xi\big)\,d\omega_m(\xi) \right)
    \,dx.
$$

As a consequence of these relations, properties of the orthogonal polynomials with respect to the weight function
\begin{equation} \label{W_kBd}
  W_\kappa (x) = \prod_{i=1}^d |x_i|^{2 \kappa_i} (1- \|x\|^2)^{k_{d+1} -\half }, \quad k_i \ge 0, 
\end{equation}
on $\BB^d$ can be derived from $h$-harmonics
\index{h-harmonics@$h$-harmonics}%
with respect to $w_\kappa(x) = \prod_{i=1}^{d+1} |x_i|^{2\kappa_i}$ 
on the sphere $\SS^{d+1}$.
In particular, following \eqref{h-harmonic-basis} and using
generalized Gegenbauer polynomials \eqref{generGegenb},
\index{Gegenbauer polynomials!generalized}%
define 
\begin{equation} \label{W_kBd-basis}
  P_\alpha(x) := [h_\alpha]^{-1} 
  \prod_{j=1}^d (1-\|\xb_{j-1}\|^2)^{\half\alpha_j}\,
  C_{\alpha_j}^{(a_j, \kappa_j)}
    \left(\big(1-\|\xb_{j-1}\|^2\big)^{-\half}x_j \right),\quad
    \alpha \in \NN_0^d,
\end{equation}
where $|\alpha^j| := \alpha_j + \ldots \alpha_d$, 
$|\kappa^j| := \kappa_j + \ldots \kappa_{d+1}$,  
$a_j := |\alpha^{j+1}|+|\kappa^{j+1}|+\frac{d-j}{2}$, and  
$$
\left[h_{\alpha}^n\right]^2 := \frac{1} {\big(|\kappa| + \half(d+1)\big)_n}
\prod_{j=1}^{d} h_{\alpha_i}^{(\a_i, \kappa_i)} (\kappa_i + a_i)_{\alpha_i}. 
$$  
Then $\{P_\alpha\}_{|\alpha| = n}$ is an orthonormal basis of $\CV_n^d$.
A second orthonormal basis 
can be given in spherical-polar coordinates,
analogous to \eqref{ball-base2} but using $h$-harmonics. 
The orthogonal polynomials in $\CV_n^d$ are eigenfunctions of a second order
differential-difference operator, 
\begin{align} \label{diff-diffEqnBd}
   \left(\Delta_h - \la x, \nabla \ra^2 - (2|\kappa| + d-1) \la x, \nabla \ra \right)
   P =  -n (n+2|\kappa|+d-1)\,P,
\end{align}
where  $\Delta_h$ is the Dunkl Laplacian
\index{Dunkl Laplacian}%
given in \S\ref{sec:h-harmonics}.
The reproducing kernel
\index{reproducing kernel}%
$P_n(\cdot,\cdot)$ of $\CV_n^d$ for $W_\kappa$ satisfies a closed formula:
\begin{multline} \label{W_kB-kernel}
 P_n(x,y) =  \frac{2n + 2|\kappa| + d-1}
     {2|\kappa| +d-1}\,  \int_{[-1,1]^{d+1}} 
  C_n^{|\kappa| + \half(d-1)}  (x_1 y_1 t_1 + \cdots + x_{d+1} y_{d+1} t_{d+1} )\\
 \times\left({\textstyle\prod\limits_{i=1}^d} (1+t_i)\right)\,
\left({\textstyle\prod\limits_{i=1}^{d+1}}
    c_{\kappa_i}\big(1-t_i^2\big)^{\kappa_i-1}\right)\, d t,
\end{multline}
where $x_{d+1} = \sqrt{1-\|x\|^2}$ and $y_{d+1} = \sqrt{1-\|y\|^2}$. 
\paragraph{Further results and references}\quad
The relation between orthogonal polynomials on the 
sphere and on the ball was explored in \cite{Xu98b, Xu01a}. The orthonormal
basis \eqref{W_kBd-basis} and the closed form of formula \eqref{W_kB-kernel} were given 
in \cite{Xu01a}.
The monic orthogonal basis for $W_\kappa$ was studied in~\cite{Xu05c}.

\subsection{Orthogonal polynomials on the ball and on the simplex}
\index{orthogonal polynomials of several variables!on classical domains!
on unit ball}%
\index{orthogonal polynomials of several variables!on classical domains!
on simplex}%
Let $W$ be a weight function defined on the simplex $T^d$. Define 
$$
   W_T(x) := {W(x)}/{\sqrt{x_1 \cdots x_d}} \quad \hbox{and} \quad W_B(x) := W(x_1^2,\ldots,x_d^2)
$$
on $T^d$ and on the ball $\BB^d$, respectively. There is a close  relation between 
orthogonal polynomials for $W_T$ and those for $\BB^d$. Let $\CV_n^d(W)$ denote the
space of orthogonal polynomials of degree $n$ with respect to $W$. Furthermore, let 
$\CV_n(W_B,\ZZ_2^d)$ denote the subspace of $\CV_n^d(W_B)$ which contains 
polynomials that are invariant under $\ZZ_2^d$, that is, even in each variable. Then
the mapping 
$$
\psi\colon T^d\to\BB^d\colon
(x_1,\ldots,x_d) \mapsto (x_1^2, \ldots, x_d^2)
$$
induces a ono-to-one correspondence between $R\in \CV_n^d(W_T)$ and
$R\circ \psi \in \CV_n(W_B,\ZZ_2^d)$. This is based on the relation
\begin{equation} \label{integral-Td-Bd}
\int_{B^d} f(x_1^2,\ldots, x_d^2)\,dx = \int_{T^d} f(u_1,\ldots, u_d)\,
\frac{du}{\sqrt{u_1 \cdots u_d}}. 
\end{equation}
Let $P_n(W; \cdot,\cdot)$ denote the reproducing kernel of $\CV_n^d(W)$. Then the 
correspondence also extends to the reproducing kernel: 
\begin{equation} 
 P_n(W_T;x,y) = 2^{-d} \sum_{\varepsilon \in \ZZ_2^d} P_n(W_B; x^{1/2}, \varepsilon y^{1/2}), 
\end{equation}
where $x^{1/2} := \big(x_1^{1/2}, \ldots, x_d^{1/2}\big)$. In
particular, the identity \eqref{reprod-simplex1}
can be deduced from \eqref{W_kB-kernel} by this relation. 

Essentially all properties of orthogonal polynomials for $W_T$ can be deduced from the 
corresponding results for $W_B$. In particular, all results in
\S\ref{sec:OPsimplex}
can be deduced from the corresponding results with respect to $W_\kappa$ in \eqref{W_kBd}. 
In combination with \S\ref{sec:ball-sphere}
 there is also a correspondence between 
orthogonal polynomials on the simplex and those on the sphere. 
\paragraph{Further results and references}\quad
The relation between orthogonal polynomials on the ball and on the
simplex was studied in \cite{Xu98c}.
The details on orthogonal systems for the classical weight functions 
were worked out in \cite{Xu01b}. The connection extends to other
aspects of analysis, including 
orthogonal expansion, approximation, and numerical integration.

\subsection{Limit relations} 
Two limit relations between orthogonal polynomials on two different
domains are worth mentioning.
\index{orthogonal polynomials of several variables!on classical domains!
limits between OPs on different domains}%
\paragraph{Limit of orthogonal polynomials on the ball}\quad
\hskip-0.2cm\cite[Theorem 8.3.5]{DunklXu}\\
\index{orthogonal polynomials of several variables!on classical domains!
on unit ball}%
\index{orthogonal polynomials of several variables!on classical domains!
generalized Hermite}%
Let $P_\alpha(W_\kappa)$ be the orthogonal 
polynomials \eqref{W_kBd-basis} on the ball for $W_\kappa$ in
\eqref{W_kBd}. Then
$$
   \lim_{\kappa_{d+1} \to 0} P_\alpha(W_\kappa, x/ \kappa_{d+1}) =  c_\alpha^{-1} 
     H_{\alpha_1}^{\kappa_1}(x_1)\ldots H_{\alpha_d}^{\kappa_d}(x_d),
$$    
where the right-hand side is the generalized Hermite polynomial
\eqref{GeneralHermitebasis} on $\RR^d$, a product of generalized
Hermite polynomials in one variable.
\paragraph{Limit of orthogonal polynomials on the simplex}\quad
\hskip-0.2cm\cite[Theorem 8.4.5]{DunklXu}\\
\index{orthogonal polynomials of several variables!on classical domains!
on simplex}%
\index{orthogonal polynomials of several variables!on classical domains!
Laguerre}%
Let $P_\alpha(W_\kappa)$ be the orthogonal 
polynomials in \eqref{eq:OPonT}  on the simplex for $W_\kappa$ in
\eqref{eq:weight-Td}.
Then
$$
   \lim_{\kappa_{d+1} \to 0} P_\alpha(W_\kappa, x/ \kappa_{d+1}) =  L_{\alpha_1}^{\kappa_1}(x_1) \cdots
   L_{\alpha_d}^{\kappa_d}(x_d)S_{\!n}
$$    
where the right-hand side is a Laguerre polynomial
\eqref{Laguerre-product} on $\RR_+^d$,
a product of Laguerre polynomials of one variable.
 
\section{Orthogonal expansions and summability} \label{OrthoExpan}
\index{orthogonal polynomials of several variables!on classical domains!
orthogonal expansions}%
As long as polynomials are dense in $L^2(d\mu)$, the standard Hilbert space theory shows that 
the partial sum $S_{\!n} f$ in \eqref{eq:Sn(f)} converges to $f$ in $L^2(d\mu)$ norm. The convergence 
does not hold in general for $S_{\!n} f$ in other norms. The summability of the orthogonal expansions 
is often studied via the Ces\`aro means. For $\delta > 0$,
\index{Cesaro means@Ces\`aro means}%
the \emph{Ces\`aro $(C,\delta)$ means} of the 
orthogonal expansion \eqref{eq:Foruier} is defined by 
\begin{equation}\label{eq:Cesaro}
 S_{\!n}^\d f(x) := {\textstyle\binom{n+\delta}{n}^{-1}}
 \sum_{k=0}^n  \binom{n -k+\delta -1}{n-k}\,S_k f(x). 
\end{equation} 
There are many results for orthogonal expansions for classical type orthogonal polynomials.
Below is a list of highlights with references.  
\paragraph{Orthogonal expansions on $\SS^{d-1}$}\quad
\index{orthogonal polynomials of several variables!on classical domains!
on sphere}%
The results are stated in terms of $w_\k$ defined in \eqref{wk-h-harmonic} on the sphere. The
case $\kappa = 0$ gives the result for ordinary spherical harmonics expansions. Let 
$\|\cdot \|_{\kappa,p}$ denote the $L^p(w_\kappa)$ norm; for $p = \infty$,  the norm is the 
uniform norm of $C(\sph)$.  Let 
$$
\s_\k : = \tfrac{d-2}2 + |\k| -\min_{1\le j \le d} \k_j, \quad 
\d_\k(p):=\max\Big((2\s_\k+1)\,\big|\tfrac1p-\thalf\big|-\thalf,0\Big).
$$

\begin{enumerate}[label=\arabic*.]

\item For $p=1$ or $\infty$, the norm of partial sum operator and the projection operator satisfy 
$\|S_{\!n}\|_{\k,p} \sim \|\proj_n\|_{\k,p}  \sim n^{\s_\k }$.

\item The $(C,\delta)$ means $S_{\!n}^\delta f \ge 0$ for all $f \ge 0$
if and only if $\delta \ge 2 |\k|+ d-1$.

\item If $f \in C(\SS^{d-1})$ then $S_{\!n}^\d f$ converges to $f$
pointwise in 
$\{x \in \sph\mid x_1\ldots x_d\ne0\}$ if $\d > \frac{d-1}{2}$;
if $f \in L^1(w_\k)$ then 
$S_{\!n}^\delta f$ converges almost everywhere to $f$ on $\SS^{d-1}$ if
    $\d > \s_\k$. 

\item For $ p =1$ or $\infty$,  $S_{\!n}^\d f$ converges to $f$ in $L^p(w_\k)$ if and only if $\delta > \sigma_\k$. 

\item If $f\in L^p(w_\k)$, $1\leq p\leq \infty$,
$\big|\frac1p-\half\big|\ge \frac 1{2\s_\k+2}$ and 
$\d>\d_\k(p)$, then $S_{\!n}^\d f$ converges to $f$ in $L^p(w_\k^2)$. 
\end{enumerate}
\vskip\medskipamount

For expansions of ordinary spherical harmonics, these were mostly proved in \cite{BC1973, So}, 
see \cite{Er2} for earlier results. For $h$-harmonics, these results were established in 
\cite{DaiXu3,DaiXu2,LiXu2,Xu97c}.
A comprehensive account that includes many more
results is given in the monograph \cite{DaiXu4}.  
\paragraph{Orthogonal expansions on $\BB^d$ and $T^d$}\quad
\index{orthogonal polynomials of several variables!on classical domains!
on unit ball}%
\index{orthogonal polynomials of several variables!on classical domains!
on simplex}%
For orthogonal expansions with respect to $W_\kappa$ in \eqref{W_kBd} on the ball
$\BB^d$ and $W_\kappa$ in \eqref{eq:weight-Td} on the simplex $T^d$, analogues 
of the above results hold. In fact, if we replace $\s_\k$ by
$$
    \s_\k: = \tfrac{d-1}2 + |\k| -\min_{1\le j \le d+1} \k_j,  
$$
then all five properties hold with obvious modification.
This is no accident, the three cases are intimately 
connected: some of the results on one domain can be deduced from the
corresponding results on one of the other two domains.
The study of summability on the ball started in \cite{Xu99a}. 
The connection between the three domains was first applied in the 
study of orthogonal expansions in \cite{Xu01a,Xu01b}, and next
in full power in 
\cite{DaiXu3,DaiXu2,LiXu2}.
See \cite[Ch.~12]{Er2} for earlier results on orthogonal 
expansions and \cite{DaiXu4} for further results. 
\paragraph{Orthogonal expansions on $[-1,1]^d$.}\quad
\index{orthogonal polynomials of several variables!on classical domains!
product Jacobi}%
On $[-1,1]^d$ we consider the orthogonal expansions in the product Jacobi polynomials 
for $W_{a,b}$ in \eqref{Jacopbi-weight}. 
\vskip\medskipamount
\begin{enumerate}[label=\arabic*.]
\item Let $a_j > -1$, $b_j > -1$ and $a_j + b_j \ge -1$ for $1 \le j \le d$. For $f$ in $L^p(W_{a,b}; [-1,1]^d)$, 
$1 \le p<\infty$, or in $C([-1,1]^d)$, the $(C,\delta)$ means $S_{\!n}^\delta f$ converges to $f$ in norm 
as $n \to \infty$ if 
$$
\textstyle
\delta > \delta_0:=\sum\limits_{j=1}^d \max(a_j,b_j) + \thalf d + 
  \max \left(0, -\sum\limits_{j=1}^d \min(a_j,b_j) -\thalf(d+2) \right).
$$
\item Let $a_j \ge  -1/2$, $b_j  \ge -1/2$ and $a_j + b_j \ge -1$ for $1 \le j \le d$. Then
$S_{\!n}^\delta f \ge 0$ whenever $f \ge 0$ if and only if $\delta \ge \sum_{i=1}^d (a_i + b_i) +3d-1$.
\end{enumerate}
\vskip\medskipamount
These were established in \cite{LiXu1}. In comparison with the results
on the sphere, ball and simplex,
far less is known for orthogonal expansions on $[-1,1]^d$.
The difficulty lies in the lack of a 
closed form of the reproducing kernel.  
\paragraph{Product Hermite and Laguerre expansions}\quad
\index{orthogonal polynomials of several variables!on classical domains!
Hermite}%
\index{orthogonal polynomials of several variables!on classical domains!
Laguerre}%
Hermite expansions on $\RR^d$ for $W_H$ in \eqref{Hermite-weight} and Laguerre expansions
on $\RR_+^d$ for $W_\kappa$ in \eqref{Laguerre-weight} have been extensively
studied. We mention just two results.
\vskip\medskipamount
\begin{enumerate}[label=\arabic*.]
\item The Riesz means $S_R^\delta$ of the product Hermite expansions 
converges in the norm for $f \in L^p(\RR^d)$ ($1 \le p < \infty$)
if $\delta > \half(d-1)$. 
For every $f \in L^p(\RR^d)$, the means $S_R^\d f$ converges to $f$ 
almost everywhere if $\d > \half\big(d -\tfrac13\big)$. 
\item Let $\k_j \ge 0$, $1 \le j \le d$. The $(C,\delta)$ means $S_{\!n}^\delta f$ converges to 
$f \in L^p(W_\kappa)$ in norm if $\delta \ge \sum_{i=1}^d (a_i + b_i) +3d-1$.
For $p = 1$ or $\infty$ the condition on $\delta$ is also necessary.
\end{enumerate}
\vskip\medskipamount
For extensive study on these expansions, see \cite{Thanga} and the
references therein.

\section{Discrete orthogonal polynomials of several variables} \label{sect:discreteOP}
\index{orthogonal polynomials of several variables!discrete}%
Let $V$ be a finite or countable set of points in $\RR^d$ and let
the \emph{weight} $W$ be a positive function on~$V$.
\emph{Discrete orthogonal polynomials} are orthogonal polynomials
with respect to the discrete inner product
\begin{equation}\label{discrete_ip}
  \la f, g\ra := \sum_{x \in V} f(x) g(x) W(x). 
\end{equation}
In case of infinite $V$ assume that $W$ decays fast enough such that
the sum in \eqref{discrete_ip} converges absolutely
for all polynomials $f,g$.

It should be noted that $\la \cdot, \cdot \ra$ in \eqref{discrete_ip}
is an inner product only on 
$\Pi^d / I(V)$, where $I(V)$ is the polynomial ideal of polynomials
which vanish on $V$. 
To be more specific, fix a monomial order and let
$\Lambda(V)$ consist of all $\alpha \in \NN_0^d$ such that 
$c \xb^\alpha$ is not a leading monomial of any polynomial in $I(V)$.
Then the inner product 
\eqref{discrete_ip} is well defined on
$\Pi_V:= \sspan \{x^\beta\}_{\beta \in \Lambda(V)}$
and there exists a sequence 
of orthonormal polynomials $\{P_\a\}_{\a \in \Lambda(V)}$ where
$P_\a$ is a polynomial in $\Pi_V$ with leading monomial $c \xb^\alpha$.
See \cite{Xu04} for details and further discussions. 

Notice, in particular, that the Gram--Schmidt method of generating an
orthonormal basis can 
only be performed within $\Pi_V$.
Beside the above precaution, general properties of
discrete orthogonal polynomials are analogous to the continuous case, 
modulo $I(V)$ if necessary.  

\subsection{Classical discrete orthogonal polynomials}
These polynomials are expressed in terms of the classical discrete
orthogonal polynomials of one
variable, which are listed in Table \ref{table-classical-discrete}. Here
the Krawtchouk and Hahn polynomials have a finite domain
$[0,N] := \{0,1,\ldots, N\}$, and accordingly $n=0,1,\ldots,N$.
The polynomials in the table
can be expressed by hypergeometric functions, and their normalization 
means that the coefficient in front of the hypergeometric function is
the constant 1. 

\begin{table}[h]
\caption{Classical discrete orthogonal polynomials of one variable}
\label{table-classical-discrete} 
\addtolength\tabcolsep{2pt}
\begin{tabular}{@{}c@{}ccccl@{}} 
\hline \hline 
Name & notation & weight & domain & constraint & normalization     \\ 
\hline 
Charlier & $C_n(x;s)$ & $s^x / x!$ & $\NN_0$ & $s>0$ &  $C_n(0;s)=1$
\\[3pt] 
Meixner & $M_n(x;\beta, c)$ & $ \beta^x (c)_x /x!$  & $\NN_0$
&$\beta>0$, $0<c<1$&$M_n(0;\beta, c)=1$ \\[3pt] 
Krawtchouk & $K_n(x; p, N)$ & $p^x(1-p)^{N-x} \binom{N}{x}$ &
$[0,N]$  &  $0<p<1$ &  $K_n(0; p, N)=1$ \\[3pt] 
Hahn & $Q_n(x; \a, \b, N)$ & $\binom{\a+x}{x}\binom{\b+N-x}{N-x}$  &
$[0,N]$&$\alpha,\beta>-1$ or $<-N$&$Q_n(0; \a, \b, N)=1$ \\
\hline \hline 
\end{tabular} 
\end{table} 
\index{orthogonal polynomials of one variable!classical discrete}%
\index{orthogonal polynomials of one variable!classical discrete!Charlier}%
\index{Charlier polynonials}%
\index{orthogonal polynomials of one variable!classical discrete!Meixner}%
\index{Meixner polynomials}%
\index{orthogonal polynomials of one variable!classical discrete!Krawtchouk}%
\index{Krawtchouk polynomials}%
\index{orthogonal polynomials of one variable!classical discrete!Hahn}%
\index{Hahn polynomials}%
There are several ways to extend classical discrete orthogonal polynomials to 
several variables. One way is to consider those families for which all orthogonal 
polynomials of degree exactly $n$ eigenfunctions of a difference
operator of specific form:
\begin{equation}\label{discreteDE}
    \left(\sum_{1\leq i, j\leq d}A_{i,j}(x)\,\fs_i\bs_j  +
    \sum_{i=1}^dB_i(x)\,\fs_i\right) \psi(x)  = \l_n\,\psi(x), 
\end{equation}
where $\fs_i$ and $\bs_i$ denote the forward and backward
difference operators 
$$
\fs_i f(x):=f(x+e_i)-f(x) \quad \hbox{and} \quad \bs_i f(x):=f(x)-f(x-e_i)
\quad\mbox{($e_1,\ldots,e_d$ standard basis)},
$$ 
It follows readily that the $A_{i,j}(x)$ are
necessarily quadratic polynomials and the $B_i(x)$ are necessarily linear polynomials in $x$. 
Some of these families are tensor products of classical polynomials of one variable, 
which will be discussed in the next subsection. Below are several families that do not come 
from products.  

For a vector $x=(x_1,x_2,\dots,x_d)\in\RR^d$ we will denote 
\begin{equation}\label{x^jX_j}
   x^j: =(x_j,x_{j+1},\dots,x_d) \quad \text{ and } \quad X_j:=(x_1,x_2,\dots,x_j),
\end{equation}
with the convention that $x^{d+1}:=0$ and $X_0:=0$.  Set $|x| := x_1 + \cdots + x_d$. 
\paragraph{Meixner polynomials}\hskip0.2cm\cite[\S6.1.2]{IX}\\
\index{orthogonal polynomials of several variables!classical discrete!
Meixner}%
Let $0< c_i < 1$ ($1 \le i \le d$) such that $|c|<1$ and let
$s>0$. For $\nu \in \NN_0^d$
define the polynomials
\begin{align}\label{Meixner}
  M_\nu(x; s, c) =   \prod_{j=1}^d  (\delta_j)_{\nu_j} 
               M_{\nu_j}\left(x_j; \delta_j , \frac{c_j}{1-|c^{j+1}|} \right),
               \quad\mbox{where $\delta_j := s+|\nu^{j+1}|+|X_{j-1}|$.}
\end{align}
They satisfy the orthogonality relation
\begin{equation}\label{MeixnerOrtho}
  \sum_{x \in \NN_0^d} M_\nu(x; s,c)M_\mu(x;s,c)\,
    (s)_{|x|}\,{\textstyle\prod\limits_{i=1}^{d}}\,\frac{c_i^{x_i}}{x_i!}= 
\frac{(s)_{|\nu|}}{(1-|c|)^s}\,
\left({\textstyle\prod\limits_{j=1}^d}\,\nu_j!\,
 \left(\frac{1-|c^{j+1}|}{c_j}\right)^{\nu_j}\right)
          \delta_{\nu,\mu}. 
\end{equation}
Furthermoere, $\psi_\nu :=   M_\nu(\cdot ; s, c)$ is an eigenfunction
of a difference operator: 
\begin{equation} \label{MeixnerDQ}
D \psi_\nu =  - |\nu|\,\psi_\nu,\quad D := \sum_{1\leq i,j\leq d}
\left(\delta_{i,j}+\frac{c_i}{1-|c|}\right)x_j \fs_i\bs_j 
+ \sum_{i=1}^d \left(- x_i + \frac{c_i}{1-|c|} s\right)\fs_i . 
\end{equation}

\paragraph{Krawtchouk polynomials}\hskip0.2cm\cite[\S6.1.1]{IX}\quad
\index{orthogonal polynomials of several variables!classical discrete!
Krawtchouk}%
Let $0< p_i < 1$ ($1 \le i \le d$) such that $|p|<1$ and let
$N \in \NN$.
For $\nu\in\NN_0^d$, $|\nu| \le N$, define the polynomials
\begin{multline}\label{Krawtchouk}
  K_\nu(x; p, N) :=  \frac{(-1)^{|\nu|} }{(-N)_{|\nu|}}\,
    \prod_{j=1}^d  \frac{p_j^{\nu_j}} {(1- |P_j|)^{\nu_j}}
       (-N+|X_{j-1}|+|\nu^{j+1}|)_{\nu_j}\\
       \times K_{\nu_j}\left(x_j; \frac{p_j}{1- |P_{j-1}|}, 
         N-|X_{j-1}|- |\nu^{j+1}|\right).
\end{multline}
They satisfy the orthogonality relation
\begin{equation}\label{KrawtchoukOrtho}
 \sum_{|x| \le N} K_\nu(x;p,N)K_\mu(x;p,N) \,
{\textstyle\prod\limits_{i=1}^{d+1}}\, \frac{p_i^{x_i}}{x_i!}=\frac{(-1)^{|\nu|}}{(-N)_{|\nu|}\,N!}\,
\left({\textstyle\prod\limits_{j=1}^d}\, \frac{ \nu_j!\,p_j^{\nu_j} } 
{(1-|P_j|)^{\nu_j-\nu_{j+1}} }\right) \delta_{\nu,\mu}, 
\end{equation}
where $x_{d+1}=N-|x|$, $p_{d+1}=1-|p|$ and $\nu_{d+1}=0$.
Furthermore,
$\psi_\nu :=   K_\nu(\cdot; p, N)$ is an eigenfunction of a
difference operator: 
\begin{equation} \label{KrawtchoukDQ}
D \psi_\nu =  - |\nu|\,\psi_\nu,\qquad
 D:= \sum_{1\leq i,j\leq d}(\delta_{i,j}-p_i)x_j \fs_i\bs_j 
 + \sum_{i=1}^d (p_iN-x_i) \fs_i. 
\end{equation}

\paragraph{Hahn polynomials on the parallelepiped}
\index{orthogonal polynomials of several variables!classical discrete!
Hahn on parallelepiped}%
\hskip0.2cm\cite[\S5.1.2]{IX}\\
Let $l_i \in \NN$ ($1 \le i \le d$), $\beta > -1$ and $r > 0$.
For $\nu \in \NN_0^d$, $\nu_i\le l_i$, define the 
polynomials
\begin{align} \label{5.25}
\phi_\nu(x;\beta,r,l)&:=
   \prod_{i=1}^d (\alpha_{1,j}+1)_{\nu_j} Q_{\nu_j}(x_j; \alpha_{1,j},   
            \alpha_{2,j},l_j),\\
&\alpha_{1,j}  := \beta+ |\nu^{j+1}|+|X_{j-1}| \quad \hbox{and}\quad 
   \alpha_{2,j}  :=|L_{j-1}| - |X_{j-1}|+|\nu^{j+1}|+ r -1\nonumber
\end{align}
(recall the notation \eqref{x^jX_j}).
They satisfy the orthogonality relation
\begin{multline}\label{HanhVnParallep}
 \sum_{x\le l} \phi_{\nu}(x;\beta,r,l)\,\phi_{\mu}(x;\beta,r,l)\,\left({\textstyle\prod\limits_{i=1}^d}\,
   \frac{(- l_i)_{x_i} }{x_i!}\right)\,
   \frac{(\beta+1)_{|x|}}{(-|l |-r+1)_{|x|}}\\
=  \frac{(-1)^{|\nu|}\,(1+\beta)_{|\nu|}}{(r+|\nu|)_{|l|-|\nu|}}\,
\left({\textstyle
\prod\limits_{j=1}^d}\,\frac{\nu_j!\,(\beta+r+2|\nu^{j+1}|+\nu_j+|L_{j-1}|)_{l_j+1}}
{(-l_j)_{\nu_j}(\beta+r+2|\nu^{j}|+|L_{j-1}|)}\right)
\,\delta_{\nu,\mu}.
\end{multline}
Furthermore, they satisfy the difference equation
$D \phi_\nu =  - |\nu|(|\nu|+\beta+ r) \phi_\nu$ with 
\begin{align} \label{HahnVndDQ}
D=\sum_{1\leq i,j\leq d}x_j[l_i - x_i + r \delta_{i,j}]
\fs_i\bs_j
 - \sum_{i=1}^d[x_i(r+\beta +1) - l_i (1+\beta)]\fs_i.  
\end{align}
Furthermore, these relations also hold if $\beta < - |l|$ and
$r < - |l|+1$. 
\paragraph{Hahn polynomials on the simplex}\hskip0.2cm
\index{orthogonal polynomials of several variables!classical discrete!
Hahn on simplex}%
\cite[\S5.2.1]{IX}\\
Let $\sigma_i > -1$ ($1 \le i \le d+1$)
and $N \in \NN$. For $\nu \in \NN_0^d$ such that $|\nu|\le N$,
define the polynomials 
 \begin{equation} \label{HahnVnd}
 Q_\nu(x;\sigma, N) :=   \frac{(-1)^{|\nu|}}{(-N)_{|\nu|}}
  \prod_{j=1}^d  \frac{(\sigma_j+1)_{\nu_j}}{(a_j+1)_{\nu_j}}\,
    (-N + |X_{j-1}|+|\nu^{j+1}|)_{\nu_j}\,
Q_{\nu_j}(x_j; \sigma_j, a_j, N-|X_{j-1}|-|\nu^{j+1}|),
\end{equation}
where $a_j := |\sigma^{j+1}|+2 |\nu^{j+1}| + d-j$. They satisfy the
orthogonality relation
\begin{multline} \label{HanhVndOrtho}
\sum_{|x|\le N} Q_\nu(x; \sigma,N)\,Q_\mu(x; \sigma, N)\,
\left({\textstyle\prod\limits_{i=1}^d}\,\binom{x_i+\sigma_i}{x_i}\,\right)\,\binom{N-|x|+\sigma_{d+1}}{N-|x|}\\
=\frac{(-1)^{|\nu|}\,(|\sigma|+d+2|\nu|+1)_{N-|\nu|}}
   {(-N)_{|\nu|}\, N!}\,
\left({\textstyle\prod\limits_{j=1}^d}\, \frac{(\sigma_j+a_j+\nu_j+1)_{\nu_j}\,(\sigma_j+1)_{\nu_j}\,\nu_j!}
       {(a_j+1)_{\nu_j}}\right)\,\delta_{\nu,\mu}.
\end{multline}
Furthermore, they are eigenfunctions of a difference operator:
\goodbreak
\begin{align} &D \psi_\nu =  - |\nu|(|\nu|+|\sigma|+d)\,\psi_\nu,\quad
D := \sum_{i=1}^d x_i (N-x_i+|\sigma|-\sigma_i+d)\fs_i\bs_i
\nonumber\\[-0.2cm]
&\qquad\quad-\sum_{1 \le i\ne j \le d}x_j(x_i+\sigma_i+1)\fs_i\bs_j
+\sum_{i=1}^d \Big((N-x_i)(\sigma_i+1)- x_i(|\sigma|-\sigma_i+d)\Big)
 \fs_i.
\label{HahnSimplexDQ}
\end{align}
These relations also hold if $\sigma_i < - N$ for $i = 1,2,\ldots, d+1$. 
\paragraph{Hahn polynomials on the simplex-parallelepiped}
\index{orthogonal polynomials of several variables!classical discrete!
Hahn on simplex-parallelepiped}%
\hskip0.2cm\cite[\S5.2.2]{IX}\\
Let $S$ be a nonempty set in $\{1,2,\ldots, d\}$,
let $l_i \in \NN_0$ ($i \in S$) and let $N \in \NN$. Define 
$V_{N,S}^d: = \{x\mid |x| \le N\} \cap \{x\mid x_i \le \ell _i \,
    \hbox{for $i \in S$} \} 
$
and set $\sigma_i = -l_i -1$ ($1 \le i \le d$). 
For $\nu \in V_{N,S}^d$ the polynomials $Q_\nu(\cdot; \sigma, N)$ in 
\eqref{HahnVnd} satisfy the orthogonality relation
\begin{equation} \label{HanhVnS}
\sum_{x\in V_{N,S}^d} Q_\nu(x; \sigma,N)\,Q_\mu(x; \sigma, N)\,
\left({\textstyle\prod\limits_{i=1}^{d}}\,\binom{x_i+\sigma_i}{x_i}\,
\right)\,\binom{N-|x|+\sigma_{d+1}}{N-|x|}
 = A_\nu\,\delta_{\nu,\mu}, 
\end{equation}
where $A_\nu$ is the coefficient of $\delta_{\nu,\mu}$ in
\eqref{HanhVndOrtho}, and they satisfy the
same difference equation as in \eqref{HahnSimplexDQ}. 
\paragraph{Further results and references}\quad
The multivariate Krawtchouk polynomials were first studied in
\cite{M} and
the Hahn polynomials on the simplex were first studied in
\cite{KaMc}. Both classes of polynomials are
associated with linear growth model of birth and death process.
Biorthogonal systems
of Hahn polynomials were found in \cite{Trat-a}. 
The Meixner and Krawtchouk polynomials were studied in \cite{Trat-c} 
They can be deduced as limits of biorthogonal or orthogonal
Hahn polynomials (\cite{Trat-c, Trat2}).
For example, for $Q_\nu$ in \eqref{HahnVnd} and $K_\nu$ in
\eqref{Krawtchouk},
$$
  \lim_{t \to \infty} Q_\nu(x; p_1 t, \ldots, p_d t, (1-p_1-\cdots - p_d)t, N)
    = K_\nu(x; p, N) 
$$
follows from the one-variable case. 

All polynomials in this section were also studied in \cite{IX} in
connection with 
difference equations. There is one more family of discrete orthogonal polynomials 
$\{R_\nu\}$ that resemble the Hahn polynomials studied in \cite{IX}. They satisfy
the orthogonal relation
\begin{equation} \label{Hahn-likeOP}
\sum_{x \in \NN_0^d}  R_\nu(x; \sigma,\beta,\gamma)\,
                    R_\mu(x; \sigma,\beta,\gamma)\,
\left({\textstyle\prod\limits_{i=1}^d} \frac{(\sigma_i+1)_{x_i}}{x_i!}
\right)\,
   \frac{(\beta+1)_{|x|}}{(\gamma+1)_{|x|}} 
 = A_\nu\,\delta_{\nu,\mu},
\end{equation}
where $R_\nu(\cdot; \sigma, \beta, \gamma)$ is defined for
$\nu \in \NN_0^d$ such that
$2|\nu| < \gamma - |\sigma|-\beta - d -1$, and they are also eigenfunctions of a 
second order difference operator. Furthermore, together with product type polynomials 
to be discussed in the following subsection, the discrete orthogonal polynomials in this 
subsection yield all orthogonal polynomial eigenfunctions of a fairly
general class of difference operators 
\eqref{discreteDE}.  

\subsection{Product orthogonal polynomials} 
\index{orthogonal polynomials of several variables!discrete product}%
By taking products of classical discrete orthogonal polynomials
in one variable
one can generate many different products of orthogonal polynomials in
several
variables. Below is a list of such polynomials when $d=2$, of the form
$u(x_1,x_2)=p_k(x_1)\,q_{n-k}(x_2)$ ($0 \le k \le n$), that satisfy the 
difference equation \eqref{discreteDE} with eigenvalue $\lambda_n=n$.. 
\paragraph{Charlier--Charlier}\quad
\index{orthogonal polynomials of two variables!discrete product
orthogonal polynomials!Charlier--Charlier}%
The polynomials $C_k(x_1;a_1)\,C_{n-k}(x_2;a_2)$ ($0 \le k \le n$)
are orthogonal with respect to the weights
$\frac{a_1^{x_1} }{x_1!}\,\frac{a_2^{x_2} }{x_2!}$
($a_1,a_2>0$) on $\NN_0^2$.
They satisfy
\begin{align} 
    x_1 \fs_1 \bs_1 u  +  x_2 \fs_2 \bs_2 u 
 + (a_1-x_1)\fs_1 u  + (a_2-x_2)\fs_2 u  = - n u. 
\end{align}
\paragraph{Charlier--Meixner}\quad
\index{orthogonal polynomials of two variables!discrete product
orthogonal polynomials!Charlier--Meixner}%
The polynomials $M_k(x_1;\beta,c)\,C_{n-k}(x_2;a)$  ($0 \le k \le n$)
are orthogonal 
with respect to the weights
$ \frac{(\beta)_{x_1}}{x_1!} c^{x_1} \frac{(a)_{x_2}}{x_2!}$ 
($a,\beta>0$, $0<c<1$) on $\NN_0^2$.
They satisfy
\begin{equation} 
(c-1)^{-1}x_1 \fs_1 \bs_1 u  -  x_2 \fs_2 \bs_2 u
+ (c-1)^{-1}\big(c(x_1 + \beta)- x_1\big)\fs_1 u 
 -  (a-x_2) \fs_2 u 
 =  n u.
\end{equation}
\paragraph{Charlier--Krawtchouk}\quad
\index{orthogonal polynomials of two variables!discrete product
orthogonal polynomials!Charlier--Krawtchouk}%
The polynomials  $K_k(x_1;p, N)\,C_{n-k}(x_2;a)$  ($0 \le k \le n$)
are orthogonal with respect to the weights
$ \binom{N}{x_1} p^{x_1}(1-p)^{N-x_1} \frac{a^{x_2}}{x_2!}$
($a>0$, $0<p<1$) on  $[0,N]\times\NN_0$.
They satisfy 
\begin{equation} 
(1-p) x_1 \fs_1 \bs_1 u  +    x_2 \fs_2 \bs_2 u
+ \big(p(N-x_1)- (1-p)x_1 \big)\fs_1 u 
 + (a-x_2) \fs_2 u  = - n u.
\end{equation}

\paragraph{Meixner--Meixner}\quad The polynomials
\index{orthogonal polynomials of two variables!discrete product
orthogonal polynomials!Meixner--Meixner}%
$M_k(x_1; \beta_1,c_1)\,M_{n-k}(x_2; \beta_2,c_2)$ ($0 \le k \le n$)
are orthogonal with respect to the weights
$\frac{(\beta_1)_{x_1}}{x_1!}c_1^{x_1}\frac{(\beta_2)_{x_2}}{x_2!}c_2^{x_2}$
($\beta_1,\beta_2>0$, $c_1,c_2\in(0,1)$)
on $\NN_0^2$. They satisfy
\begin{multline} 
(c_1-1)^{-1} \fs_1 \bs_1 u  + 
   (c_2-1)^{-1} x_2 \fs_2 \bs_2 u  \\
 +(c_1-1)^{-1} \big(c_1(x_1 + \beta_1)- x_1\big)\fs_1 u 
 + (c_2-1)^{-1}\big(c_2(x_2 + \beta_2)- x_2\big)\fs_2 u 
   =  n u.
\end{multline}

\paragraph{Meixner--Krawtchouk}\quad
\index{orthogonal polynomials of two variables!discrete product
orthogonal polynomials!Meixner--Krawtchouk}%
The polynomials $M_k(x_1;\beta,c)\,K_{n-k}(x_2;p, N)$ ($0 \le k \le n$)
are orthogonal with respect to the weights
$\frac{(\beta)_{x_1}}{x_1!} c^{x_1}
\binom{N}{x_2} p^{x_2} (1-p)^{N-x_2}$
($\beta>0$,\;\;$c,p\in(0,1)$) on $\NN_0\times[0,N]$. 
They satisfy
\begin{multline} 
(c-1)^{-1} x_1 \fs_1 \bs_1 u  - 
      (1-p) x_2 \fs_2 \bs_2 u\\
  +(c-1)^{-1} \big(c(x_1 + \beta)- x_1\big)\fs_1 u 
 - \big(p (N-x_2) - (1-p) x_2\big)\fs_2 u 
 =  n u. \notag
\end{multline}

\paragraph{Krawtchouk--Krawtchouk}\quad
\index{orthogonal polynomials of two variables!discrete product
orthogonal polynomials!Krawtchouk--Krawtchouk}%
The polynomials $ K_k(x_1;p_1,N_1)K_{n-k}(x_2;p_2,N_2)$ ($0 \le k \le n$)
are orthogonal with respect to the weights
$\binom{N_1}{x_1} p_1^{x_1}(1-p_1)^{N_1-x_1}
 \binom{N_2}{x_2} p_2^{x_2}(1-p_2)^{N_2-x_2}$ ($p_1,p_2\in(0,1)$)
on $[0,N_1]\times[0,N_2]$. They satisfy
\begin{multline} 
 (1-p_1) x_1 \fs_1 \bs_1 u  + 
    (1-p_2) x_2 \fs_2 \bs_2 u  \\
 + \big(p_1(N_1-x_1)- (1-p_1)x_1 \big)\fs_1 u 
 + \big(p_2(N_2-x_2)- (1-p_2)x_2 \big)\fs_2 u = - n u.
\end{multline}
 
For $d > 2$, there are many more product discrete orthogonal polynomials that satisfy 
the second order difference equations. In fact, besides the product of classical one 
variable polynomials, there are also product of classical polynomials of one variable
and other lower dimensional orthogonal polynomials. For example, the product 
of Meixner polynomials on the simplex with either Charlier, Meixner, or Krawtchouk 
polynomials are discrete orthogonal polynomials, of three variables, and they satisfy
difference equations of the form $D u = \l_n u$, where $n$ is the total degree 
of the orthogonal polynomials. For further discussions and details, see \cite{IX, Xu04, Xu05b}. 

There are also product orthogonal polynomials that are given by
products which have a Hahn
polynomial as one of their factors. Such polynomials, however, are
eigenfunctions of a difference operator with eigenvalues not just
depending on the total degree $n$ but also (in the two-variable case)
on $k$.

\vspace{-0.8cm}
\subsection{Further results on discrete orthogonal polynomials}

{\it Racah polynomials}.\quad
\index{orthogonal polynomials of several
variables!classical discrete!Racah}%
These are defined via ${}_4F_3$ functions 
and are orthogonal with respect to weights on $[0,N]$. 
They have bee extended to several variables in \cite{Trat2} for the
weights
\begin{align*}
  &w(x) =  w(x; c_1,\ldots,c_{d+1},\gamma,N):= 
  \frac{N!\,\Gamma(|C_d| + N+1)}{\Gamma(c_{d+1}+N)\,\Gamma(|c|  + N)}\,
      \frac{(c_1)_{x_1} (\gamma+1)_{x_1}}{x_1!\,(c_1 - \gamma)_{x_1}} \\
         &\times \prod_{k=1}^d \frac{\Gamma(c_{k+1}+x_{k+1}-x_k)\,
         \Gamma(|C_{k+1}|+x_{k+1}+x_k)}
          {(x_{k+1} - x_k)! \,\Gamma(|C_k| + x_{k+1}+x_k+1)}\,
          \frac{|C_k| + 2 x_k}{|C_k|}\quad
          (x \in \NN_0^d,\;0\le x_1 \le \cdots \le x_d \le N),
\end{align*}
where $c = (c_1,\ldots,c_{d+1})$, $C_k$ is defined as in \eqref{x^jX_j},
and $x_{d+1} = N$.

In \cite{Trat2}, multivariable dual Hahn polynomials are defined as
limit cases of 
Racah polynomials. The Hahn polynomials on the simplex are also
contained as 
a limit case of the Racah family.
The multivariable Racah polynomials are studied  in view of 
bispectrality in \cite{GI}.

Griffiths \cite{Griffiths} used a generating function to define polynomials in $d$ variables 
orthogonal with respect to the multinomial distribution, which gives a family of Krawtchouk polynomials
that satisfy several symmetric relations among their variables and parameters. These polynomials
are related to  character algebras and the Aomoto--Gel'fand
hypergeometric function in \cite{MiTa} (see also this volume,
Chapter 4, \S4.5).
The recurrence 
relations as well as the reductions which lead to the polynomials defined by Milch \cite{M} and 
Hoare--Rahman \cite{HR} can be found in~\cite{Iliev}. Some of these properties are explored in \cite{IT},
in which these polynomials are interpreted in terms of the Lie algebra $\mathrm{sl}{}_3(\CC)$. 
For applications of these bases of polynomials in probability,
see \cite{DGr}.

Orthogonal polynomials for
the negative multinomial distribution are Meixner polynomials.
A general family of these polynomials 
was defined in terms of generating functions in \cite{Griffiths2}
and their properties were studied in \cite{Iliev2}. 

\section{Other orthogonal polynomials of several variables} \label{sec:OtherOP}
This section contains several families of orthogonal polynomials of several variables 
that are not classical type but can be constructed explicitly. 

\subsection{Orthogonal polynomials from symmetric functions} \label{sec:SymmFuncOP}
A polynomial $f \in \Pi^d$ is called \emph{symmetric} if $f$ is invariant
under any permutation 
of its variables. \emph{Elementary symmetric polynomials} are given by
\index{elementary symmetric polynomials}
\begin{equation*}
  E_k(x_1, \ldots x_d) := \sum_{1\le i_1 < \cdots < i_k \le d}
x_{i_1}\ldots x_{i_k}, \qquad k = 1, 2, \ldots, d.
\end{equation*}
They generate the algebra of symmetric polynomials. 

The mapping
\begin{equation}\label{eq:symm-map}
 x \mapsto u,\quad  u_i:= E_i(x_1, \ldots x_d)\quad (i=1,\ldots,d)
\end{equation} 
is a bijection from the region
$S:=\{x \in \RR^d\mid x_1< x_2 < \cdots <x_d\}$
onto its image $\Omega$.
The Jacobian of this mapping is $J(x): =  \prod_{1\le i <j \le d} (x_i - x_j)$.
The square
of $J(x)$ becomes a polynomial $\Delta(u)$
in $u$ under the mapping \eqref{eq:symm-map},
$$
\Delta(u) = \prod_{1\le i <j \le d} (x_i - x_j)^2, \qquad u \in \Omega.
$$ 
Let $d\mu$ be a nonnegative measure on $\RR$. Define the measure
$d\nu$ on 
$\Omega$ as the image of the product measure
$d\mu(x_1)\ldots  d\mu(x_d)$
under the mapping \eqref{eq:symm-map}.
The orthogonal polynomials with respect to the measures 
$\big(\Delta(u)\big)^{\pm1/2} d \nu$
can be given in terms of orthogonal polynomials with respect to
$d\mu$ on $\RR$, as we will describe now.

\index{orthogonal polynomials of several variables!symmetrized products}%
Let $\{p_n\}$ be orthonormal polynomials for the measure
$d\mu$ on 
$\RR$. For $n \in \NN_0$ and $\alpha \in \NN_0^d$ such 
that $n = \alpha_d \ge \cdots \ge \alpha_1 \ge 0$, define
\begin{equation}\label{eq:-1/2OP}
    P_{\alpha}^{n,-\half}(u) : = \sum_{\beta}
 p_{\alpha_1}(x_{\beta_1})\ldots p_{\alpha_{d}}(x_{\beta_{d}}),  \qquad u 
  \in \Omega, 
\end{equation} 
where the summation is performed 
over all permutations $\beta$ of $\{1,2,\ldots,d\}$. These are polynomials
of degree $n$ and satisfy 
\begin{equation}\label{eq:-1/2OP-ortho}
\int_\Omega  P_{\alpha}^{n, -\half}(u)\,P_{\beta}^{m, -\half}(u)\,
  \big(\Delta(u)\big)^{-\half}\,d\nu(u) = m_1! \ldots m_{d'}! \,\delta_{n,m}\,\delta_{\a,\b},
\end{equation}
where $d'$ is the number of distinct elements in $\a$ and $m_i$ is the number
of occurrences of the $i$th distinct element in $\alpha$. 

For $n \in \NN_0$ and $\alpha \in \NN_0^d$ such 
that $n = \alpha_d \ge \cdots \ge \alpha_1 \ge 0$, define
\begin{equation}\label{eq:1/2OP}
P_{\alpha}^{n,\half}(u) =  \frac{J_{\alpha}^n(x)}{J(x)}, \quad 
\hbox{where} \quad J_{\alpha}^n(x) := \det\, [p_{\alpha_i+d-i}(x_j)]_{i,j=1}^d. 
\end{equation}
These are indeed polynomials of degree $n$ in $u$ under \eqref{eq:symm-map}
and satisfy 
\begin{equation}\label{eq:1/2OP-ortho}
\int_\Omega  P_{\alpha}^{n, \half}(u)\,P_{\beta}^{m, \half}(u)\,
 \left(\Delta(u)\right)^{\half}\,d\nu(u) = \delta_{n,m}\,\delta_{\a,\b}.
\end{equation}

Both these families of orthogonal polynomials satisfy a striking property 
that the polynomials in $\CV_n^d$ have $\dim \Pi_{n-1}^d$ distinct real 
common zeros. In other words, the Gaussian cubature formula exist for
$ \big(\Delta(u)\big)^{\pm1/2}\,d\nu(u)$ by Theorem \ref{thm:GaussCuba}.
\index{cubature formula!Gaussian}%
For $d =2$, these are Koornwinder's polynomials as discussed in
\S\ref{sec:Koornwinder}. For $d > 2$ they were studied in \cite{BSX2}. 

\subsection{Orthogonal polynomials associated with root system $\mathrm A_d$}
\label{sec:AdOP}
Using homogeneous coordinates $\tb= (t_1,\ldots, t_{d+1}) \in \RR^{d+1}$
satisfying the relation $t_1 + \ldots + t_{d+1} =0$, the space $\RR^d$ can be 
identified with the hyperplane 
$$
\RR^{d+1}_H : = \{\tb \in \RR^{d+1}\mid t_1 + \cdots +t_{d+1} = 0\}
$$ 
in $\RR^{d+1}$. The reflection group $\CA_d$ for the root system $\mathrm A_d$
is generated
by the reflections $\sigma_{ij}$, under homogeneous coordinates,
 defined by
$$
 \tb \sigma_{ij} := \tb - 2 \frac{\la \tb, \eb_{i,j}\ra}{\la \eb_{i,j},\eb_{i,j}\ra}\,\eb_{i,j}
     = \tb - (t_i-t_j)\eb_{i,j},
\quad \hbox{where $\eb_{i,j} := e_i- e_j$}.
$$
This group can be identified with the symmetric group of $d+1$ elements. 
Define the operators $\CP^+$ and $\CP^-$ by 
\begin{align*}
  \mathcal{P}^{\pm}f(\tb) := \frac1{(d+1)!}\,\left(\,{\textstyle\sum\limits_{\sigma\in \CA^{+}}}
   f(\tb \sigma)\pm
   {\textstyle\sum\limits_{\sigma\in \CA^{-}}} f(\tb \sigma) \right),
\end{align*}
where $\CA^+$ (or $\CA^-$) condists of an even (or odd) number of
products of reflections  $\sigma_{ij}$. They 
map $f$ to $\CA_d$-invariant ($\CP^+$) or anti-invariant ($\CP^-$) functions,
respectively. 

Let $\HH := \{\kb \in \ZZ^{d+1} \cap \RR_H^{d+1}\mid 
k_1 \equiv \cdots \equiv k_{d+1} \mod {(d+1)}\}$. 
The functions
\[
\phi_\kb(\tb) : =  e^{2\pi i\,(d+1)^{-1} \langle\kb,\tb\rangle}, \qquad
\kb \in \HH, 
\]
are periodic functions:\;\;$\phi_\kb (\tb) = \phi_\kb( \tb+ \jb)$\;\;
($\jb \in \ZZ^{d+1} \cap \RR^{d+1}_H$). 
Let
\[
\Lambda: = \big\{\kb \in \HH\mid k_1 \ge  k_2 \ge\ldots \ge k_{d+1}\big\}\quad
{\rm and}\quad
\Lambda^\circ: =\big\{\kb \in \HH\mid k_1 >   k_2 > \ldots > k_{d+1}\big\}.
\]
Then the functions defined by
$$
    \TC_\kb(\tb) :=  \CP^+ \phi_\kb(\tb)\quad(\kb \in \Lambda)
    \qquad\hbox{and}\qquad 
    \TS_\kb(\tb): = i^{-1}\,\CP^- \phi_\kb(\tb)\quad(\kb \in \Lambda^\circ)   
$$
are invariant and anti-invariant functions, respectively, and they are 
analogues of cosine and sine functions that are orthogonal on the simplex
\begin{align*}
  \Delta := \left\{ \tb \in \RR_H^{d+1}\mid
  0\leq t_i-t_j \le 1\;(1\le i<j\le d+1)\right\}.
\end{align*}
These functions become,
under the change of variables $\tb \mapsto z$, orthogonal polynomials, where
$z_1,\ldots, z_d$ denote the first $d$ elementary symmetric functions 
of $e^{2\pi i t_1}, \ldots, e^{2\pi i t_{d+1}}$. Indeed, for the index $\a \in \NN_0^d$ 
associated to $\kb \in \Lambda$ by 
\begin{align*}
     \a_{i} =\a_i(\kb):= \frac{k_i-k_{i+1}}{d+1},\qquad 1\le i\le d,
\end{align*}
we define under the change of variables $\tb \mapsto z =(z_1,\ldots,z_d)$, 
\begin{align*}
  T_{\a}(z) : = \TC_{\kb}(\tb) \quad \hbox{and} \quad
    U_{\a}(z) := \frac{\TS_{\kb+\vb^{\circ}}(\tb)}{\TS_{\vb^{\circ}}(\tb)}, 
\end{align*}
where $\vb^{\circ}: =(d+1)\,\big(\half d,\half d-1, \ldots,-\half d\big)$.
Then $T_\a$ and $U_\a$
are polynomials in $z$ of degree $|\alpha|$.
\index{Chebyshev polynomials!$\mathrm A_d$}%
\index{orthogonal polynomials of several variables!$\mathrm A_d$}%
These polynomials are analogues 
of Chebyshev polynomials of the first kind and the second kind, respectively. 
In particular, they are orthogonal on the domain $\Delta^*$,  the image of
$\Delta$ under $\tb \mapsto z$,
\begin{equation*}
 \Delta^* := \Big\{x = x(\tb)\in \RR^d\mid\tb \in \RR_H^{d+1},
 {\textstyle\prod\limits_{1\le i< j\le d+1}}
  \sin\big(\pi (t_i-t_j)\big) \ge 0 \Big\}
\end{equation*}
with respect to the weight function $W_{-1/2}$ and $W_{1/2}$,
respectively, where
$$
   W_\alpha(z) := {\textstyle\prod\limits_{1\le \mu < \nu \le d+1}}\,
   \big|\sin\big(\pi (t_{\mu}-t_{\nu})\big)\big|^{2\alpha}. 
$$
For $d =2$, these are the second family of Koornwinder's polynomials
 in \S\ref{sec:Koornwinder}. 

These polynomials satisfy simple recurrence relations and the relation
\begin{equation*}
  \overline{P_{\a}(z)}  = P_{\a_d,\a_{d-1},\dots, \a_1}(z), \quad \hbox{$P_\a = T_\a$ or $U_\a$},
    \quad   \a \in \NN^{d}_0.
\end{equation*}
Together with the fact that $z_k = z_{d-k+1}$, one can derive a sequence of real orthogonal
polynomials from either $\{T_\a\}$ or $\{U_\a\}$. The set of orthogonal polynomials 
$\{U_\a\}_{|\a| =n}$ of degree $n$ has $\dim \Pi_{n-1}^d$ distinct real
common zeros in $\Delta^*$, 
so that the Gaussian cubature formula exists for $W_{1/2}$ on $\Delta^*$
by Theorem \ref{thm:GaussCuba}.
\index{cubature formula!Gaussian}%
The Gaussian cubature, however, does not exist for $W_{-1/2}$.  

\paragraph{Further results and references}\quad
These orthogonal polynomials were studied systematically 
in \cite{Beer}, which extended earlier work in \cite{Koorn74b} for $d =2$ and partial results 
in \cite{DuLi, EiLi74, EiLi}. The presentation here follows \cite{LhyXu1}, which studied these 
polynomials from the point of view of tiling and discrete Fourier analysis
and, in particular,  studied their common zeros.
The Chebyshev polynomials of the second kind are closely related
to the \emph{Schur functions}.
\index{Schur functions}%
In fact, 
$$
   \TS_{\kb + \vb^\circ}(\tb) = \det \left(z_j^{\l_k + \beta} \right)_{1 \le j,k \le d+1}, \quad   z_j = e^{2\pi i t_j},
$$
where $\l : = (k_1-k_{d+1}, \ldots, k_d - k_{d+1}, 0)$ and $\beta = (d, d-1,\ldots, 1,0)$. In terms of 
symmetric polynomials in $z_1,\ldots, z_d$, they are related to the
$\mathrm{BC}_n$ type orthogonal polynomials; see this volume,
Chapter 8, \cite{BeOp, Vret} and the references therein. 

\subsection{Sobolev orthogonal polynomials}
\index{Sobolev orthogonal polynomials}%
\index{orthogonal polynomials of several variables!Sobolev}%
\index{orthogonal polynomials of several variables!on classical
domains!on unit ball}%
Despite extensive studies of Sobolev orthogonal polynomials
in one variable, there are few
results in several variables until now, and what is known is mostly
on the unit ball $\BB^d$. 
Let $\CH_n^d$ be the space of harmonic polynomials of degree $n$ as in
\S\ref{sec:harmonics}
and let $\CV_n^d(W_\mu)$ denote the space of orthogonal polynomials on $\BB^d$ with 
respect to $W_\mu(x) := (1-\|x\|^2)^{\mu}$, which differs from
\eqref{ball-weight} by a shift of $\half$ in the index. 
\paragraph{First family on $\BB^d$}\quad
Let $\Delta := \partial_1^2 + \ldots + \partial_d^2$ be the Laplace operator. 
Define the inner product on the unit ball $\BB^d$ by
\begin{equation}\label{Sobolev1}
  \la f, g \ra_\Delta :=  \frac1{4 d^2 \operatorname{vol}(\BB^d)}
  \int_{\BB^d} \Delta\Big((1-\|x\|^2) f (x)\Big)\,
           \Delta\Big((1-\|x\|^2) g (x)\Big)\,dx, 
\end{equation}
which is normalized such that $\la 1,1\ra_\Delta =1$. The space 
$\CV_n^d$ of orthogonal polynomials of degree~$n$ for $\la \cdot, \cdot\ra_{\Delta}$
satisfies an orthogonal decomposition 
\begin{equation}\label{VDelta}
\CV_n^d  = \CH_n^d \oplus (1- \|x\|^2\, \CV_{n-2}^d(W_2),  
\end{equation}
from which explicit orthonormal bases can be derived easily. 
\paragraph{Second family on $\BB^d$}\quad
Let $\nabla := (\partial_1, \ldots, \partial_d)$. Define 
the inner product 
\begin{equation}\label{Sobolev2}
  \la f, g\ra_{\nabla} := \frac{\l }{\omega_d} 
      \int_{\BB^d} \nabla f(x) \cdot \nabla g(x)\,dx 
     + \frac{1}{\omega_d} \int_{\SS^{d-1}} f(x) g(x)\,d\sigma_d(x), 
\end{equation}
where $\lambda>0$ such 
that  $\la 1, 1 \ra_{\nabla} =1$. The space of 
orthogonal polynomials of degree $n$ for $\la \cdot, \cdot\ra_{\nabla}$ 
satisfies an orthogonal decomposition 
\begin{equation}\label{V-nabla}
\CV_n^d = \CH_n^d \oplus (1- \|x\|^2)\,\CV_{n-2}^d(W_1).  
\end{equation}
Moreover, the polynomials in $\CV_n^d$ are eigenfunctions of a 
second order differential operator that is exactly the limiting case of 
\eqref{eq:Bdiff} with $\mu = -1/2$. 
\paragraph{Third family on $\BB^d$}\quad
Define the inner product by 
\begin{equation} 
   \la f, g \ra : = \frac{\l}{\omega_d} \int_{\BB^d} \Delta f(x)\,\Delta g(x)\,dx + 
   \frac{1}{\omega_d}  \int_{S^{d-1}} f(x) g(x)\,d\omega, 
\end{equation}
where $\lambda>0$ such  that   $\la 1,1\ra =1$. The space 
$\CV_n^d$ of orthogonal polynomials of degree $n$ satisfies an
orthogonal decomposition
\begin{equation}\label{V-third}
\CV_n^d =  \CH_n^d  \oplus  (1-\|x\|^2)\,\CH_{n-2}^d \oplus
       (1-\|x\|^2)^2\,\CV_{n-4}^d(W_2).  
\end{equation}

\paragraph{Further results and references}\quad
The first family was studied in \cite{Xu06b}.
The motivation for the inner product \eqref{Sobolev1} came from a
Galerkin method 
in the numerical solution of the Poisson equation on the disk.
The second family was studied in \cite{Xu08}.
That reference also considered the inner product where the second
integral in \eqref{Sobolev2} is replaced by $f(0)g(0)$.
The case where  the second
integral in \eqref{Sobolev2} is replaced by an integral over
the ball was studied in \cite{PPX}.
The third family was studied in \cite{PiXu}, where
the connection of orthogonal polynomials with the eigenfunctions of the differential
operator was explored.  Finally, Sobolev orthogonal polynomials with higher
order derivatives in the inner product are studied in \cite{LX14}.
They are used in connection with simultaneous
approximation by polynomials on the unit ball. A first study of Sobolev orthogonal polynomials
on the simplex was conducted in \cite{AX}. A further reference is \cite{LL},
which, however, contains few concrete examples. 

\subsection{Orthogonal polynomials with additional point masses}
\index{orthogonal polynomials of several variables!with additional
point masses}%
Let $\la p, q \ra_\mu : = \int_{\RR^d} p(x)q(x)\,d\mu(x)$ be an inner 
product, for which orthogonal 
polynomials exist. Let $\{\xi_1, \xi_2, \ldots, \xi_N\}$ be a set of distinct points 
in $\mathbb{R}^d$ and let $\Lambda$ be a positive definite matrix of size $N\times N$. 
With the notation
$\mathbf{p}(\xi) = \left \{ p(\xi_1), p(\xi_2), \ldots, p(\xi_N) \right \}$,
considered as a column vector, we define a new inner product 
\begin{equation} 
 \la p, q \ra_\nu := \la p, q \ra_\mu +\mathbf{p}(\xi)^\tran \, \Lambda \,\mathbf{q}(\xi).  
\end{equation}
When $\Lambda  = \mathrm{diag} \{\l_1,\ldots,\l_N\}$, 
the inner product $\la \cdot,\cdot \ra_\nu$ takes the form 
\begin{equation} \label{sum-mass}
  \la p, q \ra_\nu = \la p, q \ra_\mu + \sum_{j=1}^N \l_j p(\xi_j) q(\xi_j).
\end{equation}
The orthogonal polynomials with respect to $\la \cdot,\cdot\ra_{\nu}$
and their kernels 
can be expressed in terms of quantities associated with
$\la \cdot,\cdot\ra_{\mu}$. 

Let $\PP_n$ denote a basis of orthogonal polynomials for $\CV_n^d$ with respect to 
$\la \cdot,\cdot \ra_{\mu}$, as in \eqref{eq:PPdef}, and let $P_n(\mu; \cdot, \cdot)$ 
and $K_n(\mu; \cdot,\cdot)$ denote the reproducing kernel of $\CV_n^d$ and $\Pi_n^d$, 
respectively, with respect to $\la \cdot,\cdot \ra_{\mu}$, as defined in \eqref{eq:reprod} 
and \eqref{eq:ReprodK}. Let  $\mathsf{P}_n(\xi)$ be the matrix that has $\PP_n(\xi_i)$
as columns,
\begin{equation*}
 \mathsf{P}_n(\xi):= \big(\mathbb{P}_n(\xi_1) \mid \mathbb{P}_n(\xi_2) \mid
  \ldots \mid \mathbb{P}_n(\xi_N)   \big) \in \mathcal{M}_{r_n^d \times N},
\end{equation*}
let $\Kb_{n-1}$ be the matrix whose entries are $K_{n-1}(\mu;\xi_i,\xi_j)$, 
\begin{equation*}
\mathbf{K}_{n-1} := \big(K_{n-1}(\mu;\xi_i,\xi_j) \big)_{i,j=1}^N \in \mathcal{M}_{N\times N},
\end{equation*}
and, finally, let $\mathbb{K}_{n-1}(\xi,x)$ be the column vector of functions
\begin{equation*} 
  {\mathbb{K}_{n-1}(\xi,x)} = \left\{K _{n-1}(\mu; \xi_1,x),  K_{n-1}(\mu;\xi_2,x),
\ldots, K_{n-1}(\mu;\xi_N,x)\right \}. 
\end{equation*}
Then the orthogonal polynomials $\QQ_n$ associated with
$\la \cdot,\cdot \ra_{\nu}$ are given by
\begin{equation} \label{ex-expl}
\mathbb{Q}_n(x) = \mathbb{P}_n(x) - {\mathsf{P}_n(\xi)}\,(I_N +
\Lambda\,{\mathbf{K}_{n-1}})^{-1}\,\Lambda \,
{\mathbb{K}_{n-1}(\xi,x)}, \quad n\ge 1,
\end{equation}
and the reproducing kernel of $\Pi_n^d$ associated with 
$\la \cdot,\cdot \ra_{\nu}$ is given by
\index{reproducing kernel}%
\begin{equation}\label{kernel}
K_n(\nu; x,y)  = K_n(\mu; x,y)  -
\mathbb{K}_n^\tran(\xi,x) \,  (I_N + \,\Lambda \,\mathbf{K}_{n})^{-1}  \Lambda \, \mathbb{K}_n(\xi,y).
\end{equation}

These results were developed in \cite{DFPPX}, where the Jacobi weight on the simplex with 
mass points on its vertexes was studied as an example. The results can be
modified to allow 
derivatives at the point masses; for example,  
$$
   \la p, q \ra_\nu := \la p, q \ra_\mu +
   \sum_{j =0}^N \l_j p(\xi_j) q(\xi_j) + \sum_{j=0}^N \l_j' \nabla p(\xi_j) \cdot \nabla q(\xi_j).
$$

\subsection{Orthogonal polynomials for radial weight functions}
\index{weight function!radial}%
\index{orthogonal polynomials of several variables!for radial weight
function}%
Let $w$ be a nonnegative function on the real line with
support setl $[a, b]$, where $0 \le a \le b \le \infty$.
For a \emph{radial weight function} $W(x) := w(\|x\|)$ the orthogonal 
polynomials can be constructed explicitly in polar coordinates.
Indeed, let $p_{2n}^{(2n-4j+d-1)}$
denote the orthonormal polynomials with respect to the weight function $|t|^{2n-4j+d-1} w(t)$
and, for $0 \le j \le n/2$, let $\{Y_{n-2j,\beta}\}$ denote an orthonormal basis for $\CH_{n-2j}^d$ 
of ordinary spherical harmonics. Then the polynomials 
\begin{equation} \label{OPradial_weight}
P_{\beta,j}(x):= p_{2j}^{(2n-4j+d-1)}(\|x\|) Y_{\beta,n-2j}(x)
\end{equation}
form an orthonormal basis of $\CV_n^d$ with $W(x) = w(\|x\|)$. 

The classical examples of radial weight functions are $W_\mu$ in \eqref{ball-weight} on
the unit ball $\BB^d$ and the Hermite weight function $W_H$ in \eqref{Hermite-weight}.
The orthogonal polynomials in \eqref{OPradial_weight} appeared in \cite{Xu05} and they
were used in \cite{Waldron2}.
\index{orthogonal polynomials of several variables!monic
basis|seealso{Appell biorthogonl polynomials}}%
%
\acknow
{\quad\\
I would like to thank Tom Koornwinder for his numerous comments and
corrections.}

\addcontentsline{toc}{section}{\emph{References}}
\bibliography{AskBatVol2}
\bibliographystyle{cambridgeauthordate}

\end{document}